\input amssym.def
\input amssym.tex

%%%%%%%%%%%%%%%%%%%%%%%%%%%%%%%%%%
%%                              %%
%%        IMPAGINAZIONE         %%
%%                              %%
%%%%%%%%%%%%%%%%%%%%%%%%%%%%%%%%%%

\font \title=cmbx10 scaled \magstep2                      
\font \titlepart=cmbx9 scaled\magstep2
\font \ti=cmbx10 scaled\magstep1                        
\font \tit=cmbx9 scaled\magstep1                         

\font \titleit=cmti10 scaled \magstep2                     
\font \titlepartit=cmti9 scaled\magstep2
\font \tiit=cmti10 scaled\magstep1                          
\font \titit=cmti9 scaled\magstep1

\magnification 1000
\hsize 15,8 true cm
\vsize 24 true cm

%%%%%%%%%%%%%%%%%%%%%%%
%                     %
%  DATA DEL PREPRINT  %
%%%%%%%%%%%%%%%%%%%%%%%

\def\today{\ifcase\month\or
  January\or February\or March\or April\or May\or June\or
  July\or August\or September\or October\or November\or December\fi
  \space\number\day, \number\year}

%%%%%%%%%%%%%%%%%%%%%%%%%%%%%%%%%%%%%%%%%%%%%%
%%                                          %%
%%      Definizioni varie                   %%
%%                                          %%
%%%%%%%%%%%%%%%%%%%%%%%%%%%%%%%%%%%%%%%%%%%%%%

\def\b{\Bbb}  \def\c{\cal} 
\def\u{\underline} \def\o{\overline} \def\ovec{\overrightarrow}
\def\p{\prime}\def\l{\left}\def\r{\right}\def\lg{\log}\def\lgp{\log^{+}}
  
\def\[[{\l[\!\l[} \def\]]{\r]\!\r]}
\def\a{\alpha}\def\be{\beta}\def\ga{\gamma}
\def\la{\lambda}\def\La{\Lambda}\def\na{\nabla}
\def\sg{\sigma}\def\Sg{\Sigma}
\def\veps{\varepsilon}

\def\upla#1#2#3{{#1}_{#2},\dots,{#1}_{#3}}

\def\min{\mathop{\rm min}\limits}
\def\max{\mathop{\rm max}\limits}

\def\complessi{{\b C}}\def\reali{{\b R}}\def\raz{{\b Q}}
\def\nat{{\b N}}\def\interi{{\b Z}}

%%%%%%%%%%%%%%%%%%%%%%%%%%%%%%%%%%
%                                %
%        FRECCE                  %
%                                %
%%%%%%%%%%%%%%%%%%%%%%%%%%%%%%%%%%

\def\rightto#1{ \hbox to#1 pt {\rightarrowfill} }
\def\leftto#1{ \hbox to#1 pt {\leftarrowfill} }

\def\hookleftto#1{ \lhook\mkern -8mu\rightto{#1} }

\def\lrarrow{\longrightarrow}
\def\hrarrow{\hookleftto{15}}

%%%%%%%%%%%%%%%%%%%%%%%%%%%%%%%%%%%%%%%%%%%%%%%%%%%%%%%%%
%                                                       %
%  NUMERI E TITOLI DEI CAPITOLO E PARAGRAFI/ETICHETTE   %
%                                                       %
%%%%%%%%%%%%%%%%%%%%%%%%%%%%%%%%%%%%%%%%%%%%%%%%%%%%%%%%%

\newif\ifpps          
\def\void{}
\def\li{\hskip 0pt}

\newcount\chapno\chapno=0                                                   
\def\chapter{\global\advance\chapno by 1}                   
\def\parag{\the\chapno}                 
\def\chaplabel#1{\def\test{#1}
  {\ifx\test\void\errmessage
  {SVEGLIA! NON HAI DATO IL NOME AL CAPITOLO!}
  \else\csname c#1\endcsname
  \expandafter\xdef\csname c#1\endcsname{\S\parag}\fi}
  \ifpps{\par\noindent\line{\null\hskip -65pt{\it c#1}\hfill}}\else\void\fi}

\newcount\subchapno\subchapno=0
\def\subchapter{\global\advance\subchapno by 1}
\def\subparag{\the\subchapno}
\def\subchaplabel#1{\def\test{#1}
  {\ifx\test\void\errmessage {SVEGLIA! NON HAI DATO IL NOME AL PARAGRAFO!}
  \else\csname sc#1\endcsname
  \expandafter\xdef\csname sc#1\endcsname{{\parag.\li\subparag}}\fi}
  \ifpps{\par\noindent\line{\null\hskip -65pt{\it sc#1}\hfill}}\else\void\fi}

\newcount\sectno\sectno=0
\def\section{\global\advance\sectno by1}
\def\form{\the\sectno}
\def\label#1{\def\test{#1}
  {\ifx\test\void\errmessage {SVEGLIA! NON HAI DATO IL NOME AL PARAGRAFO!}
  \else\csname p#1\endcsname
  \expandafter\xdef\csname p#1\endcsname{{\parag.\li\subparag.\li\form}}\fi}
  \ifpps{\par\noindent\line{\null\hskip -65pt{\it p#1}\hfill}}\else\void\fi}

\newcount\equano\equano=0
\def\equation{\global\advance\equano by 1}
\def\numero#1{\equation\def\test{#1}
   {\ifx\test\void\errmessage
   {SVEGLIA! NON HAI DATO IL NOME ALLA FORMULA!}
   \else\csname f#1\endcsname
   \expandafter\xdef\csname f#1\endcsname{\parag.\li\subparag.\li\form.\li\the\equano}\fi}
   \ifpps\leqno{
   \hskip -65pt\hbox to 65pt{\it f#1\hfill}{\rm(\parag.\li\subparag.\form.\the\equano)}}
   \else\leqno{\rm(\parag.\li\subparag.\form.\the\equano)}\fi}

%%%%%%%%%%% CAPITOLI

\def\newchapter#1{\vskip 20pt\chapter
\centerline{\ti\parag. {#1}}{\vskip 5pt}{\global\subchapno=0}}

%%%%%%%%%%% SOTTO-CAPITOLI

\def\newsubchapter#1{\vskip 15pt\subchapter
 \centerline{\tit \parag.\subparag. {#1}}{\vskip 0pt}{\global\sectno=0}}

%%%%%%%%%%% PARAGRAFI

\def\vpar{\vskip 10pt}

\def\newpar#1{\vskip 10pt\section{\global\equano=0}\noindent{\bf #1}{\bf \parag.\subparag.\form.}}
\def\C#1{\vskip 10pt\section{\global\equano=0}\noindent{\bf Corollary  \parag.\subparag.\form.} {\sl #1}}
\def\L#1{\vskip 10pt\section{\global\equano=0}\noindent{\bf Lemma \parag.\subparag.\form.} {\sl #1}}
\def\P#1{\vskip 10pt\section{\global\equano=0}\noindent{\bf Proposition \parag.\subparag.\form.} {\sl #1}}
\def\T#1{\vskip 10pt\section{\global\equano=0}\noindent{\bf Theorem \parag.\subparag.\form.} {\sl #1}}
\def\Def#1{\vskip 10pt\section{\global\equano=0}\noindent{\bf Definition \parag.\subparag.\form.} {\sl #1}}

\def\Oss{\vskip 10pt\section{\global\equano=0}\noindent{\bf Remark \parag.\subparag.\form. }}
\def\Ex{\vskip 10pt\section{\global\equano=0}\noindent{\bf Example \parag.\subparag.\form. }}

\def\Dim{\vskip 5pt\noindent{\bf Proof. }}
\def\tparnoind{\hfill\break}      
\def\tpar{\hfill\break\indent}

\font\tenrsfs=rsfs10
\font\sevenrsfs=rsfs7
\font\fiversfs=rsfs5
\newfam\rsfsfam
\textfont\rsfsfam=\tenrsfs
\scriptfont\rsfsfam=\sevenrsfs
\scriptscriptfont\rsfsfam=\fiversfs
\def\cors#1{{\fam\rsfsfam\relax#1}}

%\overfullrule=0pt

%%%%%%%% DEFINIZIONI DELL'ULTIMA ORA %%%%%%%%%%%%%%%%%%%%%

\def\dq{d_q}
\def\Dq{\Delta_q}
\def\nq{[n]_q}
\def\iq{[i]_q}
\def\ord{{\rm ord}}

\def\QED{\null\ \hfill{$\blacksquare$}}

\def\boxit#1{\vbox{\hrule\hbox{\vrule\kern3pt{\global\sectno=0}
     \vbox{\kern3pt #1\kern3pt}\kern 3pt\vrule}\hrule}}

\def\bigdisplay#1{\hbox{\centerline{$\displaystyle #1$}}}
\def\V#1{\mathaccent 20 #1}

%%%%%%%%%%%%% Numeri di Teoremi importanti

\def\nilpotenza{\S 2}
\def\diff{\S 3}
\def\padico{\S 4}
\def\nilpred{\S 5}
\def\regolarita{\S 6}
\def\groth{\S 7}
\def\prova{\S 8}
\def\esempi{\S 11}

\def\analitico{4.2.6}
\def\transfert{4.3.3}

\def\lemmarith{6.1.2}

\def\grothendieck{7.1.1}

\def\stepb{8.1}
\def\stepc{8.2}
\def\smalldivisor{8.3}
\def\stepd{8.4}

\def\ce#1{\setbox0=\hbox{#1}\ifdim\ht0=1ex\accent'30#1\else
         {\ooalign{\hidewidth\char'30\hidewidth\crcr\unbox0}}\fi}

%\nopagenumbers\ \ \vfill\eject\pageno=1

%%%%%%%%%%%%%%%%%%%%%%%%%
%%%%%%%%%%%%%%%%%%%%%%%%% FOOTLINE/HEADLINE
%%%%%%%%%%%%%%%%%%%%%%%%%

\def\makeheadline{\vbox to 0pt{\vskip -22.5pt
   \line{\vbox to 8.5pt{}\the\headline}\vss}\nointerlineskip}

\headline={\ifnum\pageno=1\hfill\else
\hfil{\rm The $q$-analogue of Grothendieck-Katz's conjecture on $p$-curvatures}\hfil\fi}

\footline{\hfill\rm\folio\hfill}

%%%%%%%%%%%%%%%%%%%%%%%%%
%%%%%%%%%%%%%%%%%%%%%%%%% FOOTNOTE
%%%%%%%%%%%%%%%%%%%%%%%%%

\footnote {}{2000 Mathematics Subject Classification: 12H99 (33D15, 39A13).}

%%%%%%%%%%%%%%%%%%%%%%%%%
%%%%%%%%%%%%%%%%%%%%%%%%% INTESTAZIONE
%%%%%%%%%%%%%%%%%%%%%%%%%
\null\vskip 10pt
\centerline{\title
Arithmetic theory of {\titleit q}-difference equations}
\vskip 15pt
\centerline{\ti
The {\tiit q}-analogue of Grothendieck-Katz's conjecture on {\tiit p}-curvatures}
\vskip 20pt
\centerline{\tit by}
\vskip 10pt
\centerline{\tit Lucia Di Vizio}

\vpar
\centerline{\sevenrm Laboratoire Emile Picard, Universit\'e Paul Sabatier}
\centerline{\sevenrm 118, route de Narbonne, 31062 Toulouse, France}
\centerline{\sevenrm e-mail: divizio@picard.ups-tlse.fr}

%%%%%%%%%%%%%%%%%%%%%%%%%
%%%%%%%%%%%%%%%%%%%%%%%%% INDICE
%%%%%%%%%%%%%%%%%%%%%%%%%

\vskip 15pt\noindent
{\ti Table of contents}
\chapno=0

\vpar\noindent Introduction

\vpar\noindent Part I. Generalities on $q$-difference modules

\par\noindent\chapter\S\parag. $q$-difference modules 
\par\subchapter{\parag.\subparag.} Summary of $q$-difference algebra
\par\subchapter{\parag.\subparag.} The $q$-analogue of the Wronskian lemma
\par\subchapter{\parag.\subparag.} The $q$-analogue of the cyclic vector lemma
\par\subchapter{\parag.\subparag.} Formal classification of $q$-difference modules

\subchapno=0
\par\noindent\chapter\S\parag. Unipotent $q$-difference modules
\par\subchapter{\parag.\subparag.} Trivial $q$-difference modules
\par\subchapter{\parag.\subparag.} Extensions of trivial $q$-difference modules

\vpar\noindent Part II. $p$-adic methods

\subchapno=0
\par\noindent\chapter\S\parag. Considerations on the differential case

\subchapno=0
\par\noindent\chapter\S\parag. Introduction to $p$-adic $q$-difference modules
\par\subchapter{\parag.\subparag.} $p$-adic estimates of $q$-binomials
\par\subchapter{\parag.\subparag.} The Gauss norm  and the invariant $\chi_v({\c M})$
\par\subchapter{\parag.\subparag.} $q$-analogue of the Dwork-Frobenius theorem

\subchapno=0
\par\noindent\chapter\S\parag. $p$-adic criteria for unipotent reduction
\par\subchapter{\parag.\subparag.} $q$-difference modules having unipotent reduction modulo $\varpi_v$
\par\subchapter{\parag.\subparag.} $q$-difference modules having unipotent reduction modulo $1-q^{\kappa_v}$

\vpar\noindent Part III. A $q$-analogue of Grothendieck's conjecture on $p$-curvatures

\subchapno=0
\par\noindent\chapter\S\parag. Arithmetic $q$-difference modules and regularity
\par\subchapter{\parag.\subparag.} On cyclic subgroups of $\o\raz^\times$ and their reduction modulo almost every prime
\par\subchapter{\parag.\subparag.} Unipotent reduction and regularity

\subchapno=0
\par\noindent\chapter\S\parag. Statement of the $q$-analogue of Grothendieck's conjecture on $p$-curvatures
\par\subchapter{\parag.\subparag.} Statement of the theorem
\par\subchapter{\parag.\subparag.} Idea of the proof 

\subchapno=0
\par\noindent\chapter\S\parag. Proof of (\grothendieck)
\par\subchapter{\parag.\subparag.} Finiteness of size of $\c M$
\par\subchapter{\parag.\subparag.} Finiteness of size of a fundamental matrix of solutions
\par\subchapter{\parag.\subparag.} How to deal with the problem of archimedean small divisors
\par\subchapter{\parag.\subparag.} Conclusion of the proof: a criterion for rationality
\par\subchapter{\parag.\subparag.} A corollary

\vpar\noindent Part IV. A $q$-analogue of Katz's conjectural
description of the generic Galois group

\subchapno=0
\par\noindent\chapter\S\parag. Definition of the generic $q$-difference Galois group
\par\subchapter{\parag.\subparag.} Some algebraic constructions
\par\subchapter{\parag.\subparag.} Definition of the Galois group of a $q$-difference module
\par\subchapter{\parag.\subparag.} Definition of the generic Galois group of a $q$-difference module

\subchapno=0
\par\noindent\chapter\S\parag. An arithmetic description of the generic Galois group
\par\subchapter{\parag.\subparag.} Algebraic groups ``containing $\Phi_q^{\kappa_v}$ for almost all $v$''
\par\subchapter{\parag.\subparag.} Statement of the main theorem
\par\subchapter{\parag.\subparag.} Proof of the main theorem

\subchapno=0
\par\noindent\chapter\S\parag. Examples of calculation of generic Galois groups

\subchapno=0
\vpar\noindent Appendix. A $q$-analogue of Schwarz's list
\par\subchapter{A.\subparag.} Logarithmic singularities
\par\subchapter{A.\subparag.} The case of algebraic solutions

\vpar\noindent References

\chapno=0
\subchapno=0

%%%%%%%%%%%%%%%%%%%%%%%%%%%%%%%%%%%%%%%%%%%%%%%%%%%%%%%%%%%%%%%%%%%%%%%%%%%%
%%%%%%%%%%%%%%%%%%%%%%%%%%%%%%%%%%   0   %%%%%%%%%%%%%%%%%%%%%%%%%%%%%%%%%%%
%%%%%%%%%%%%%%%%%%%%%%%%%%%%%%%%%%%%%%%%%%%%%%%%%%%%%%%%%%%%%%%%%%%%%%%%%%%%

\vskip 30pt\noindent
\centerline{\titlepart Introduction}
\vpar
%%%%%%%%%%%%%%%%%%%%%%%%%%%%%%%%%%%%%%%%%%%%%%%%%%%%%%%%%%%%%%%%%%%%%%%%%%%%

\vpar\vpar
In 1941 Birkhoff and Guenther 
wrote: ``{\it Up to the present time the theory of linear $q$-difference equations
has lagged noticeably behind the sister theories of linear difference
and differential equations. In the opinion of the authors the use of
the canonical system [...] is
destinated to carry the theory of $q$-difference equations to a
comparable degree of completeness}''.   
In the same paper they announced a program which they did not develop further, and 
$q$-difference equations theory remains less advanced today than difference and differential 
equation theories. 
\par
In recent years, mathematicians have reconsidered $q$-difference equations 
for their links with other branches of mathematics such as 
quantum algebras and $q$-combinatorics, and 
Birkhoff and Guenther's program has been continued. 
Now there are theories of divergence (J-P. B\'ezivin, J-P. Ramis) and of 
$q$-summation (Ch. Zhang).  There is also a good analogue of the concept 
of monodromy and thus a description of the $q$-difference Galois 
group in the regular case (P.I. Etingof, M. van der Put and M. Singer, J. Sauloy). 
\par
B\'ezivin and Ramis's results on divergent series have applications to rationality 
criteria for solutions of systems of $q$-difference equations ({\it cf.}\ [BB]) 
and for systems of $q$-difference and differential equations ({\it cf.}\ [Ra]), 
which provides an answer to the old problem of finding criteria to establish whether a 
formal power series is the Taylor expansion of an algebraic or a rational function.
\par
The question was first raised by Schwarz, who established an exhaustive
list of hypergeometric differential equations having a full set of algebraic solutions. 
Grothendieck's conjecture on $p$-curvatures tries to give a complete answer to this problem. 
More precisely, when we consider a differential equation 
$$
{\c L}y=a_\mu(x){d^\mu y\over dx^\mu}+
a_{\mu-1}(x){d^{\mu-1} y\over dx^{\mu-1}}+
\dots+a_0(x)y=0\ ,
$$ 
with coefficients in the field $\raz(x)$, we can reduce the equation
${\c L}y=0$ modulo $p$ for almost all primes $p\in\interi$.  
Then Grothendieck's conjecture predicts:

\vpar\noindent
{\bf Grothendieck's conjecture on $p$-curvatures.}
{\sl The equation ${\c L}y=0$ has a full set of algebraic solutions 
if and only if for almost all primes $p\in\interi$ the reduction modulo
$p$ of ${\c L}y=0$ has a full set of solutions in ${\b F}_p(x)$.}

\vpar\noindent
In spite of numerous papers dedicated to this conjecture in which
some particular cases are proved (we recall [Ho], [CC], [K2], [K3], [A2], [Bo]), 
the conjecture remains open.

\vpar
In this paper we give a proof of an analogous statement for 
$q$-difference equations.  Following [K3], this allows for an arithmetic description 
of the generic Galois group of a $q$-difference equation. 
In fact, in [K3] N. Katz proposes a conjectural arithmetic 
description of the generic Galois group of a differential equation which 
is equivalent to Grothendieck's conjecture:

\vpar\noindent
{\bf Katz's conjectural description of the generic Galois group.}
{\sl The Lie algebra of $Gal(M)$ is the smallest algebraic Lie sub-algebra of 
${\rm End}_{\raz(x)}(M)$ whose reduction modulo $p$ contains the $p$-curvature $\psi_p$
for almost all $p$.}
 
\vpar
Let us briefly explain his statement. 
Let ${\c M}=(M,\na)$ be a $\raz(x)$-vector space with a $\raz(x)/\raz$-connection. 
We define the generic Galois group $Gal({\c M})$ of ${\c M}$ to be the algebraic subgroup of
$GL(M)$ stabilizing all the sub-quotients of the mixed tensor spaces 
$\oplus_{i,j}(M^{\otimes^i}\otimes_{\raz(x)}(M^\ast)^{\otimes^j})$. 
We can consider a lattice $\widetilde M$ of $M$ over a finite type algebra over
$\interi$, stable under the connection, and we can 
reduce $\widetilde M$ modulo $p$, for almost all primes $p$.
The operator $\psi_p=\na\big({d\over dx}\big)^p$ acting over
$\widetilde M\otimes_\interi {\b F}_p$, is called the {\it $p$-curvature}. 
Moreover it makes sense to consider the
reduction modulo $p$ for almost all $p$ of $Gal(M)$ and its Lie algebra. 

\vpar
\centerline{$\ast\ast\ast$}
\vpar

The present paper contains proofs of the analogues of the 
conjectures above for $q$-difference equations. 
More precisely, let $q$ be a nonzero rational number.
We consider the $q$-difference equation 
$$
{\c L}y=a_\mu(x)y(q^\mu x)+a_{\mu-1}(x)y(q^{\mu-1}x)+
\dots+a_0(x)y(x)=0\ ,\ a_0(x)\neq 0\neq a_\mu(x)\ ,
$$
with $a_j(x)\in\raz(x)$, for all $j=0,\dots,\mu$. 
For almost all rational primes $p$ the image $\o q$ of 
$q$ in ${\b F}_p$ is nonzero and generates a cyclic subgroup of 
${\b F}_p^\times$ of order $\kappa_p$, and 
there exists a positive integer $\ell_p$ such that
$1-q^{\kappa_p}=p^{\ell_p}{h\over g}$, with $h,g\in\interi$ prime with respect to
$p$. We denote by ${\c L}_py=0$ the reduction of ${\c L}y=0$ 
modulo $p^{\ell_p}$.
Let us consider a $\interi$-algebra 
${\c A}=\interi\l[x,{1\over P(q^ix)}, i\geq 0\r]$, with
$P(x)\in\interi[x]\smallsetminus\{0\}$, such that 
$a_j(x)\in \interi\l[x,{1\over P(q^ix)}, i\geq 0\r]$, for all $j=0,\dots,\mu$. 
Our result is ({\it cf.}\ (\grothendieck) below):

\vpar\noindent
{\bf Theorem 1.}
{\sl The $q$-difference equation ${\c L}y=0$ has a full set of
solutions in $\raz(x)$ if and only if for almost all rational primes $p$ the 
set of equations ${\c L}_py=0$ has a full set of solutions in  
${\c A}\otimes_\interi\interi/p^{\ell_p}\interi$.}
\def\grothendieckint{1}

\vpar
Theorem 1 is a partial answer to a B\'ezivin's conjecture ({\it cf.}\ [Be1]).
\vpar
Let $M$ be a finite dimensional $\raz(x)$-vector space equipped with a 
$q$-difference operator $\Phi_q:M\lrarrow M$, {\it i.e.}, with a
$\raz$-linear invertible morphism such that
$\Phi_q(fm)=f(qx)\Phi_q(m)$ for all $f(x)\in\raz(x)$ and all $m\in M$.
As in the differential case, it is equivalent to consider a $q$-difference equation 
or a couple $(M,\Phi_q)$. 
\par
One can attach to ${\c M}=(M,\Phi_q)$ 
an algebraic closed subgroup $Gal({\c M})$ of $GL(M)$,  
that we call the $q$-difference generic Galois group. 
It is the stabilizer of all $q$-difference sub-modules 
of all finite sums of the form
$\oplus_{i,j}(M^{\otimes^i}\otimes_{\raz(x)}(M^\ast)^{\otimes^j})$, 
equipped with the operator induced by $\Phi_q$.
We consider the reduction modulo $p^{\ell_p}$ of $M$ 
for almost all $p$, by reducing a lattice $\widetilde M$ of $M$, defined over a
$\interi$-algebra and stable by $\Phi_q$. 
The algebraic group $Gal({\c M})$ can also be reduced modulo $p^{\ell_p}$ for
almost all $p$.
Then our description of $Gal({\c M})$ is the following:

\vpar\noindent
{\bf Theorem 2.}
{\sl The algebraic group $Gal({\c M})$ is the smallest algebraic subgroup of 
$GL(M)$ whose reduction modulo $p^{\ell_p}$ contains the reduction of
$\Phi_q^{\kappa_p}$ modulo $p^{\ell_p}$ for almost all $p$.}

\vpar\noindent
Taking into account the fact that $\Phi_q$ is a semi-linear
endomorphism (which is easier to handle than the higher derivations
occurring in the differential case), sometimes it happens that one can
calculate all $\Phi_q^n$ at once and therefore determine the 
generic Galois group.

\vpar
\centerline{$\ast\ast\ast$}
\vpar
The techniques employed in the proof of theorem \grothendieckint\ are borrowed from the theory of
$G$-functions. There are essentially two properties of
arithmetic $q$-difference equations which allow us to obtain stronger
results than in the differential case:
\par\noindent
1) A formal power series with a nonzero radius of convergence,  
which is a solution of a $q$-difference equation, 
has infinite radius of meromorphy whenever $|q|>1$: we say that solutions of 
$q$-difference equations have good meromorphic uniformization.
If the algebraic number $q$ is not a root of unity, 
one can always find a place, archimedean or not,  
such that the associated norm of $q$ is greater than 1. This is the 
key-point of the proof: 
if we had good meromorphic uniformization of solutions of arithmetic
differential equations, Grothendieck's conjecture would become a
corollary of $G$-function theory. 
\par\noindent
2) An arithmetic differential equation whose reduction modulo $p$ can be written as a product 
of trivial factors for almost all $p$ 
is regular singular and has rational exponents ({\it cf.}\ [K1, 13.0]). 
A $q$-difference equation whose reduction modulo $p$ can be written as a product 
of trivial factors for almost every $p$ is not only
regular singular, but its ``exponents'' are in $q^\interi$: 
this means that the equation 
has a complete set of solutions in $K((x))$.

\vpar
In one instance the techniques used in $G$-function theory give a weaker
result in the $q$-difference case: the $q$-analogue of the Katz estimates for the $p$-adic
generic radius of convergence is very unsatisfactory ({\it cf.}\ \nilpred\ below). This is at
the origin of many complications in the text 
({\it cf.}\ (\stepb)): actually the naive $q$-analogue of 
the notion of nilpotent reduction does not allow us to conclude the
proof of theorem 1.
A deeper analysis of the definition of $p$-curvatures for
arithmetic differential equation shows that we can define two 
$q$-analogues of the notion of trivial reduction ({\it cf.}\ \diff). 
Both of them are natural and useful. The first one permits us only to obtain the
triviality over $K((x))$, the second one leads to the triviality
over $K(x)$. 

\vpar
Finally, we want to stress the fact that we have very poor
information on the sequence of integers $(\kappa_p)_p$ and no control at
all over $(\ell_p)_p$. 
We are just able to prove ({\it cf.}\ (\lemmarith)) that the
sequence $(\kappa_p)_p$ completely determines the set $\{q,q^{-1}\}$.  
The difficulties linked to these numbers and
their distribution are the arithmetical counterpart of the 
classical (archimedean) problem of small divisors. 
This becomes clearer if we translate the definition of $\kappa_p$ and
$\ell_p$ as follows:
$$
\kappa_p=\min\{m\in\interi:~m>0,~|1-q^m|_p<1\}
$$
and $p^{-\ell_p}=|1-q^{\kappa_p}|_p$, where $|~|_p$ 
is the $p$-adic norm over $\raz$ such that
$|p|_p=p^{-1}$. 
\par
Concerning the archimedean problem of small divisors, 
it would seem natural 
to assume that for all embeddings $\o\raz\hrarrow\complessi$, the
image of $q$ in $\complessi$ does not have complex norm $1$.  
Actually, this assumption is not needed since $q$ is an algebraic
number. 
In fact in (\smalldivisor) 
we need to show that a formal power series $y(x)\in \complessi\[[x\]]$ 
which is a solution of 
a regular singular $q$-difference equation with coefficients in
$\complessi(x)$ is convergent. 
In [Be2], the author gives some technical sufficient conditions on the
estimate of $|1-q^n|_\complessi$
to assure the convergence of the power series $y(x)$. 
It is a consequence of Baker's
theorem on linear forms in logarithms 
that these conditions are always satisfied when $q$ is an algebraic number: 
the idea is already present in [Be1].
It is possible that the techniques of (\smalldivisor) can be applied to more
general problems of small divisors. 

\vpar
\centerline{$\ast~\ast~\ast$}
\vpar

This paper is organized as follows. 
\par
In the first part we introduce some basic
properties of $q$-difference modules, in particular a $q$-analogue of
the cyclic vector lemma. We then recall some results on the
formal classification of $q$-difference modules.
In \nilpotenza\ we prove a characterization of
trivial $q$-difference modules and of $q$-difference modules which are 
extensions of trivial modules when $q$ is a root of unity and $K$ 
is a commutative ring. This degree of generality is motivated by theorem
\grothendieckint, where we consider a $q$-difference equation over
a $\interi/p^{\ell_p}\interi$-algebra.
\par
Part II is devoted to the $p$-adic situation. 
Section \diff\ contains some considerations on arithmetic
differential modules, with the purpose of motivating the choice of 
considering two different $q$-analogues of the notion of nilpotent
reduction. 
In \padico\ we introduce $p$-adic $q$-difference
modules and we establish their primary properties. In particular we prove a
$q$-analogue of the Dwork-Frobenius-Young theorem. In \nilpred\ we 
introduce the two notions of nilpotent reduction and  
revisit and translate some classical estimates for differential modules having 
nilpotent reduction in the $q$-difference setting ({\it cf.}\ [DGS, page 96]). 
The results
of this section are crucial for the proof of the main theorem
(\grothendieck), together with the results of \regolarita. 
\par
In Part III we consider the
arithmetic situation. In \regolarita\ we prove a $q$-analogue of [K1, 13.0]: 
as we have already pointed out we obtain a stronger result than
in the differential setting.
Section \groth\ contains the statement of theorem 1 and \prova\ 
its proof. 
\par
In Part IV we introduce the generic $q$-difference Galois group and prove theorem 
2. In \esempi\ we show how theorem 2 
can sometimes be an effective instrument for calculating Galois groups. 
\par
Finally, the appendix contains an analogue of Schwarz's list
for basic hypergeometric series.

\vpar\noindent
{\bf Acknowledgements. }
I am grateful to Y. Andr{\'e} for suggesting that I work on an arithmetic 
theory of $q$-difference equations and
for assisting me during all stages of preparing this work.  
I would like to thank M. van der Put for providing me with a manuscript containing
the proof of (6.2.2, 1), with some notes on the reduction modulo $p$ of
$q$-difference equations of rank 1 which he gave me in April 1999, and also
for many conversations we had on the same occasion.  
This paper is a part of my Ph.D. thesis, that was referred by 
J-P. B\'ezivin et J-P. Ramis: 
I would like to thank them for their remarks.  
Finally, I am indebted
to P. Colmez for his suggestion concerning (\lemmarith), to  J. Sauloy 
for our discussion on the appendix, to the referee of this paper for 
his many useful remarks, in particular for 
suggesting a simpler proof of (\stepd.1, Step 6) and to S. Dumas 
for correcting many English language inaccuracies in the text. 
Of course, I'm alone responsible for all 
the deficiencies that remain.

%%%%%%%%%%%%%%%%%%%%%%%%%%%%%%%%%%%%%%%%%%%%%%%%%%%%%%%%%%%%%%%%%%%%%%%%%
%%%%%%%%%%%%%%%%%%%%%%%% PART I %%%%%%%%%%%%%%%%%%%%%%%%%%%%%%%%%%%%%%%%%
%%%%%%%%%%%%%%%%%%%%%%%%%%%%%%%%%%%%%%%%%%%%%%%%%%%%%%%%%%%%%%%%%%%%%%%%%
%%%%%%%%%%%%%%%%%%%%%%%%%%%%%%%%%%%%%%%%%%%%%%%%%%%%%%%%%%%%%%%%%%%%%%%%%

\vskip 20pt 
\centerline{\titlepart Part I. Generalities on {\titlepartit q}-difference modules}

%%%%%%%%%%%%%%%%%%%%%%%%%%%%%%%%%%%%%%%%%%%%%%%%%%%%%%%%%%%%%%%%%%%%%%%%%
%%%%%%%%%%%%%%%%%%%%%%%%%%%%%%%%%%%%%%%%%%%%%%%%%%%%%%%%%%%%%%%%%%%%%%%%%
%%%%%%%%%%%%%%%%%%%%%%%%%%%%%%%%%%%%%%%%%%%%%%%%%%%%%%%%%%%%%%%%%%%%%%%%%
%%%%%%%%%%%%%%%%%%%%%%%%%%%%%%%%%%%%%%%%%%%%%%%%%%%%%%%%%%%%%%%%%%%%%%%%%

%%%%%%%%%%%%%%%%%%%%%%%%%%%%%%%%%%%%%%%%%%%%%%%%%%%%%%%%%%%%%%%%%%%%%%%%%
%%%%%%%%%%%%%%%%%%%%%%%%%%%%  1  %%%%%%%%%%%%%%%%%%%%%%%%%%%%%%%%%%%%%%%%
%%%%%%%%%%%%%%%%%%%%%%%%%%%%%%%%%%%%%%%%%%%%%%%%%%%%%%%%%%%%%%%%%%%%%%%%%

\newchapter{{\tiit q}-difference modules}
\chaplabel{definizioni}

%%%%%%%%%%%%%%%%%%%%%%%%%%%%%%%%%%%%%%%%%%%%%%%%%%%%%%%%%%%%%%%%%%%%%%%%%

%%%%%%%%%%%%%%%%%%%%%%%%%%%%%%%%%%%%%%%%%%%%%%%%%%%%%%%%%%%%%%%%%%%%%%%%%
%%%%%%%
\newsubchapter{Summary of {\titit q}-difference algebra}
%%%%%%%
%%%%%%%%%%%%%%%%%%%%%%%%%%%%%%%%%%%%%%%%%%%%%%%%%%%%%%%%%%%%%%%%%%%%%%%%%

\vpar
{\it Let $R$ be a commutative ring.} 

\newpar{}
{\bf $q$-binomials.} 
For any $a,q\in R$ and any integer $n\geq 1$, we shall use the following standard notation:
$$\eqalign{
&[0]_q=0\ ,\ \nq=1+q+\dots+q^{n-1}\ ,\cr
&[0]_q!=1\ ,[n]_q!=1_q\cdots\nq\ ,\cr
&(x-a)_0=1\ ,\ (x-a)_n=(x-a)(x-qa)\cdots(x-q^{n-1}a)\ ,\cr
&(a,q)_0=1\ ,\ (a;q)_n=(1-a)_n\ .}
$$
If $q\neq 1$, we have $\nq={1-q^n\over 1-q}$.
\par
The {\it $q$-binomial coefficients} ${n\choose i}_q$ are the
elements  of $R$ defined by the polynomial identity 
$$
(1-x)_n=\sum_{j=0}^n(-1)^j{n\choose j}_q q^{j(j-1)/2}x^j\ .
\numero{binomio}$$ 
It was already known to Gauss that these are polynomials in $q$ 
which have the following properties:
$$\eqalign{
&{n\choose 0}_q={n\choose n}_q=1\cr
&{n\choose i}_q={[n]_q!\over [n-i]_q![i]_q!}={\nq [n-1]_q\cdots[n-i+1]_q\over [i]_q!}\ ,\cr
&{n\choose i}_q
={n-1\choose i-1}_q+{n-1\choose i}_q q^i
={n-1\choose i-1}_q q^{n-i}+{n-1\choose i}_q\ ,
\ \hbox{for $n\geq i\geq 1$.}}
\numero{tartaglia}
$$ 

\newpar{}
{\bf $q$-dilatation.}
We fix a {\it unit} $q$ in $R$.
We shall consider several rings of functions 
of one variable $x$ and uniformly denote by $\varphi_q$ 
the automorphism ``of dilatation" induced by 
$x\mapsto qx$.  We shall denote this automorphism 
either by $f(x)\longmapsto f(qx)$ or by 
$f\longmapsto \varphi_q(f)$. 
\vpar
We shall informally refer to an $R$-algebra $\c F$ of functions endowed with the 
operator $\varphi_q$ as a {\it $q$-difference algebra over $R$}.  A morphism 
${\c F}\lrarrow{\c F}^\p$ of $q$-difference algebras is a morphism of $R$-algebras commuting with 
the action of $\varphi_q$. 
Moreover, we shall say that a $q$-difference algebra $\c F$ over $R$ is {\it essentially of finite type} 
if there exist $P_1,\dots,P_n\in{\c F}$ such that ${\c F}=R[P_1(q^ix),\dots,P_n(q^ix);i\geq 0]$.
\label{dilatationdefinizione}

\vpar\noindent
{\bf Examples.}
Typical examples of $q$-difference algebras are:
\par\noindent
{\it (i)} $R((x))$, with the obvious action of $\varphi_q$.
\par\noindent 
{\it (ii)} When $R$ is a field, the subfield $R(x)$ of $R((x))$ is a $q$-difference algebra over $R$.
\par\noindent 
{\it (iii)} The $R$-algebra 
$$
R\[[x-a\]]_q=\l\{
\sum_{n=0}^\infty a_n(x-a)_n:a_n\in R\r\}\ , 
\hbox{\ for $a\in R$, $a\neq 0$,}
$$
with $\varphi_q(x-a)_n=q^n(x-a)_n+q^{n-1}(q^n-1)a(x-a)_{n-1}$. 
\par\noindent
{\it (iv)} Let us consider the $q$-difference algebra $R[x]$ and 
$P_1(x),\dots,P_n(x)\in R[x]$. Then the $R$-algebra 
$$
R\l[x,{1\over P_1(q^ix)},\dots,{1\over P_n(q^ix)},i\geq 0\r]
$$
is a $q$-difference algebra essentially of finite type over $R$. 

\Def
{The ring $C=\{f\in{\c F}:\varphi(f)=f\}$ is the subring of constants ${\c F}$.}

\Ex
Let ${\c F}=R((x))$. 
If $q$ is not a root of unity, 
then $\varphi_q(f)(x)=0$ if and only if $f\in R$.
If $q$ is a primitive root of unity of order $\kappa$, 
then $\varphi_q(f)(x)=f$ if and only if $f\in R((x^\kappa))$.

\Def
{A {\it $q$-difference module ${\c M}=(M,\Phi_q)$ over  a
$q$-difference algebra $\c F$} is a free $\c F$-module $M$ 
of finite rank together with an $R$-linear automorphism:
$$
\Phi_q:M\lrarrow M
$$
satisfying the rule
$$
\Phi_q(f(x)m)=f(qx)\Phi_q(m)\ ,\ 
\hbox{for every $f(x)\in{\c F}$ and every $m\in M$.} 
$$}

\vpar\noindent
{\bf Remark.}
The operator $\Phi_q$ is nothing but a $\varphi_q$-semilinear automorphism 
of the $\c F$-module $M$.

\Def
{A morphism $\psi:(M,\Phi_q)\lrarrow(M^\p,\Phi_q^\p)$ is an $R$-linear morphism $M\lrarrow M^\p$
which commutes with the semilinear automorphisms $\Phi_q$ and $\Phi_q^\p$.}

\vpar
Let us consider a morphism ${\c F}\lrarrow{\c F}^\p$ of $q$-difference algebras and  
a $q$-difference module ${\c M}=(M,\Phi_q)$ over $\c F$. 

\Def
{The $q$-difference module ${\c M}_{{\c F}^\p}$ obtained from ${\c M}$ by extension of
coefficients from $\c F$ to ${\c F}^\p$ is the ${\c F}^\p$-module $M\otimes_{\c F}{\c F}^\p$ equipped with the operator 
$\Phi_q\otimes\varphi_q$.}

\newpar{}
{\bf $q$-derivations.} 
Let $q\neq 1$. 
Until the end of this subsection, 
we assume that $\c F$ is stable with respect to
the operator 
$$\matrix{
\dq :& {\c F}&\lrarrow &{\c F}\cr\cr
&f(x)&\longmapsto &
\displaystyle {\varphi_q-id\over (q-1)x}f(x)&
\displaystyle ={f(qx)-f(x)\over (q-1)x}\cr}\ .
$$
\label{derivazionibendefinite}

\vpar\noindent
{\bf Remark.}
The operator $\dq$ satisfies the twisted Leibniz rule: 
$$
\dq(fg)(x)=\dq(f)(x)g(x)+f(qx)\dq(g)(x)\ .
$$ 
More generally, for any 
positive integer $n$, we have 
$$
\dq^n(fg)(x)=\sum_{j=0}^n {n\choose j}_q\dq^{n-j}(f)(q^jx)\dq^j(g)(x)\ .
\numero{leibniz}$$

\vbox{
\Ex
\smallskip\noindent
1) Let us consider the $q$-difference algebra 
${\c F}=R((x))$. For any positive integer $n$,  
$\dq x^n=\nq x^{n-1}$. More generally}
$$
{\dq^s\over [s]_q!}x^n=\cases{
0&if $n<s$\cr\cr
\displaystyle{n\choose s}_q x^{n-s}& otherwise}\ .
$$
\smallskip\noindent
2)Let ${\c F}=R\[[x-a\]]_q$; then $\dq(x-a)_n=\nq(x-a)_{n-1}$. 
\label{primoesempio}

\vpar 
We have the following relations between $\dq$ and $\varphi_q$:

\L
{We set $\dq^0=\varphi_q^0=1$. For any integer $n\geq 1$ we obtain:
$$
\varphi_q^n=\sum_{i=0}^n{n\choose i}_q(q-1)^i q^{i(i-1)/2}x^i\dq^i
$$
and 
$$
\dq^n={\l(\varphi_q-1\r)_n\over(q-1)^nq^{n(n-1)/2}x^n}
={(-1)^n\over(q-1)^n x^n}
\sum_{j=0}^n (-1)^j{n\choose j}_{q^{-1}}q^{-{j(j-1)\over 2}}\varphi_q^j\ .
$$}
\label{fidelta}

\Dim
We remark that 
$x\dq\circ x^i\dq^i=q^ix^{i+1}\dq^{i+1}+[i]_qx^i\dq^i$,
for all $i\geq 1$.
For $n=2$ one has:
$$
\varphi_q^2=(q-1)^2qx^2\dq^2+[2]_q(q-1)x\dq+1\ .
$$
It follows by induction that 
$$
\bigdisplay{\eqalign{\varphi_q^{n+1}
&=\l((q-1)x\dq+1\r)\varphi_q^n\cr
&=\sum_{i=0}^n{n\choose i}_q(q-1)^i q^{i(i-1)/2}
  \l((q-1)q^ix^{i+1}\dq^{i+1}+(q-1)[i]_qx^i\dq^i+x^i\dq^i\r)\cr
&=(q-1)^{n+1}q^{n(n+1)/2}x^{n+1}\dq^{n+1}+
  \sum_{i=1}^n\l({n\choose i}_q q^i+{n\choose i-1}_q\r)(q-1)^i
  q^{i(i-1)/2}x^i\dq^i+1\cr 
&=\sum_{i=0}^{n+1}{n+1\choose i}_q(q-1)^iq^{i(i-1)/2}x^i\dq^i\ .}}
$$
The second formula in (\pfidelta) holds for $n=1$ by definition of $\dq$. 
By induction we obtain 
$$\eqalign{\dq^{n+1}
&={\varphi_q-1\over (q-1)x}\circ{\l(\varphi_q-1\r)_n\over(q-1)^nq^{n(n-1)/2}x^n}\cr
&={\l(\varphi_q-q^n\r)\l(\varphi_q-1\r)_n\over(q-1)^{n+1}q^{n(n+1)/2}x^{n+1}}\cr
&={\l(\varphi_q-1\r)_{n+1}\over(q-1)^{n+1}q^{n(n+1)/2}x^{n+1}}\cr
&={(-1)^n\over(q-1)^{n+1}x^{n+1}}
  \l(1-\varphi_q\r)\l(1-q^{-1}\varphi_q\r)\cdots\l(1-q^{-n}\varphi_q\r)\ .}
$$
We conclude by using (\fbinomio).
\QED

\vpar\noindent
{\bf Remark.}
Let $M$ be a free $\c F$-module of finite rank. Let $\Dq:M\lrarrow M$ be 
an $R$-linear endomorphism satisfying the twisted Leibniz rule:
$$
\Dq(f(x)m)=f(qx)\Dq(m)+\dq(f)(x)m\ , 
\hbox{\ for every $f(x)\in {\c F}$ and every $m\in M$.}
\numero{rel}$$
Then $\Phi_q=(q-1)x\Dq+1$ is $\varphi_q$-semilinear. 
Therefore, if it is invertible, it defines a $q$-difference module. 
Conversely, if $(q-1)x$ is a unit in $\c F$, any $\Phi_q$ gives rise to a twisted 
derivation $\Dq$ as before. 
We remark that $\Dq$ satisfies the generalized Leibniz formula:
$$
\Dq^n(f(x)m)=\sum_{i=0}^n {n\choose i}_q\dq^{n-i}(f)(q^ix)\Dq^i(m)\ ,
\hbox{\ for all $f\in {\c F}$ and $m\in M$.}
\numero{grothleibniz}
$$

\L
{The analogue of the formulas in (\pfidelta) holds:
$$
\Phi_q^n=\sum_{i=0}^n{n\choose i}_q(q-1)^i q^{i(i-1)/2}x^i\Dq^i
\numero{fiD}$$
and
$$\Delta_q^n={\l(\Phi_q-1\r)_n\over(q-1)^nq^{n(n-1)/2}x^n}
={(-1)^n\over(q-1)^n x^n}
\sum_{j=0}^n (-1)^j{n\choose j}_{q^{-1}}q^{-{j(j-1)\over 2}}\Phi_q^j\ .
\numero{Dfi}
$$}

\Dim
The proof is similar to the proof of (\pfidelta).
\QED

%%%%%%%%%%%%%%%%%%%%%%%%%%%%%%%%%%%%%%%%%%%%%%%%%%%%%%%%%%%%%%%%%%%%%%%%%
%%%%%%%
\newsubchapter{The {\titit q}-analogue of the Wronskian lemma}
%%%%%%%
%%%%%%%%%%%%%%%%%%%%%%%%%%%%%%%%%%%%%%%%%%%%%%%%%%%%%%%%%%%%%%%%%%%%%%%%%

\L
{We assume that $q$ is {\it not} a primitive root of unity of order $\leq\mu$ 
and that the ring of constants $C=\{f\in {\c F}:\varphi_q(f)=f\}$ is a field.  
Let $u_0,\dots, u_{\mu-1}\in{\c F}$, then 
$$
\bigdisplay{\dim_C\sum_{i=0}^{\mu-1}Cu_i
={\rm rank}~Cas(u_0,\dots, u_{\mu-1})\ ,}
$$
where $Cas(u_0,\dots, u_{\mu-1})$ is the so-called Casorati matrix 
$$
Cas(u_0,\dots, u_{\mu-1})=\pmatrix{
u_0&\cdots &u_{\mu-1}\cr
\varphi_q u_0&\cdots &\varphi_q u_{\mu-1}\cr
\vdots&\ddots&\vdots\cr
\varphi_q^{\mu-1}u_0&\cdots &\varphi_q^{\mu-1}u_{\mu-1}}\ .
$$}
\label{wronskiano}

\vpar\noindent
{\bf Remark.}
Of course, if $(q-1)x$ is a unit of $\c F$, lemma \pfidelta\ 
implies that 
$$
{\rm rank}\pmatrix{
u_0&\cdots &u_{\mu-1}\cr
\dq u_0&\cdots &\dq u_{\mu-1}\cr
\vdots&\ddots&\vdots\cr
\dq^{\mu-1}u_0&\cdots &\dq^{\mu-1}u_{\mu-1}}
={\rm rank}\pmatrix{
u_0&\cdots &u_{\mu-1}\cr
\varphi_q u_0&\cdots &\varphi_q u_{\mu-1}\cr
\vdots&\ddots&\vdots\cr
\varphi_q^{\mu-1}u_0&\cdots &\varphi_q^{\mu-1}u_{\mu-1}}\ .
$$

\Dim
Obviously we have 
$$
\dim_C\sum_{i=0}^{\mu-1}Cu_i
\geq{\rm rank}\pmatrix{\upla{\varphi_q^j u}{0}{\mu-1}}_{j=0,\dots,\mu-1}\ .
$$
Let us suppose that the rank of $Cas(u_0,\dots, u_{\mu-1})$ is $<\mu$.
Changing the order of $\upla{u}{0}{\mu-1}$, we may assume that 
$$
r={\rm rank}\pmatrix{\upla{\varphi_q^j u}{0}{\mu-1}}_{j=0,\dots,\mu-1}=
{\rm rank}\pmatrix{\upla{\varphi_q^j u}{0}{r-1}}_{j=0,\dots,\mu-1}\ ,
\numero{rango}
$$
with $r<\mu$.
It is enough to show that $u_r$ is in $\sum_{i=0}^{r-1}Cu_i$.
By (\frango), there exists
$(\upla{a}{0}{r-1})\in{\c F}^r$, such that:
$$
\pmatrix{\upla{\varphi_q^j u}{0}{r-1}}_{j=0,\dots,\mu-1}
\pmatrix{a_0\cr\vdots\cr a_{r-1}}=
\pmatrix{\varphi_q^j u_r}_{j=0,\dots,\mu-1}\ .
\numero{Appunogroth}
$$
If we apply $\varphi_q$ to (\fAppunogroth) and substract the
expression obtained from (\fAppunogroth) we get 
$$
\pmatrix{\varphi_q^j u_0(x),\dots,\varphi_q^j u_{r-1}(x)}_{j=1,\dots,\mu-1}
\pmatrix{\varphi_q(a_0)-a_0\cr\vdots\cr\varphi_q(a_{r-1})-a_{r-1}}=0
$$
Hence $(\varphi_q a_0,\dots,\varphi_q
a_{r-1})=(\upla{a}{0}{r-1})$ 
and therefore $a_i\in C$.
\QED

%%%%%%%%%%%%%%%%%%%%%%%%%%%%%%%%%%%%%%%%%%%%%%%%%%%%%%%%%%%%%%%%%%%%%%%%%
%%%%%%%
\newsubchapter{The {\titit q}-analogue of the cyclic vector lemma}
%%%%%%%
%%%%%%%%%%%%%%%%%%%%%%%%%%%%%%%%%%%%%%%%%%%%%%%%%%%%%%%%%%%%%%%%%%%%%%%%%

\vpar
The following $q$-analogue of the classical cyclic vector lemma 
for differential modules is a classical result 
({\it cf.} [S2, annexe B.2] and the very old references cited therein). 
It can also be deduced from the theory of skew fields 
({\it cf.}\ for instance [Ch]).  We prefer to give an elementary proof here, 
following [DGS, III, 4.2].

\L
{Let us assume that ${\c F}$ is a field of characteristic zero and 
that $q$ is not a root of unity.
Let $(M,\Phi_q)$ be a $q$-difference module of rank $\mu$
over ${\c F}$. Then there exists a cyclic
vector $m\in M$, {\it i.e.} an element $m$ such that $(m,\Phi_q(m),\dots,\Phi_q^{\mu-1}(m))$ is 
an $\c F$-basis of $M$.}
\label{vettorciclico}

\vpar\noindent
{\bf Remark.}
By (\fDfi), if $m$ is a cyclic vector for $\Phi_q$, 
then it is also a cyclic vector with respect to the operator
$\Dq$. 

\Dim
Let us denote the exterior product by $\wedge$.
Let 
$$
\nu=\max\{l\in\interi:\exists m\in M~\hbox{s.t.}~
m\wedge\Phi_q(m)\wedge\dots\wedge\Phi_q^{l-1}(m)\neq 0\}\ ;
$$
we suppose that $\nu$ is smaller than $\mu$ and we choose $m\in M$ such that 
$$
m\wedge\Phi_q(m)\wedge\dots\wedge\Phi_q^{\nu-1}(m)\neq 0\ .
$$
For all $\la\in C$, $s\in \interi$, $s\geq 1$, and $m^\p\in M$,  
there exists $m_i\in\wedge^\nu M$, for $i=0,\dots,\nu$, such that we have 
$$\eqalign{0
&=(m+\la x^s m^\p)\wedge\Phi_q(m+\la x^sm^\p)
    \wedge\dots\wedge\Phi_q^\nu(m+\la x^s m^\p)\cr 
&=m_0+m_1\la +\dots+m_\nu\la^\nu\ .} 
$$
Since the field of constants $C$ is
infinite, we have $m_0=\dots=m_\nu=0$; in particular 
$$
m_1=x^s\l(\sum_{i=0}^\nu q^{si}
    m\wedge\Phi_q(m)\wedge\dots\wedge\Phi_q^{i-1}(m)\wedge\Phi_q^i(m^\p)\wedge\Phi_q^{i+1}(m)
    \wedge\dots\wedge\Phi_q^{\nu}(m)\r)=0 
$$
for all positive integers $s$. 
It follows that for all $m^\p\in M$ and all $i=0,\dots,\nu$ we have 
$$
m\wedge\Phi_q(m)\wedge\dots\wedge\Phi_q^{i-1}(m)\wedge\Phi_q^i(m^\p)\wedge\Phi_q^{i+1}(m)
    \wedge\dots\wedge\Phi_q^{\nu}(m)=0\ .
$$ 
In particular, for $i=\nu$ we obtain 
$$
m\wedge\Phi_q(m)\wedge\dots\wedge\Phi_q^{\nu-1}(m)\wedge
\Phi_q^{\nu}(m^\p)=0\ ,\ \forall m^\p\in M, 
$$
which implies that 
$m\wedge\Phi_q(m)\wedge\dots\wedge\Phi_q^{\nu-1}(m)=0$. 
This contradicts the premises and hence $\nu=\mu$.
\QED

%%%%%%%%%%%%%%%%%%%%%%%%%%%%%%%%%%%%%%%%%%%%%%%%%%%%%%%%%%%%%%%%%%%%%%%%%
%%%%%%%
\newsubchapter{Formal classification of {\titit q}-difference modules}
%%%%%%%
%%%%%%%%%%%%%%%%%%%%%%%%%%%%%%%%%%%%%%%%%%%%%%%%%%%%%%%%%%%%%%%%%%%%%%%%%

\vpar
We recall the definition of regular singularity in the $q$-difference case, 
when {\it $K=R$ is a field of characteristic zero 
and $q$ is not a root of unity}.
Let ${\c M}=(M,\Phi_q)$ be a $q$-difference module over $K((x))$ 
of finite rank $\mu$.

\Def
{One says that ${\c M}$ is regular singular 
if there exists a $K((x))$-basis $\u e$ of $M$ in which the matrix $A(x)$ of $\Phi_q$ 
({\it a priori} an element of $Gl_\mu\l(K((x))\r)$) 
belongs to $Gl_\mu(K\[[x\]])$. }
\label{defreg}

\Oss
It is in fact equivalent to require the existence of 
a basis $\u e$ in which the matrix of $\Phi_q$ is a constant matrix 
({\it cf.}\ {\rm [PS, Ch. 12]}). 
\par 
One can say more: 
if $\Phi_q(\u e)=\u eA(x)$ with 
$A(x)\in Gl_\mu(K\[[x\]])$ 
and any couple $\a,\be$ of eigenvalues of $A(0)$ is such that either $\a=\be$
or $\a\be^{-1}\not\in q^\interi$, 
then we can find a basis $\u f$ of $M$ over $K((x))$ such that 
$\Phi_q(\u f)=\u fA(0)$.
Observe that if ${\c M}$ is regular singular it is always possible to find a
basis $\u e$ satisfying these properties ({\it cf.} [S2, 1.1.1] for a ditailed
proof). This remark will be useful in \regolarita. 
\label{Ao}

\Def
{The {\it exponents of a regular singular $q$-difference module ${\c M}$, with respect to 
a given basis $\u e$} as in the definition above,  
are the $q$-orbits $q^\interi a$ of the eigenvalues $a$ of $A(0)$.}

\vpar
Let us consider an extension of $K((x))$ of the form $L((t))$, 
where $x=t^d$ and $L$ is a finite extension of $K$ containing 
a root $\widetilde q$ of $q$ of order $d$. 
Then $\varphi_q$ extends to $L((t))$ in the following way:
$$\matrix{
\varphi_{\widetilde q }:& L((t))&\lrarrow & L((t))\cr\cr 
& t & \longmapsto & \widetilde q t}\ .
$$
The module $L((t))\otimes_{K((x))}M$, equipped with the operator  
$$\matrix{
\Phi_{\widetilde q }:&L((t))\otimes_{K((x))}M
   &\lrarrow &L((t))\otimes_{K((x))}M\cr\cr 
& f(t)\otimes m& \longmapsto & 
   \varphi_{\widetilde q }\l(f(t)\r)\otimes \Phi_q(m)} 
$$
is a $\widetilde q $-difference module over $L((t))$.
We recall the following result that will be useful in the sequel:

\T{{\rm [P, Cor. 9 and \S 9, 3)]}
Let $K$ be a field of characteristic zero, $q$ not a root of unity, 
${\c M}$ a $q$-difference module over $K((x))$ of rank $\mu$. 
Then there exists a divisor $d$ of $\mu!$ and a finite extension $L((t))$ of $K((x))$ 
as above, such that the 
$\widetilde q $-difference module $L((t))\otimes_{K((x))}M$ 
has an $L((t))$-basis $\u e$ with the following property:
the matrix $A(t)$ defined by $\Phi_{\widetilde q }(\u e)=\u eA(t)$ 
is a diagonal block matrix and each block has the form 
$$
t^{1-\la_i}
\pmatrix{\a_i& 0 & 0 & \cdots & 0\cr
         1 & \a_i& 0 & \cdots & 0\cr
         0 & \ddots &\ddots &\ddots & \vdots \cr
         \vdots & \ddots &\ddots & \ddots & 0\cr
         0 & \cdots & 0 & 1 &\a_i\cr}\ \ ,
$$
where 
$\la_i\in(1/d)\interi$, $\a_i\in\sum_{h=0}^d{\a_{i,h}\over t^h}$, 
with $\a_{i,h}\in L$ and $\a_{i,d}\neq 0$.
\tpar
The matrix $A(t)$ is unique up to permutation of the blocks.} 
\label{formclass}

\Oss
One can prove that a $q$-difference submodule of a regular singular $q$-difference
module is regular singular ({\it cf.}\ [P]).
\label{newton}

%%%%%%%%%%%%%%%%%%%%%%%%%%%%%%%%%%%%%%%%%%%%%%%%%%%%%%%%%%%%%%%%%%%%%%%%%
%%%%%%%%%%%%%%%%%%%%%%%%%%%%  2  %%%%%%%%%%%%%%%%%%%%%%%%%%%%%%%%%%%%%%%%
%%%%%%%%%%%%%%%%%%%%%%%%%%%%%%%%%%%%%%%%%%%%%%%%%%%%%%%%%%%%%%%%%%%%%%%%%

\newchapter{Unipotent {\tiit q}-difference modules}
\chaplabel{curva}

%%%%%%%%%%%%%%%%%%%%%%%%%%%%%%%%%%%%%%%%%%%%%%%%%%%%%%%%%%%%%%%%%%%%%%%%%

\vpar
In this section $R$ is again an arbitrary commutative ring 
and {\it $q\in R$ is a root of unity}. 
Let $\kappa$ denote its order:
$$
\kappa=\min\{m\in\interi:m>0,\ q^m=1\}\ . 
\numero{kappa}$$
Let $\c F$ be a $q$-difference algebra over $R$ and 
$C=\{f\in{\c F}:\varphi_q(f)=f\}$ the ring of constants of $\cal F$. 
We notice that $\varphi_q^\kappa=id_{\c F}$.

\vpar\noindent
{\bf Remark.}
Let ${\c M}=(M,\Phi_q)$ be a $q$-difference module over $\c F$. 
Then the operator $\Phi_q^\kappa$ is an $\c F$-linear automorphism of $M$.

%%%%%%%%%%%%%%%%%%%%%%%%%%%%%%%%%%%%%%%%%%%%%%%%%%%%%%%%%%%%%%%%%%%%%%%%%
%%%%%%%
\newsubchapter{Trivial {\titit q}-difference modules}
%%%%%%%
%%%%%%%%%%%%%%%%%%%%%%%%%%%%%%%%%%%%%%%%%%%%%%%%%%%%%%%%%%%%%%%%%%%%%%%%%

\Def
{The $q$-difference module ${\c M}=(M,\Phi_q)$ over $\c F$ is {\it trivial} if 
it is isomorphic to a $q$-difference module  
of the form $(N\otimes_C{\c F},id_N\otimes\varphi_q)$, 
where $N$ is a free $C$-module. }
\label{deftriviale}

\P
{\tparnoind
1) If $\c M$ is trivial over $\c F$ then 
$\Phi_q^\kappa$ is the identity morphism. 
\tparnoind
2) Let $R$ be a field and ${\c F}=R(x)$. If $\Phi_q^\kappa$ is the
identity, then $\c M$ is trivial over $\c F$.}
\label{grothradici}

\Dim
\par\noindent
{\sl 1)} Let us suppose that  $\c M$ is trivial over $\c F$.
By hypothesis there exists a basis $\u e$ of $M$ over $\c F$ 
such that $\Phi_q(\u e)=\u e$, 
which implies that $\Phi_q^\kappa=1$.
\par\noindent
{\sl 2)}  
Since $R$ is a field we can consider the operator
$\Dq=(\Phi_q-1)/(q-1)x$ on $M$. 
When $q$ is a root of unity of order $\kappa$ the formula (\ffiD) simplifies to 
$$
\Phi_q^\kappa=1+(q-1)^\kappa x^\kappa\Dq^\kappa\ .
\numero{unipotenti}$$
Therefore, under the assumption $\Phi_q^\kappa=1$, 
the operator $\Dq$ is a $C$-linear nilpotent morphism of order $\kappa$. 
Let $\mu=\dim_{\c F}M$. 
There exists a basis $\u m=(\upla{m}{1}{\mu\kappa})$ of $M$ over $C=R(x^\kappa)$ such that 
the matrix of $\Dq$ with respect to $\u m$ is an upper triangular nilpotent matrix in canonical form. 
In particular, this implies that $\Dq(m_1)=0$. 
To conclude, it is enough to prove that there exists a basis 
$\u m^\p$ of $M$ over $\cal F$ such that $\Dq(\u m^\p)=0$, 
which is equivalent to $\Phi_q(\u m^\p)=\u m^\p$.
\par
If $\mu=1$, it is enough to choose $\u m^\p=(m_1)$. 
Let us suppose $\mu>2$.   
If for all $i=2,\dots,\mu\kappa$ we have 
$\Dq(m_i)=m_{i-1}$, then $\Dq$ would be a nilpotent $C$-linear morphism 
of order $\mu\kappa>\kappa$, therefore there exists $j\in\{2,\dots,\mu\kappa\}$ 
such that $\Dq(m_j)=0$. We can suppose $j=2$. 
Repeating the reasoning we find that $\Dq(m_1)=\dots=\Dq(m_\mu)=0$.
We want to show that $(m_1,\dots,m_\mu)$ is a basis of $M$ over ${\c F}$.
Let us suppose that $\sum_{i=1}^\mu a_i(x)m_i=0$, with $a_i(x)\in R[x]$, 
$a_1(x)\neq 0$, and 
that the degree $\deg_x a_1(x)$ of $a_1(x)$ with respect to $x$ is minimal. Then 
$\Dq\l(\sum_{i=1}^\mu a_i(x)m_i\r)
=\sum_{i=1}^\mu \dq(a_i)(x)m_i=0$, with $\deg_x \dq(a_1)(x)\leq\deg_x a_1(x)-1$, 
so we get a contradiction.
Finally, we have found a basis $\u m^\p=(m_1,\dots,m_\mu)$ of $M$ over ${\c F}$ such that 
$\Dq(\u m^\p)=0$.
\QED

\vpar
If $(q-1)x$ is a unit of $\c F$, the operator 
$\Dq$ is defined over $M$ and 
(\funipotenti) shows that:

\C
{The operator $\Phi_q^\kappa$ is
unipotent if and only if $\Dq^\kappa$ is nilpotent.}
\label{unipotenza}

%%%%%%%%%%%%%%%%%%%%%%%%%%%%%%%%%%%%%%%%%%%%%%%%%%%%%%%%%%%%%%%%%%%%%%%%%
%%%%%%%
\newsubchapter{Extensions of trivial {\titit q}-difference modules}
%%%%%%%
%%%%%%%%%%%%%%%%%%%%%%%%%%%%%%%%%%%%%%%%%%%%%%%%%%%%%%%%%%%%%%%%%%%%%%%%%

\P
{\tparnoind
1) If the $q$-difference module $\c M$ is an extension of 
trivial $q$-difference modules then 
the ${\c F}$-linear morphism $\Phi_q^\kappa$ is unipotent. 
\tparnoind
2) If $R$ is a field , ${\c F}=R(x)$  
and $\Phi_q^\kappa$ is unipotent, then $\c M$ is an
extension of trivial $q$-difference modules.}
\label{trivialext}

\Dim
\par\noindent
{\sl 1)} We have to prove that $\Phi_q^\kappa-1$ is a nilpotent
endomorphism. 
If $\c M$ is an extension of trivial $q$-difference modules over
${\c F}$, by (\pgrothradici) we can find a basis $\u e$ of $M$ over ${\c F}$ such that 
$\Phi_q^\kappa\u e=\u e({\b I}+H_\kappa(x))$, 
where ${\b I}$ is the identity matrix and $H_\kappa(x)$ is a block
matrix of the form 
$$
H_\kappa(x)=\pmatrix{
\matrix{0 \cr \hrulefill \cr 0\cr \hrulefill\cr 0}&\vrule&
\matrix{ \ast\cr  \cr \ddots \cr \cr 0}&\vrule&
\matrix{\ast \cr \hskip -60 pt\hrulefill \cr \ast    \cr 
     \hskip -60pt\hrulefill\cr 0}}\ .
$$ 
The matrix $H_\kappa(x)$ is nilpotent, hence $\Phi_q^\kappa$ is a 
unipotent endomorphism.
\par\noindent
{\sl 2)} If $R$ is a field we can consider the operator $\Dq$
associated to $\Phi_q$. 
Since $\Phi_q^\kappa$ is unipotent,
the $C$-linear morphism $\Dq$ is nilpotent ({\it cf.}\ (\punipotenza)), 
therefore there exists $m_1\in M$ such that $\Dq m_1=0$. 
Let $\mu=\dim_{\c F}M$. 
If $\mu=1$ there is nothing more to prove. 
Let $\mu$ be greater than $1$. The operator $\Dq$ induces a structure
of a $q$-difference module over the quotient ${\c F}$-vector
space $M/{\c F}m_1$ which satisfies the hypothesis. 
By induction we can find a filtration of 
$M/{\c F}m_1$ 
$$
\widetilde M_0=\{0\}\subset \widetilde M_1\subset\dots\subset
\widetilde M_l=M/{\c F}m_1\ ,
$$
such that:
\par\noindent
1) for all $i=0,\dots,l$ the sub-vector space $\widetilde M_i$ is
stable by the operator induced by $\Dq$;
\par\noindent 
2) for all $i=1,\dots,l$ the quotient module $\widetilde M_i/\widetilde M_{i-1}$ equipped
with its natural structure of a $q$-difference module is trivial over
${\c F}$.
\par\noindent 
Let $\iota: M\lrarrow M/{\c F}m_1$ be the canonical projection. Then  
$$
M_{-1}=\{0\}\subset M_0=\iota^{-1}(\widetilde M_0)={\c F}m_1\subset
M_1=\iota^{-1}(\widetilde M_1)
\subset\dots\subset
M_l=\iota^{-1}(\widetilde M_l)=M
$$
satisfies the conditions:
\par\noindent
1) for all $i=-1,0,\dots,l$ the sub-vector space $M_i$ is
stable under the operator induced by $\Dq$;
\par\noindent 
2) for all $i=0,\dots,l$ the quotient module 
$M_i/M_{i-1}\cong \widetilde M_i/\widetilde M_{i-1}$ equipped
with its natural structure of a $q$-difference module is trivial over
${\c F}$.
\QED

\Oss
In [H, Ch. 6] we find a classification of $q$-difference modules 
over $R(x)$ when $q$ is a root of unity and $R$ is a field of characteristic zero. 
The author defines the Galois group
associated to a linear $q$-difference module and proves that 
it is the smallest algebraic group over $R(x^\kappa)$ containing
$\Phi_q^\kappa$. 

%%%%%%%%%%%%%%%%%%%%%%%%%%%%%%%%%%%%%%%%%%%%%%%%%%%%%%%%%%%%%%%%%%%%%%%%%
%%%%%%%%%%%%%%%%%%%%%%%%% PART II     %%%%%%%%%%%%%%%%%%%%%%%%%%%%%%%%%%%
%%%%%%%%%%%%%%%%%%%%%%%%%%%%%%%%%%%%%%%%%%%%%%%%%%%%%%%%%%%%%%%%%%%%%%%%%
%%%%%%%%%%%%%%%%%%%%%%%%%%%%%%%%%%%%%%%%%%%%%%%%%%%%%%%%%%%%%%%%%%%%%%%%%

\vskip 20pt 
\centerline{\titlepart Part II.  {\titlepartit p}-adic methods}
%%%%%%%%%%%%%%%%%%%%%%%%%%%%%%%%%%%%%%%%%%%%%%%%%%%%%%%%%%%%%%%%%%%%%%%%%
%%%%%%%%%%%%%%%%%%%%%%%%%%%%%%%%%%%%%%%%%%%%%%%%%%%%%%%%%%%%%%%%%%%%%%%%%
%%%%%%%%%%%%%%%%%%%%%%%%%%%%%%%%%%%%%%%%%%%%%%%%%%%%%%%%%%%%%%%%%%%%%%%%%
%%%%%%%%%%%%%%%%%%%%%%%%%%%%%%%%%%%%%%%%%%%%%%%%%%%%%%%%%%%%%%%%%%%%%%%%%

%%%%%%%%%%%%%%%%%%%%%%%%%%%%%%%%%%%%%%%%%%%%%%%%%%%%%%%%%%%%%%%%%%%%%%%%%
%%%%%%%%%%%%%%%%%%%%%%%%%%%%  3  %%%%%%%%%%%%%%%%%%%%%%%%%%%%%%%%%%%%%%%%
%%%%%%%%%%%%%%%%%%%%%%%%%%%%%%%%%%%%%%%%%%%%%%%%%%%%%%%%%%%%%%%%%%%%%%%%%

\newchapter{Considerations on the differential case}
\chaplabel{differenziale}

\vpar
We would like to recall some properties of arithmetic differential modules 
that are supposed to motivate the structures we will introduce in the sequel. 
In particular the considerations below show that the two notions of
nilpotent reduction introduced in  
\nilpred\ are both natural $q$-analogues of the notion of nilpotent reduction 
for differential modules. 

\vpar
Let us consider the field of rational numbers $\raz$. For all prime
$p\in\interi$ we consider the $p$-adic norm $|~|_p$ over $\raz$,
normalized so that $|p|_p=p^{-1}$. 
By the Gauss lemma, the norm $|~|_p$ can be extended to
the $p$-adic Gauss norm $|~|_{p,Gauss}$ over $\raz(x)$, by setting  
$$
\l|{\sum_{i=0}^na_ix^i\over \sum_{j=0}^mb_jx^j}\r|_{p,Gauss}=
{\sup_{i=0,\dots,n}|a_i|_p\over\sup_{j=0,\dots,m}|b_j|_p}\ .
$$
Let us consider a differential module $(M,\Delta)$ over $\raz(x)$, 
{\it i.e.} a $\raz(x)$-vector space $M$ of finite dimension $\mu$
equipped with a $\raz$-linear morphism $\Delta:M\lrarrow M$ such that 
$\Delta(fm)={d f\over dx}m+f\Delta(m)$, for all $f\in\raz(x)$ and $m\in M$. 
We fix a basis $\u e$ of $M$ over $\raz(x)$ and set 
$\Delta^n\u e=\u eG_n(x)$, for all $n\geq 1$, 
where $G_n(x)\in M_{\mu\times\mu}(\raz(x))$ is a
square matrix of order $\mu$ with coefficients in $\raz(x)$. 
For $n=0$ we set $G_0(x)={\b I}_\mu$. 
The matrix 
$\sum_{n\geq 0}{G_n(t)\over n!}(x-t)^n\in Gl_\mu(\raz(t)\[[x-t\]])$
is a formal solution of ${dY\over dx}=YG_1(x)$ at any $t$ in the
algebraic closure of $\raz$, such that $G_1(x)$ has no pole at $t$.
\par
For almost all primes $p\in\interi$ we have $|G_1(x)|_{p,Gauss}\leq 1$
and we can consider the image of $G_1(x)$ in 
$M_{\mu\times\mu}({\b F}_p(x))$. One usually says that $(M,\Delta)$ has
{\it $p$-adic nilpotent reduction of order $n$} if
$|G_{np}(x)|_{p,Gauss}<1$ or, equivalently, if 
$G_{np}(x)\equiv 0$ modulo $p$. If $n=1$ we say that $(M,\Delta)$ has 
{\it $p$-curvature zero}. 
\par
Another equivalent statement of the Grothendieck conjecture 
({\it cf.}\ Introduction)  is:

\vpar\noindent
{\bf Grothendieck's conjecture. }
{\sl If $(M,\Delta)$ has $p$-curvature zero for almost all primes $p$,
then $(M,\Delta)$ becomes trivial over an algebraic extension of $\raz(x)$.}

\vpar
Let us consider the 
{\it $p$-adic generic radius of convergence}
$$
R_p(M)=\inf\l(1, 
\liminf_{n\rightarrow\infty}\l|{G_n(x)\over n!}\r|_{p,Gauss}^{-1/n}\r)\ .
$$
If $(M,\Delta)$ has $p$-adic nilpotent reduction of order $n$ we have:
$$
R_p(M)\geq p^{1/np}p^{-1/(p-1)}\ ;
$$
in particular if $(M,\Delta)$ has $p$-curvature zero, the previous
inequality reduces to $R_p(M)\geq p^{-1/p(p-1)}$.
An important and useful property of arithmetic differential modules is that 
({\it cf.}\ (\stepb))
$$
\sum_{\hbox{\sevenrm $\scriptstyle p$-curvature zero}}\log{1\over R_p(M)}\leq 
\sum_{\hbox{\sevenrm $\scriptstyle p$-curvature zero}}{\log p\over p(p-1)}<\infty\ .
$$
For the case of $q$-difference
modules, a naive translation of these definitions gives deceiving results.
A more accurate analysis of the
case of $p$-curvature zero leads to the following remark. Let
$\o G_1(x)$ be the image of $G_1(x)$ in 
$M_{\mu\times\mu}({\b F}_p(x))$. By imposing 
that $G_p(x)\equiv 0$ modulo $p$ we are actually requiring that the differential
system in positive characteristic 
$$
{dY\over dx}=Y\o G_1(x)
$$
has a fundamental matrix of solutions $Y(x)\in Gl_\mu({\b F}_p(x))$. 
Since the derivation ${1\over p!}{d^p\over dx^p}$ makes sense in
characteristic $p$ this implies that 
${1\over p!}{d^pY\over dx^p}\equiv {G_p(x)\over p!}Y$ modulo p, with 
$|G_p(x)|_{p,Gauss}\leq p^{-1}=|p!|_p$.
\par
The problem is that in the $q$-difference case, one can define some
$q$-analogue of factorials, but they generally are $p$-adically
smaller than the uniformizer $p$. It turns out that the $q$-analogue of the condition
$G_p(x)\equiv 0$ modulo $p$ is not equivalent to the $q$-analogue of
the condition $\l|{G_p(x)\over p!}\r|_{p,Gauss}\leq 1$, but both of
them are linked to the property of a suitable $q$-difference system having
a fundamental solution matrix in some polynomial ring 
over a quotient of $\interi$.
Therefore, for a $q$-difference module, we have two natural notions  of
nilpotent reduction: from a {\it local} point of view the notions are not
equivalent ({\it cf.}\ \nilpred), but we conjecture that they are {\it globally} equivalent. 

%%%%%%%%%%%%%%%%%%%%%%%%%%%%%%%%%%%%%%%%%%%%%%%%%%%%%%%%%%%%%%%%%%%%%%%%%
%%%%%%%%%%%%%%%%%%%%%%%%%%%%  4  %%%%%%%%%%%%%%%%%%%%%%%%%%%%%%%%%%%%%%%%
%%%%%%%%%%%%%%%%%%%%%%%%%%%%%%%%%%%%%%%%%%%%%%%%%%%%%%%%%%%%%%%%%%%%%%%%%

\newchapter{Introduction to {\tiit p}-adic {\tiit q}-difference modules}

%%%%%%%%%%%%%%%%%%%%%%%%%%%%%%%%%%%%%%%%%%%%%%%%%%%%%%%%%%%%%%%%%%%%%%%%%

\vpar
Let $K_v$ be a field of characteristic zero, 
complete with respect to a non-archimedean norm $|~|_v$.
Let ${\c V}_v$ be the ring of integers of $K_v$, 
$\varpi_v$ the uniformizer of ${\c V}_v$, $k_v$ its residue field of characteristic
$p>0$. 
\par
{\it We fix a nonzero element $q\in K_v$, such that $q$ is not a root of
unity and $|q|_v=1$.}
Let $\o q$ be the image of $q$ in $k_v$; 
we suppose that $\o q$ is algebraic over the prime field ${\b F}_p$ and we set 
$$
{\kappa_v}=\min\{m\in\interi:m>0,\ \o q^m=1\}\geq 1\ .
$$
We notice that $\o q\in k_v$ satisfies the assumption of \nilpotenza.
\par
In addition, we assume that 
$$
|1-q^{\kappa_v}|_v<|p|_v^{1/(p-1)} 
$$  
If $q$ is an element of $\raz_p\subset K_v$ and $p>2$, this holds
automatically in fact $|1-q^{\kappa_v}|_v\leq|p|_v<|p|_v^{1/(p-1)}$.

%%%%%%%%%%%%%%%%%%%%%%%%%%%%%%%%%%%%%%%%%%%%%%%%%%%%%%%%%%%%%%%%%%%%%%%%%
%%%%%%%
\newsubchapter{{\titit p}-adic estimates of {\titit q}-binomials}
%%%%%%%
%%%%%%%%%%%%%%%%%%%%%%%%%%%%%%%%%%%%%%%%%%%%%%%%%%%%%%%%%%%%%%%%%%%%%%%%%

\L
{Let $n\geq i\geq 0$ be two integers. 
We have 
$$                
|\nq!|_v=|[{\kappa_v}]_q|_v^{\l[n\over{\kappa_v}\r]}
    \l|\l[n\over{\kappa_v}\r]!\r|_v\ , 
\numero{stimafattoriale}$$
where $[x]$ is the integer part of $x\in\reali$, and 
$$
\l|{n\choose i}_q\r|_v\leq 1\ .
\numero{stimabinomiale}$$}
\label{ext}

\Dim
\par\noindent
{\sl 1)} By the definition of ${\kappa_v}$, if ${\kappa_v}$ does not divide $n$, $|1-q^n|_v=1$. 
Since $|1-q^{n{\kappa_v}}|_v\leq|1-q^{\kappa_v}|_v<|p|_v^{1/(p-1)}$ 
for all $n\in\interi$, $n\geq 1$, we have 
({\it cf.}\ for instance [DGS, II, 1.1]) 
$$
|1-q^{n{\kappa_v}}|_v=|\log q^{n{\kappa_v}}|_v=|n\log q^{\kappa_v}|_v
=|n|_v|1-q^{\kappa_v}|_v\ ,
\numero{stimanq}$$
so $|[n{\kappa_v}]_q|_v=\l|{1-q^{n{\kappa_v}}\over 1-q}\r|_v
=|n|_v|[{\kappa_v}]_q|_v$.
We obtain 
$$\eqalign{|[n]_q!|_v
&=|[{\kappa_v}]_q|_v^{\l[n\over{\kappa_v}\r]}
    \prod_{i\leq n\atop {\kappa_v}\vert i}|i|_v\cr
&=|[{\kappa_v}]_q|_v^{\l[n\over{\kappa_v}\r]}
    |{\kappa_v}|_v^{\l[n\over{\kappa_v}\r]}\l|\l[n\over{\kappa_v}\r]!\r|_v\ .}
$$
Since ${\kappa_v} $ is a divisor of 
$p^s-1$ for a suitable integer $s\geq 1$, we have $({\kappa_v},p)=1$,
which implies that $|{\kappa_v}|_v=1$.
\par
\vbox{\noindent
{\sl 2)} Since $q$ is an invertible element of the 
ring of integers ${\c V}_v$, we obtain the inequality $\l|{n\choose i}_q\r|_v\leq 1$ 
using the relation 
$$
(1-x)_n=\sum_{j=0}^n(-1)^j{n\choose j}_q q^{j(j-1)/2}x^j
\in {\c V}_v[x]\ .
$$
\QED}

%%%%%%%%%%%%%%%%%%%%%%%%%%%%%%%%%%%%%%%%%%%%%%%%%%%%%%%%%%%%%%%%%%%%%%%%%
%%%%%%%
\newsubchapter{The Gauss norm  and the invariant $\chi_v({\c M})$}
%%%%%%%
%%%%%%%%%%%%%%%%%%%%%%%%%%%%%%%%%%%%%%%%%%%%%%%%%%%%%%%%%%%%%%%%%%%%%%%%%

\vpar
By the Gauss lemma, one can 
extend the norm $|~|_v$ to the so-called Gauss norm $|~|_{v,Gauss}$ over 
$K_v(x)$ by setting 
$$
\l|{\sum_{i=0}^na_ix^i\over \sum_{j=0}^mb_jx^j}\r|_{v,Gauss}=
{\sup_{i=0,\dots,n}|a_i|_v\over\sup_{j=0,\dots,m}|b_j|_v}\ .
$$
Remark that $|~|_{v,Gauss}$ is multiplicative.

\L
{For any $f(x)\in K_v(x)$ and any positive integer $n$,  
we have 
$$
\l|{\dq^n\over[n]_q!}f(x)\r|_{v,Gauss}\leq|f(x)|_{v,Gauss}\ .
$$}
\label{normaderivazioni}

\Dim
By (\pprimoesempio) and (\pext) the inequality holds for all $f(x)\in K_v[x]$.  
Furthermore we have 
$$
\l|\dq\l({1\over f(x)}\r)\r|_{v,Gauss}=
\l|{\dq f(x)\over f(x)f(qx)}\r|_{v,Gauss}\leq 
\l|{1\over f(x)}\r|_{v,Gauss}\ .
$$
By the $q$-analogue of the Leibniz formula (\fgrothleibniz) we have 
$$
{\dq^n\over[n]_q!}\l({1\over f(x)}\r)=-{1\over f(q^nx)}
\sum_{i=0}^{n-1}{\dq^i\over [i]_q!}(f)(q^{n-i}x)
{\dq^{n-i}\over[n-i]_q!}\l({1\over f(x)}\r)\ ;
$$
and therefore by induction 
$$
\l|{\dq^n\over[n]_q!}\l({1\over f(x)}\r)\r|_{v,Gauss}\leq 
\l|{1\over f(x)}\r|_{v,Gauss}\ .
$$
\par\vbox{\noindent
Finally, if $g(x)\in K_v[x]$ we obtain 
$$
\l|{\dq^n\over[n]_q!}\l({g(x)\over f(x)}\r)\r|_{v,Gauss}
=\l|\sum_{i=0}^{n}{\dq^i\over [i]_q!}(g)(q^{n-i}x)
{\dq^{n-i}\over[n-i]_q!}\l({1\over f(x)}\r)\r|_{v,Gauss}
\leq\l|{g(x)\over f(x)}\r|_{v,Gauss}\ .
$$
\QED}

\vpar
Let ${\c M}=(M,\Phi_q)$ be a $q$-difference module over $K_v(x)$. 
Since $K_v$ is a field, we have well defined operators 
$$
\dq={\varphi_q-1\over (q-1)x}\ 
\hbox{and}\ \Dq={\Phi_q-1\over (q-1)x}
$$
acting over $K_v(x)$ and $M$, respectively ({\it cf.}\ (\pderivazionibendefinite)).
\par
We fix a basis $\u e$ of $M$ over $K_v(x)$ and define a sequence of matrices 
$G_n(x)\in M_{\mu\times\mu}(K_v(x))$ for any integer $n\geq 0$ by setting 
$$
\Dq^n(\u e)=\u eG_n(x)\ .
\numero{Dn}$$
The matrices $G_n(x)$ satisfy the inductive relation 
$$
G_0(x)={\b I}_\mu,\ G(x)=G_1(x),\ G_{n+1}(x)=G_1(x)G_n(qx)+\dq G_n(x)\ .
\numero{coefficienti}$$
We call 
$$
\dq Y=YG(x) 
\leqno{({\cors S})}
$$
the {\it $q$-difference system associated 
to $\c M$ with respect to the basis $\u e$}.  
If zero is not a pole of $G_1(x)$, 
we obtain a formal solution of $({\cors S})$:  
$\sum_{n=0}^\infty{G_n(0)\over[n]_q!}x^n$.
More generally, if $q^na$ is not a pole of $G_1(x)$ for all positive integers $n$, 
the matrix 
$$
Y(x)=\sum_{n=0}^\infty{G_n(a)\over[n]_q!}(x-a)_n\in
M_{\mu\times\mu}(K_v\[[x-a\]]_q) 
\numero{solformale}$$
is a formal solution of $({\cors S})$\footnote{$\hskip -3pt{}^{(1)}$}
{\rm The ring $K_v\[[x-a\]]_q$ is neither integral nor local: in particular  
the remark that $Y(a)={\Bbb I}_\mu$ is not enought to conclude that 
$Y(x)\in Gl_\mu(K_v\[[x-a\]]_q)$ (which is true, anyway). 
Moreover there exists subalgebra of $K_v\[[x-a\]]_q$ that can be
identified to the algebra af analytic function over a $q$-invariant analytic
domain, eventually non connected, via a $q$-Taylor expansion. 
As a consequence, it is possible to generalize (\analitico) and (\transfert)
below, obtaining more precise statements. 
The proofs of
these facts, that are not relevant for the sequel,  
are in a paper in preparation by the author on the $p$-adic theory
of $q$-difference eqauations.}.
\par
One can easily check that if ${\c Y}$ is a fundamental 
matrix for $({\cors S})$ with coefficients in some fixed $q$-difference algebra ${\c F}$ 
({\it i.e.} an invertible matrix solution of $({\cors S})$ with coefficients in ${\c F}$), 
any other fundamental matrix of $({\cors S})$ in $Gl_\mu({\c F})$ 
is of the form ${\c Y}F$, where 
$F$ is an invertible matrix with coefficients in the subring of 
constants ${\c F}$.

\vpar
In the following definition, 
the sup-norm of a matrix is the maximum of the norms of its
entries: 

\Def
{$\displaystyle 
\chi_v({\c M})=\inf\l(1,\liminf_{n\rightarrow\infty}
\l|{G_n(x)\over [n]_q!}\r|_{v,Gauss}^{-1/n}\r)$.}

\L
{Let 
$$
h(n)=\sup_{s\leq n}\lgp\l|{G_s(x)\over [s]_q!}\r|_{v,Gauss}=
\sup_{s\leq n}\log\l|{G_s(x)\over [s]_q!}\r|_{v,Gauss}\ ,
$$
with $\lgp x=\log\sup(x,1)$, for all $x\in\reali$.
Then 
$$
\limsup_{n\rightarrow\infty}{h(n)\over n}=\log{1\over\chi_v({\c M})}\ .
$$
Moreover, $\chi_v({\c M})$ is independent of the choice 
of the $K_v(x)$-basis $\u e$ of $M$.}
\label{buonadef}

\Dim
We recall that $G_0(x)={\b I}_\mu$ and that therefore the two
definitions of $h(n)$ are equivalent. 
\par
Let $\u f=\u e F(x)$ be another $K_v(x)$-basis of $M$, 
with $F(x)\in Gl_\mu(K_v(x))$. 
For any integer $n\geq 0$, we set:
$$\cases{
\Dq^n(\u e)=\u eG_n(x),
  &$\displaystyle
  h_{\u e}(n)=\sup_{s\leq n}\log\l|{G_s(x)\over [s]_q!}\r|_{v,Gauss}\ ;$\cr\cr
\Dq^n(\u f)=\u fH_n(x),
  &$\displaystyle
  h_{\u f}(n)=\sup_{s\leq n}\log\l|{H_s(x)\over [s]_q!}\r|_{v,Gauss}\ .$\cr}
$$
By (\fgrothleibniz) we have 
$$\eqalign{\u f{H_n(x)\over[n]_q!}
&={\Dq^n\over[n]_q!}(\u f)={\Dq^n\over[n]_q!}(\u eF(x))\cr
&=\u fF(x)^{-1}\sum_{i=0}^n{G_i(x)\over[i]_q!}{\dq^{n-i}(F)\over [n-i]_q!}
      (q^ix)\ ,\cr}
$$
and hence it follows that  
$$
h_{\u f}(n)\leq 
\log|F(x)^{-1}|_{v,Gauss}+\log|F(x)|_{v,Gauss}+h_{\u e}(n)\ .
\numero{dish}$$
By symmetry, we deduce that 
$$
\limsup_{n\rightarrow\infty}{h_{\u e}(n)\over n}=
\limsup_{n\rightarrow\infty}{h_{\u f}(n)\over n}\ .
$$
\par\vbox{\noindent
Let $h(n)=h_{\u e}(n)$. 
It is a general fact ({\it cf.}\ for instance the proof of [DGS, VII, Lemma 4.1]) that 
$$
\limsup_{n\rightarrow\infty}{h(n)\over n}=\lg{1\over\chi_v({\c M})}\ . 
$$
\QED}

\vpar
Let $a$ be an element of $K_v$ such that 
$q^na$ is not a pole of $G_1(x)$ for any $n\geq 0$.
We want to relate $\chi_v({\c M})$ and the radius of convergence of the matrix 
$\sum_{n=0}^\infty{G_n(a)\over[n]_q!}(x-a)_n$, which solves 
the linear $q$-difference system associated to $\c M$, 
with respect to the basis $\u e$.
First of all, we notice that if $|a|_v\leq 1$, we have 
$$
\l|{G_n(x)\over [n]_q!}\r|_{v,Gauss}\geq \l|{G_n(a)\over [n]_q!}\r|_v\ ,
$$ 
and therefore we obtain 
$$
\chi_v({\c M})\leq \liminf_{n\rightarrow\infty}
\l|{G_n(a)\over [n]_q!}\r|_v^{-1/n}\ .
$$
Hence, if zero is not a pole of $G_1(x)$, the matrix 
$\sum_{n=0}^\infty{G_n(0)\over[n]_q!}x^n$ converges at least for 
$|x|_v<\chi_v({\c M})$. 

\Ex
Let us consider the analogue of the exponential series 
$$
\exp_q(x)=\sum_{n=0}^\infty{x^n\over [n]_q!}\ .
$$
Obviously $\exp_q(x)$ is the solution at zero of the 
$q$-difference equation $\dq y=y$, which is the system associated 
to the $q$-difference module $(K_v(x),\Dq)$, with $\Dq(f(x))=\dq f(x)+f(qx)$ 
for all $f(x)\in K_v(x)$.
Then $\chi_v(K_v(x),\Dq)$ coincides with the radius of convergence of $\exp_q(x)$, 
that is $|[{\kappa_v}]_q|_v^{1/{\kappa_v}}|p|_v^{1/{\kappa_v}(p-1)}$, by (\fstimafattoriale). 
\label{exp}

\vpar 
If $a\neq 0$ the situation is slightly more complicated:

\L
{Let $\sum_{n=0}^\infty a_n(x-a)_n\in K_v\[[x-a\]]_q$ and 
let $\varrho\in(0,1]$ be a real number. 
Then if
$$
\sup(\varrho,|a|_v)^{1-{1\over{\kappa_v}}}
\sup(\varrho,|a|_v|[{\kappa_v}]_q|_v)^{1\over{\kappa_v}}
<\liminf_{n\rightarrow\infty}\l|a_n\r|_v^{-1/n}
$$
$$
\hbox{(in particular, if 
$\sup(\varrho,|a|_v)<\liminf_{n\rightarrow\infty}\l|a_n\r|_v^{-1/n}$)}
$$
the series $\sum_{n=0}^\infty a_n(x-a)_n$ converges in the disk
$\{x\in K_v:|x-a|_v<\varrho\}$.}
\label{raggio}

\C
{Let $a\in{\cal V}_v$ be such that $|a|_v\leq \chi_v({\c M})$,  
then $\sum_{n=0}^\infty{G_n(a)\over[n]_q!}(x-a)_n$
converges  in the open disk
$\{x\in K_v:|x-a|_v<\chi_v({\c M})\}$.}

\vpar\noindent
{\bf Proof of lemma \praggio.}
By the Maximum Modulus Principle [DGS, VI, 1.1] and (\fstimanq) we have 
$$\eqalign{\sup_{|x-a|_v<\varrho}|(x-a)_n|_v
&=\sup_{|x-a|_v<\varrho}\l|\l(x-a\r)\l(x-a+a(1-q)\r)
   \cdots\l(x-a+a\l(1-q^{n-1}\r)\r)\r|_v\cr
&\leq\varrho\sup\l(\varrho,|a|_v\r)^{(n-1)-\l[{n-1\over{\kappa_v}}\r]}
   \prod_{i=1}^{\l[{n-1\over{\kappa_v}}\r]}\sup\l(\varrho,|ai|_v|[{\kappa_v}]_q|_v\r)\cr
&\leq \varrho\sup\l(\varrho,|a|_v\r)^{(n-1)-\l[{n-1\over{\kappa_v}}\r]}
   \sup\l(\varrho,|a|_v|[{\kappa_v}]_q)|_v\r)^{\l[{n-1\over{\kappa_v}}\r]}\cr}
$$
\par\vbox{\noindent
Finally, $\sum_{n=0}^\infty a_n(x-a)_n$ converges if 
$$
\sup(\varrho,|a|_v)^{1-{1\over{\kappa_v}}}
\sup(\varrho,|a|_v|[{\kappa_v}]_q|_v)^{1\over{\kappa_v}}
<\liminf_{n\rightarrow\infty}\l|a_n\r|_v^{-1/n}\ .
$$
\QED}

\vpar
The following characterization of $\chi_v({\c M})$ is the 
$q$-analogue of a result by Andr{\'e} ({\it cf.}\ [A, IV, \S 5]):

\P
{The sequence $\l({h(n)\over n}\r)_{n\in\nat}$ defined in (\pbuonadef)
is convergent:
$$
\lim_{n\rightarrow\infty}{h(n)\over n}=\log{1\over\chi_v({\c M})}\ .
\numero{limh}$$}

\Dim
By  (\pbuonadef) it is enough to prove the existence of the limit.
Let $s,n$ be two positive integers; we have 
$$\eqalign{{\Dq^{s+n}\over [s+n]_q!}\l(\u e\r)
&={\Dq^{n}\over [s+n]_q!}\l(\Dq^s\u e\r)
     ={\Dq^{n}\over [s+n]_q!}\l(\u eG_s(x)\r)\cr
&=\u e{1\over [s+n]_q!}\sum_{i=0}^n
     {n\choose i}_qG_i(x)\dq^{n-i}(G_s)(q^ix)\cr
&=\u e\sum_{i+j=n}{[n]_q![s]_q!\over [s+n]_q!}
     {G_i(x)\over [i]_q!}{\dq^j\over [j]_q!}\l({G_s(q^ix)\over [s]_q!}\r)\cr}
$$
It follows that 
$$
{G_{s+n}(x)\over [s+n]_q!}=\sum_{i+j=n}{[n]_q![s]_q!\over [s+n]_q!}
     {G_i(x)\over [i]_q!}{\dq^j\over [j]_q!}\l({G_s(q^ix)\over [s]_q!}\r)
$$
and hence 
$$
\log\l|{G_{s+n}(x)\over [s+n]_q!}\r|_{v,Gauss}\leq
\log\l|{G_s(x)\over [s]_q!}\r|_{v,Gauss}+h(n)-\log\l|{n+s\choose s}_q\r|_v\ .
$$
For all $k\in\nat$ and $n\geq s$, by induction we obtain 
$$\eqalign{\log\l|{G_{s+kn}(x)\over [s+kn]_q!}\r|_{v,Gauss}
&\leq\log\l|{G_{s+(k-1)n}(x)\over[s+(k-1)n]_q!}\r|_{v,Gauss}+h(n)
     -\log\l|{s+kn\choose s+(k-1)n}_q\r|_v\cr
&\leq\log\l|{G_s(x)\over [s]_q!}\r|_{v,Gauss}+kh(n)
     -\log\l|\prod_{i=1}^k{s+in\choose s+(i-1)n}_q\r|_v\ .}
$$
Let $N\in\nat$, $N\geq n$ ; then $N=\l[{N\over n}\r]n+s$, 
with $0\leq s<n$, and the previous inequality becomes 
$$
\log\l|{G_N(x)\over [N]_q!}\r|_{v,Gauss}\leq 
\l(\l[{N\over n}\r]+1\r)h(n)-
\log\l|\prod_{i=1}^{\l[{N\over n}\r]}{s+in\choose s+(i-1)n}_q\r|_v\ .
$$
Since 
$\log\l|\prod_{i=1}^{\l[{N\over n}\r]}{s+in\choose s+(i-1)n}_q\r|_v\leq 0$ 
is a decreasing function of $N$ we obtain 
$$\eqalign{{h(N)\over N}
&\leq \l({1\over n}+{1\over N}\r)h(n)
      -\log\l|\prod_{i=1}^{\l[{N\over n}\r]}{s+in\choose s+(i-1)n}_q\r|_v\cr
&\leq \l({1\over n}+{1\over N}\r)h(n)
      -\log\l|{[N]_q!\over\l([n]_q!\r)^{\l[{N\over n}\r]}[s]_q!}\r|_v\ .}
$$
Finally we deduce by (\fstimafattoriale) that 
$$\eqalign{\limsup_{N\rightarrow\infty}{h(N)\over N}
&\leq\limsup_{N\rightarrow\infty}\l(\l({1\over n}+{1\over N}\r)h(n)
   -\log\l|{[N]_q!\over\l([n]_q!\r)^{\l[{N\over n}\r]}[s]_q!}\r|_v^{1\over N}\r)\cr
&\leq {h(n)\over n}
   -\log\l({|[{\kappa_v}]_q|_v^{1/{\kappa_v}}|p|_v^{1/{\kappa_v}(p-1)}
   \over \l|[n]_q!\r|_v^{1/n}}\r)\ .\cr}
$$
Therefore 
$$
\limsup_{N\rightarrow\infty}{h(N)\over N}
\leq \liminf_{n\rightarrow\infty}{h(n)\over n}\ ,
$$
from which it follows that the sequence $\l({h(n)\over n}\r)_{n\in\nat}$ is convergent.
\QED

\vpar
Now we prove a first estimate for $\chi_v({\c M})$. In succeeding sections 
we will prove a more precise estimate linked to the notion of unipotent reduction.

\P
{We have:
$$
\chi_v({\c M})\geq {|[{\kappa_v}]_q|_v^{1/{\kappa_v}}
|p|_v^{1/{\kappa_v}(p-1)}\over\sup\l(|G(x)|_{v,Gauss},1\r)}\ .
$$}
\label{stimatriviale}

\Dim
By induction, we get 
$$
\l|{G_n(x)\over[n]_q!}\r|_{v,Gauss}\leq{\sup\l(|G(x)|_{v,Gauss},1\r)^n\over
|[n]_q!|_v}\
$$
and the conclusion follows by lemma (\fstimafattoriale).
\QED

\Oss
We have assumed that $|1-q^{\kappa_v}|_v<|p|_v^{1/(p-1)}$. Let us briefly analyze
what happens if we drop this assumption. 
\par
We notice that if $|q|_v>1$ then $|[n]_q!|_v=|q|_v^{n(n-1)\over 2}$
for all $n\geq 0$, and therefore 
$$
\chi_v({\c M})^{-1}
=\limsup_{n\rightarrow\infty}\l|{G_n(x)\over[n]_q!}\r|_{v,Gauss}^{1/n}
=\limsup_{n\rightarrow\infty}{|G_n(x)|_{v,Gauss}^{1/n}\over 
|q|_v^{n-1\over 2}}\ .
$$
This limit can be zero as well as $\infty$ or a finite value. 
On the other hand, if $|q|_v<1$ then $|[n]_q!|=1$ for any $n\geq 0$,
and therefore 
$$
\chi_v({\c M})^{-1}
=\limsup_{n\rightarrow\infty}\l|{G_n(x)\over[n]_q!}\r|_{v,Gauss}^{1/n}
=\limsup_{n\rightarrow\infty}|G_n(x)|_{v,Gauss}^{1/n}
\leq \sup\l(|G(x)|_{v,Gauss},1\r)\ .
$$

\P
{If $|q|_v=1$ and $|1-q^{\kappa_v}|_v\geq |p|_v^{1/(p-1)}$ then 
$$
{\sup\l(|G(x)|_{v,Gauss},1\r)\over 
|p|_v^{1/{\kappa_v}(p-1)}|1-q^{e{\kappa_v}}|_v^{1/e{\kappa_v}}}
\geq \chi_v({\c M})^{-1}\geq {\sup\l(|G(x)|_{v,Gauss},1\r)\over 
|[{\kappa_v}]_q|_v^{1/{\kappa_v}}}\ , 
$$
where
$$
e=\inf\{m\in\interi: m>0, |1-q^{e{\kappa_v}}|_v<|p|_v^{1/(p-1)}\}\ .
$$}
\label{qmerdosi}

\Dim
Since $|1-q^{n\kappa_v}|_v\leq |1-q^{\kappa_v}|_v$ for any positive integer $n$ 
we have: $|[n]_q!|_v\leq |[{\kappa_v}]_q|_v^{\l[{n\over{\kappa_v}}\r]}$. 
For any positive integer $n$, there exist two positive integers $r,s<e$,
such that $n=se+r$. We obtain 
$$
|1-q^{n{\kappa_v}}|_v=|1-q^{(se+r){\kappa_v}}|_v
=|1-q^{se{\kappa_v}}+q^{se{\kappa_v}}(1-q^{r{\kappa_v}})|_v
\cases{
\geq |p|_v^{1/(p-1)} & if $r\neq 0$\cr
= |s|_v|1-q^{e{\kappa_v}}|_v & otherwise}\ ;
$$
from which we infer that 
$$\eqalign{|[{\kappa_v}]_q|_v^{1/{\kappa_v}}
&\geq \limsup_{n\rightarrow\infty}|[n]_q!|_v^{1/n}\cr
&\geq \liminf_{n\rightarrow\infty}|[n]_q!|_v^{1/n}\cr
&\geq \liminf_{n\rightarrow\infty}
   \Bigg(|p|_v^{\l(\l[n\over {\kappa_v}\r]-\l[n\over e{\kappa_v}\r]\r){1\over p-1}}
   |1-q^{e{\kappa_v}}|_v^{\l[n\over e{\kappa_v}\r]}
   \l|\l[n\over e{\kappa_v}\r]!\r|_v\Bigg)^{1/n}\cr
&\geq |p|_v^{1/{\kappa_v}(p-1)}|1-q^{e{\kappa_v}}|_v^{1/e{\kappa_v}}\ .}
$$
\par\vbox{\noindent
Finally we have 
$$
{\sup\l(|G(x)|_{v,Gauss},1\r)\over 
|p|_v^{1/{\kappa_v}(p-1)}|1-q^{e{\kappa_v}}|_v^{1/e{\kappa_v}}}
\geq \chi_v({\c M})^{-1}\geq {\sup\l(|G(x)|_{v,Gauss},1\r)\over 
|[{\kappa_v}]_q|_v^{1/{\kappa_v}}}\ .
$$
\QED}

%%%%%%%%%%%%%%%%%%%%%%%%%%%%%%%%%%%%%%%%%%%%%%%%%%%%%%%%%%%%%%%%%%%%%%%%%
%%%%%%%
\newsubchapter{{\titit q}-analogue of the Dwork-Frobenius theorem}
%%%%%%%
%%%%%%%%%%%%%%%%%%%%%%%%%%%%%%%%%%%%%%%%%%%%%%%%%%%%%%%%%%%%%%%%%%%%%%%%%

\vpar
The next proposition is the $q$-analogue of the Dwork-Frobenius-Young
theorem [DGS, VI, 2.1], which establishes a relation between $\chi_v({\c M})$ and
the coefficients of the $q$-difference matrix associated to
${\c M}=(M,\Phi_q)$ with respect to a cyclic basis, {\it when ${\kappa_v}=1$}: 

\P
{We suppose that $|1-q|_v<|p|_v^{1/(p-1)}$. Let $\c M$ be a $q$-difference module over $K_v(x)$ of rank
$\mu$, $m\in M$ a cyclic vector ({\it cf.}\ {\rm (\pvettorciclico)})
such that 
$$
\Dq(m,\Dq(m),\dots,\Dq^{\mu-1}(m))=
(m,\Dq(m),\dots,\Dq^{\mu-1}(m))
\pmatrix{\matrix{0&\dots& 0 }&\vrule&{a_0(x)}\cr
                     \hrulefill&&\hskip -25 pt\hrulefill\cr\cr
                     {\b I}_{\mu-1}&\vrule&
                          \matrix{a_1(x)\cr\vdots\cr
                            a_{\mu-1}(x)}}\ .
$$
If $\sup_{i=0,\dots,\mu-1}|a_i(x)|_{v,Gauss}>1$ then 
$$
\chi_v({\c M})={|p|_v^{1/(p-1)}\over 
\sup_{i=0,\dots,\mu-1}|a_i(x)|_{v,Gauss}^{1/(\mu-i)}}\ .
$$}
\label{dworkfrobenius}

\vpar
It follows immediately from (\pstimatriviale) that:

\C
{Let $|1-q|_v<|p|_v^{1/(p-1)}$. Then $\chi_v({\c M})\geq|p|_v^{1/(p-1)}$ if and only if
$$
\sup_{i=0,\dots,\mu-1}|a_i(x)|_{v,Gauss}\leq 1\ .
$$}
\label{cordworkfrobenius}

\vpar\noindent
{\bf Proof of proposition \pdworkfrobenius.}
We recall that $\chi_v({\c M})$ is independent of the choice of the basis 
of $M$ over $K_v(x)$.
\par
Let $\gamma\in K_v$ be such that 
$|\gamma|_v=\sup_{i=0,\dots,\mu-1}|a_i(x)|_{v,Gauss}^{1/(\mu-i)}$, 
$\u e=(m,\Dq(m),\dots,\Dq^{\mu-1}(m))$ and 
$$
H=\pmatrix{\gamma^{\mu-1} & & & 0\cr
            &\gamma^{\mu-2}&\cr
            &&\ddots\cr
           0&&& 1}\ .
$$
We set $\u f=\u e H$.  By a direct calculation we obtain 
$$
\Dq(\u f)=H^{-1}\Dq(\u e)H=\u f\ga W(x)\ ,
\hbox{\ with \ }
W(x)=\pmatrix{\matrix{0&\dots& 0 }&\vrule&{a_0(x)/\ga^\mu}\cr
                     \hrulefill&&\hskip -25 pt\hrulefill\cr\cr
                     {\b I}_{\mu-1}&\vrule&
                          \matrix{a_1(x)/\ga^{\mu-1}\cr\vdots\cr
                            a_{\mu-1}(x)/\ga}}\ .
$$
We set $\Dq^n(\u f)=\u fH_n(x)$, with $H_1(x)=\gamma W(x)$. 
We want to prove by induction on $n$ that 
$H_n(x)\equiv \ga^n W(x)\cdots W(q^{n-1}x)$ mod $\ga^{n-1}$. 
We remark that $H_n(x)\equiv \ga^n W(x)\cdots W(q^{n-1}x)$ mod
$\ga^{n-1}$ implies that $|H_n(x)|_{v,Gauss}\leq |\ga|_v^n$, and therefore 
that $|\dq H_n(x)|_{v,Gauss}\leq |\ga|_v^n$ ({\it cf.}\ (\pnormaderivazioni)). 
Then we have 
$$\eqalign{H_{n+1}(x)
&=H_1(x)H_n(qx)+\dq H_n(x)\cr
&\equiv\ga^{n+1} W(x)\cdots W(q^nx)\hbox{\ mod\ }\ga^n.}
$$
We deduce that $|H_n(x)|_{v,Gauss}\leq|\ga|_v^n$, 
for all $n\geq 1$, and hence that 
$$
\chi_v({\c M})\geq {|p|_v^{1/(p-1)}\over 
\sup_{i=0,\dots,\mu-1}|a_i(x)|_{v,Gauss}^{1/(\mu-i)}}\ .
$$
\par
Let us prove the reverse inequality.  By induction on $\mu$ one proves that 
the characteristic polynomial of $W(x)$ is 
$$
X^\mu-{a_{\mu-1}(x)\over\ga}X^{\mu-1}
-{a_{\mu-2}(x)\over\ga^2}X^{\mu-2}-\dots-{a_0(x)\over\ga^\mu}\ .
\numero{polcar}$$ 
By our choice of $\ga$, the reduction modulo $\varpi_v$ of (\fpolcar) has a nonzero 
root, hence $W(x)$ has an eigenvalue of norm $1$.  Then there exist 
\par\noindent
- an extension $L$ of $K_v(x)$ equipped with an extension 
of $|~|_{v,Gauss}$, still denoted $|~|_{v,Gauss}$, 
\par\noindent
- $\La\in L$, such that $|\La|_{v,Gauss}=1$,
\par\noindent
- $\ovec V\in L^\mu$, such that $|\ovec V|_{v,Gauss}=1$, 
\par\noindent
satisfying the relation  
$$
W(x)\cdots W(q^{n-1}x)\ovec V\equiv W(x)^n\ovec V\equiv\La^n\ovec V
\hbox{\ in the residue field of $L$ with respect to $|~|_{v,Gauss}$.}
$$
We deduce that  
$$
\ga^{-n}H_n(x)\ovec V\equiv \La^n\ovec V 
\hbox{\ in the residue field of $L$ with respect to $|~|_{v,Gauss}$.}
$$
Finally we obtain 
$$
\l|\ga^{-n}{H_n(x)\over[n]_q!}\r|_{v,Gauss}
\geq\l|{\ga^{-n}\over[n]_q!}H_n(x)\ovec V\r|_{v,Gauss}
=\l|{\La^n\ovec V\over[n]_q!}\r|_{v,Gauss}
=\l|{1\over[n]_q!}\r|_v 
$$
\par\vbox{\noindent
and hence 
$$\eqalign{\chi_v({\c M})^{-1}
&=\sup\l(1,\limsup_{n\rightarrow\infty}
   \l|{H_n(x)\over [n]_q!}\r|_{v,Gauss}^{1/n}\r)\cr
&\geq\limsup_{n\rightarrow\infty}\l|{\ga^n\over[n]_q!}\r|_v^{1/n}\cr
&={\sup_{i=0,\dots,\mu-1}|a_i(x)|_{v,Gauss}^{1/(\mu-i)}\over |p|_v^{1/(p-1)}}\ .}
$$
\QED}

\vpar
In the previous proposition we assumed that ${\kappa_v}=1$. 
If $\kappa_v>1$ we have:

\P
{The $q$-difference module $(M,\Phi_q)$ 
equipped with the operator $\Phi_q^{\kappa_v}$ (and consequently with 
$\Delta_{q^{\kappa_v}}=(\Phi_q^{\kappa_v}-{\b I}_\mu)/(q^{\kappa_v}-1)x$) is a
$q^{\kappa_v}$-difference module and 
$$
\chi_v(M,\Phi_q)\leq \chi_v(M,\Phi_q^{\kappa_v})^{1/{\kappa_v}}\ .
$$}
\label{unokappa}
 
\Dim
Applying successively (\ffiD) and (\fDfi), we obtain 
\smallskip
$$
\Delta_{q^{\kappa_v}}^n={(-1)^n\over 
(q^{\kappa_v}-1)^nx^n}\sum_{i=0,\dots,n\atop j=0,\dots,i{\kappa_v}}
\l((-1)^i{n\choose i}_{q^{-{\kappa_v}}}q^{-{\kappa_v}{i(i-1)\over 2}}{i{\kappa_v} 
\choose j}_q(q-1)^jq^{j(j-1)\over 2}x^j\r)
\Dq^j\ .
\numero{schifo}$$
Let $\u f$ be a basis of $M$ over $K_v(x)$ such that 
$\Delta_{q^{\kappa_v}}^n\u f=\u f H_n(x)$ and 
$\Dq^n\u f=\u fG_n(x)$. 
We deduce using (\fschifo) that 
$$
|H_n(x)|_{v,Gauss}\leq {1\over |q^{\kappa_v}-1|_v^n}
\l(\sup_{s\leq n{\kappa_v}}|G_s(x)|_{v,Gauss}\r)\ .
$$
\par\vbox{\noindent
Recalling the estimates in (\pext) and some general properties of $\limsup$ 
({\it cf.}\ [AB, II, 1.8]) we obtain 
$$\eqalign{{1\over \chi_v(M,\Phi_{q^{\kappa_v}})}
&=\limsup_{n\rightarrow\infty}\l|{H_n(x)\over n_{q^{\kappa_v}}!}\r|_{v,Gauss}^{1/n}\cr
&\leq {\limsup_{n\rightarrow\infty}\l(\sup_{s\leq n{\kappa_v}}|G_s(x)|_{v,Gauss}\r)^{1/n}
    \over|q^{\kappa_v}-1|_v|p|_v^{1/(p-1)}}\cr
&= {1\over \chi_v(M,\Phi_q)^{\kappa_v}}\ .}
$$
\QED}

%%%%%%%%%%%%%%%%%%%%%%%%%%%%%%%%%%%%%%%%%%%%%%%%%%%%%%%%%%%%%%%%%%%%%%%%%
%%%%%%%%%%%%%%%%%%%%%%%%%%%%  5  %%%%%%%%%%%%%%%%%%%%%%%%%%%%%%%%%%%%%%%%
%%%%%%%%%%%%%%%%%%%%%%%%%%%%%%%%%%%%%%%%%%%%%%%%%%%%%%%%%%%%%%%%%%%%%%%%%

\newchapter{{\tiit p}-adic criteria for unipotent reduction}
\chaplabel{nilp}

%%%%%%%%%%%%%%%%%%%%%%%%%%%%%%%%%%%%%%%%%%%%%%%%%%%%%%%%%%%%%%%%%%%%%%%%%

\vpar
We recall that $K_v$ is a complete field with respect to the norm 
$|~|_v$ and that ${\c V}_v$ is its ring of integers, $\varpi_v$ its 
uniformizer and $k_v$ the residue field. 
\par
{\it Let $q$ be an element 
of $K_v$ such that $|q|_v=1$, $q$ is not a root of unity}, and the order 
${\kappa_v}$ of its image in the multiplicative group $k_v^\times$ is finite. 

\vpar
Let ${\c F}\subset K_v(x)$ be a $q$-difference algebra essentially of finite type 
over ${\c V}_v$ ({\it cf.}\ (\pdilatationdefinizione)). 
Let $\frak a$ be an ideal of ${\c V}_v$ and 
$\o q$ be the image of $q$ in ${\c V}_v/{\frak a}$.  
The algebra ${\c F}\otimes_{{\c V}_v}{\c V}_v/{\frak a}$ has the natural
structure of a $\o q$-difference algebra. 
Let ${\c M}=(M,\Phi_q)$ be a $q$-difference module over $\c F$.
We consider the free  
${\c F}\otimes_{{\c V}_v}{\c V}_v/{\frak a}$-module 
$M\otimes_{{\c V}_v}{\c V}_v/{\frak a}$ equipped with the morphism
$\Phi_{\o q}$ induced by $\Phi_q$: it is a $\o q$-difference module over 
${\c F}\otimes_{{\c V}_v}{\c V}_v/{\frak a}$, 
which satisfies the assumptions of \ccurva. 
\par
We are especially interested in the following two cases:
\smallskip\noindent
- ${\frak a}$ is the maximal ideal of ${\c V}_v$ generated by $\varpi_v$.  
We will refer to 
$(M\otimes_{{\c V}_v}{\c V}_v/\varpi_v{\c V}_v, \Phi_{\o q})$ 
as the {\it reduction of $\c M$ modulo $\varpi_v$ or over $k$}.
\smallskip\noindent
- $\frak a$ is the ideal of ${\c V}_v$ generated by $1-q^{{\kappa_v}}$. 
We will refer to 
$(M\otimes_{{\c V}_v}{\c V}_v/(1-q^{{\kappa_v}}){\c V}_v, \Phi_{\o q})$ 
as the {\it reduction modulo $1-q^{{\kappa_v}}$}.

\Oss
We notice that $|p!|_v=|p|_v$, therefore 
both reductions are $q$-analogues of the reduction modulo $p$ in the
differential case ({\it cf.}\ \cdifferenziale). 
In (\regolarita) we analyze the reduction modulo $\varpi_v$, while in
our main theorem (\grothendieck) we consider the reduction modulo $1-q^{{\kappa_v}}$.

\vpar
Motivated by \ccurva, we are particularly interested in $q$-difference modules 
$\c M$ over $\c F$ such that the reduction modulo 
$\varpi_v$ (resp. $1-q^{\kappa_v}$) of the operator 
$\Phi_q^{\kappa_v}$ is unipotent. 
We shall say briefly that $\c M$ has {\it unipotent reduction of
order $n$ modulo $\varpi_v$ (resp. $1-q^{\kappa_v}$)} if the 
reduction of $\Phi_q^{\kappa_v}$ modulo  $\varpi_v$ (resp. $1-q^{\kappa_v}$)
is a unipotent morphism of order $n$.

\vpar
The following example shows that 
$\c M$ can have unipotent reduction modulo $\varpi_v$  
without having unipotent reduction modulo $1-q^{\kappa_v}$:

\vpar\noindent
{\bf Example.} 
Let us consider the $q$-difference module over $\raz_p(x)$ associated to the 
$q$-difference system 
$$
\pmatrix{y_1(qx)\cr y_2(qx)}=\pmatrix{1&p\cr 0&1}\pmatrix{y_1(x)\cr y_2(x)}\ .
\numero{esempio}$$
Then 
$$
\pmatrix{y_1(q^{\kappa_p} x)\cr y_2(q^{\kappa_p} x)}
=\pmatrix{1&{\kappa_p} p\cr 0&1}\pmatrix{y_1(x)\cr y_2(x)}\ ,
$$
from which it follows that
$$
\l|\pmatrix{1&{\kappa_p} p\cr 0&1}-{\b I}_2\r|_p=|{\kappa_p} p|_p=|p|_p\ .
$$
If we choose $q=8$ and $p=3$ then ${\kappa_p}=2$ and
$|1-q^{\kappa_p}|_p=|1-8^2|_p=|3^2|_p<|3|_p$, and therefore 
$\pmatrix{1&{\kappa_p} p\cr 0&1}=\pmatrix{1& 6\cr 0&1}\equiv {\b I}_2$ mod
$3$, but $\pmatrix{1& 6\cr 0&1}\not\equiv {\b I}_2$ mod $3^2$. 
Then by (\pgrothradici) the $q$-difference module associated to (\fesempio) has trivial
reduction modulo $p=3$, but not modulo $1-q^{\kappa_p}=(-7)3^2$.

\vpar
We want to relate the property of having unipotent reduction 
modulo $\varpi_v$ (resp. $1-q^{\kappa_v}$) to 
an estimate of the invariant $\chi_v({\c M}):=\chi_v({\c M}_{K_v(x)})$. 

%%%%%%%%%%%%%%%%%%%%%%%%%%%%%%%%%%%%%%%%%%%%%%%%%%%%%%%%%%%%%%%%%%%%%%%%%
%%%%%%%
\newsubchapter{{\titit q}-difference modules having 
unipotent reduction modulo $\varpi_v$}
%%%%%%%
%%%%%%%%%%%%%%%%%%%%%%%%%%%%%%%%%%%%%%%%%%%%%%%%%%%%%%%%%%%%%%%%%%%%%%%%%

\vpar
First we consider $q$-difference modules having unipotent reduction modulo $\varpi_v$. 
The following proposition is a $q$-analogue of a classical estimate for 
$p$-adic differential modules [DGS, page 96]:

\P
{If $\c M$ has unipotent reduction modulo $\varpi_v$ of order $n$ then 
$$
\chi_v({\c M})\geq
|\varpi_v|_v^{-1/{\kappa_v} n}|[{\kappa_v}]_q|_v^{1/{\kappa_v}}
       |p|_v^{1/{\kappa_v}(p-1)}\ .
\numero{chizero}$$}
\label{chinilp}

\vpar
The proof of (\pchinilp) relies on the following lemma:

\L
{Let us assume that $\c M$ has unipotent 
reduction of order $n$ modulo $\varpi_v$. 
Let $\u e$ be a basis of $M$ over ${\c F}$ and let $\Dq^m\u e=\u eG_m(x)$, 
for any $m\geq 1$, with $G_m(x)\in M_{\mu\times\mu}(K_v(x))$. 
Then 
$$
|G_{sn{\kappa_v}}(x)|_{v,Gauss}\leq
|\varpi_v|_v^s\ ,
$$
for every integer $s\geq 1$.}
\label{lemmanilpuno}

\Dim
By (\fgrothleibniz), for all $s\in\nat$, $s>1$ we have 
$$\eqalign{\Dq^{(s+1)n{\kappa_v}}(\u e)
&=\u e G_{(s+1)n{\kappa_v}}(x)\cr
&=\Dq^{n{\kappa_v}}(\u e G_{sn{\kappa_v}}(x))\cr
&=\sum_{i=0}^{n{\kappa_v}}{n{\kappa_v}\choose i}_q
   \Dq^{n{\kappa_v}-i}(\u e)\dq^i
   \l(G_{sn{\kappa_v}}\r)(q^{n{\kappa_v}-i}x)\cr
&=\u e\sum_{i=0}^{n{\kappa_v}}{n{\kappa_v}\choose i}_q
   G_{n{\kappa_v}-i}(x)\dq^{n{\kappa_v}-i}
   \l(G_{sn{\kappa_v}}\r)(q^{n{\kappa_v}-i}x)\ ,}
$$
and hence 
$$
G_{(s+1)n{\kappa_v}}(x)=
\sum_{i=0}^{n{\kappa_v}}
{n{\kappa_v}\choose i}_qG_{n{\kappa_v}-i}(x)\dq^i
\l(G_{sn{\kappa_v}}\r)(q^{n{\kappa_v}-i}x)
\numero{coeffnilp}$$
By (\punipotenza), the definition of unipotent reduction modulo
$\varpi_v$ is equivalent to the condition  
$$
|G_{n{\kappa_v}}(x)|_{v,Gauss}\leq|\varpi_v|_v\ . 
$$
We shall prove the statement by induction on $s>1$, using (\fcoeffnilp).
We suppose that 
$$
|G_{sn{\kappa_v}}(x)|_{v,Gauss}\leq|\varpi_v|_v^s\ .
$$ 
Then 
all the terms occurring in the sum (\fcoeffnilp) are bounded by 
$|\varpi_v|_v^{s+1}$, in fact:
\par\noindent 
1) If $({\kappa_v},i)=1$, then 
$$
\l|{n{\kappa_v}\choose i}_q\r|_v=\l|{[n{\kappa_v}]_q\over \iq}{n\kappa_v-1\choose i-1}_q\r|_v
\leq|[{\kappa_v}]_q|_v\leq |\varpi_v|_v
$$ 
and the absolute value of the corresponding term in sum (\fcoeffnilp) is bounded by 
$|\varpi_v|_v^{s+1}$.
\par\noindent
2) For all $i=1,\dots,n$, we have $|\dq^{i{\kappa_v} }f(x)|_{v,Gauss}\leq |[{\kappa_v}]_q|_v|f(x)|_{v,Gauss}$ 
and therefore 
$$
\l|{n{\kappa_v} \choose i{\kappa_v} }_q
G_{n{\kappa_v} -i{\kappa_v} }(x)\dq^{i{\kappa_v} }
\l(G_{sn{\kappa_v} }\r)(q^{n{\kappa_v} -i{\kappa_v} }x)
\r|_{v,Gauss}\leq  |[{\kappa_v}]_q|_v|\varpi_v|_v^s\leq |\varpi_v|_v^{s+1}\ ,
$$
for all $i=1,\dots,n$. 
\par\noindent
3) The term of (\fcoeffnilp) corresponding to $i=0$ is 
$G_{n{\kappa_v} }(x)G_{sn{\kappa_v} }(q^{n{\kappa_v} }x)$, and
therefore it is bounded by $|\varpi_v|_v^{s+1}$, by induction.
\par
Thus we have proved that $|G_{(s+1)n{\kappa_v} }(x)|_{v,Gauss}\leq |\varpi_v|_v^{s+1}$.
\QED

\vpar\noindent
{\bf Proof of proposition \pchinilp.}
By the recursive formula (\fcoefficienti) we have 
$$
\l|G_m(x)\r|_{v,Gauss}\leq 
\l|G_{\l[{m\over n{\kappa_v} }\r]n{\kappa_v} }(x)\r|_{v,Gauss}\ .
$$
\par\vbox{\noindent
The estimate (\fchizero) follows from previous lemma, since 
$$\eqalign{\chi_v({\c M})
&\geq \inf\l(1,\liminf_{m\rightarrow\infty}
{|G_{\l[{m\over n{\kappa_v}}\r]n{\kappa_v}}(x)|_{v,Gauss}^{-1/m}
   \over |[m]_q!|_v^{-1/m}}\r)\cr
&\geq \liminf_{m\rightarrow\infty}
   |\varpi_v|_v^{-\l[{m\over n{\kappa_v}}\r]{1\over m}}|m_q!|_v^{1/m}\cr
&=|\varpi_v|_v^{-1/n{\kappa_v}}|[{\kappa_v}]_q|_v^{1/{\kappa_v}}|p|_v^{1/{\kappa_v}(p-1)}\ .}
$$
\QED}

\C
{The $q$-difference module $\c M$ has unipotent 
reduction modulo $\varpi_v$ if and only if 
$$
\chi_v({\c M})>|[{\kappa_v}]_q|_v^{1/{\kappa_v}}|p|_v^{1/{\kappa_v}(p-1)}\ .
$$}
\label{corchinilp}

\Dim
If $\c M$ has unipotent reduction modulo $\varpi_v$, we immediately deduce by
(\pchinilp) that $\chi_v({\c M})>|[{\kappa_v}]_q|_v^{1/{\kappa_v}}|p|_v^{1/{\kappa_v}(p-1)}$.
\par
On the other hand, by hypothesis we have 
$$\eqalign{|[{\kappa_v}]_q|_v^{1/{\kappa_v}}|p|_v^{1/{\kappa_v}(p-1)}<\chi_v({\c M})
&=\inf\l(1,\liminf_{n\rightarrow\infty}
   \l|{G_n(x)\over [n]_q!}\r|_{v,Gauss}^{-1/n}\r)\cr
&=\inf\l(1,|[{\kappa_v}]_q|_v^{1/{\kappa_v}}
  |p|_v^{1/{\kappa_v}(p-1)}\liminf_{n\rightarrow\infty}
   \l|G_n(x)\r|_{v,Gauss}^{-1/n}\r)\ .}
$$
We deduce that 
$$
\limsup_{n\rightarrow\infty}
   \l|G_n(x)\r|_{v,Gauss}^{1/n}<1\ .
$$
We conclude that there exists $N\in\nat$ 
such that $|G_n(x)|_{v,Gauss}<1$ for all $n>N$, 
which implies that $\c M$ has unipotent reduction modulo $\varpi_v$.
\QED

%%%%%%%%%%%%%%%%%%%%%%%%%%%%%%%%%%%%%%%%%%%%%%%%%%%%%%%%%%%%%%%%%%%%%%%%%
%%%%%%%
\newsubchapter{{\titit q}-difference modules having 
unipotent reduction modulo $1-q^{\kappa_v}$}
%%%%%%%
%%%%%%%%%%%%%%%%%%%%%%%%%%%%%%%%%%%%%%%%%%%%%%%%%%%%%%%%%%%%%%%%%%%%%%%%%

\vpar
Under the hypothesis of unipotent reduction modulo $1-q^{\kappa_v}$, we obtain a 
slight but crucial improvement of (\pchinilp) 
that will be fundamental in the proof of the $q$-analogoue of
Grothendieck's conjecture (\grothendieck) below:

\P
{Let $\c M$ be a $q$-difference module 
over $\cal F$, with unipotent reduction modulo
$1-q^{\kappa_v}$ of order $n$. Then  
$$
\chi_v({\c M})\geq|[{\kappa_v}]_q|_v^{(n-1)/n{\kappa_v}}|p|_v^{1/{\kappa_v}(p-1)}\ .
$$}

\Dim
Let $\u e$ be the basis of $M$ such that 
$\Dq^m\u e=\u eG_m(x)$, for all $m\geq 1$. Then 
$$
|G_{n{\kappa_v}}(x)|_{v,Gauss}\leq |[{\kappa_v}]_q|_v\ .
$$
\par
The estimates in (\plemmanilpuno) show that 
$$
|G_{sn{\kappa_v}}(x)|_{v,Gauss}\leq |[{\kappa_v}]_q|_v^s\ ,\ \forall\ s\geq 1, 
\numero{trivG}$$
\par\vbox{\noindent
therefore we conclude that 
$$
\chi_v({\c M})\geq|[{\kappa_v}]_q|_v^{(n-1)/n{\kappa_v}}|p|_v^{1/{\kappa_v}(p-1)}\ .
$$
\QED}

\C
{The following assertions are equivalent: 
\tparnoind
1) $\chi_v({\c M})\geq |p|_v^{1/{\kappa_v}(p-1)}$.
\tparnoind
2) There exists a cyclic basis $\u e$ of $({\cal M}_{K(x)},\phi_q^{\kappa_v})$ such that 
$\Phi_q^{\kappa_v}\u e\equiv \u e$ modulo $1-q^{\kappa_v}$.}
\label{chitrivial}

\Dim
The implication ``{\sl 2)$\Rightarrow$1)}'' is a consequence of the previous proposition. 
\par
We prove ``{\sl 1)$\Rightarrow$2)}''.
The $\c F$-module $M$ 
equipped with the operators $\Phi_q^{\kappa_v}$ is a $q^{\kappa_v}$-difference module. 
It follows by (\punokappa) that $\chi_v(M,\Phi_q^{\kappa_v})\geq
|p|_v^{1/(p-1)}$. 
We know by (\pvettorciclico) that ${\cal M}_{K(x)}$ 
admits a cyclic basis $\u e$ over $K(x)$. Let $\Phi_q^{\kappa_v}\u e=\u eA_{\kappa_v}(x)$. 
We deduce from (\pdworkfrobenius) that 
$$
\l|{A_{\kappa_v}(x)-{\b I}_\mu\over (q^{\kappa_v}-1)x}\r|_{v,Gauss}\leq 1 
$$ 
and hence that $|A_{\kappa_v}(x)-{\b I}_\mu|_{v,Gauss}\leq|1-q^{\kappa_v}|_v$.
\QED

%%% Local Variables: 
%%% mode: plain-tex
%%% TeX-master: "q-groth"
%%% End: q-groth"
%%% End: 

%%%%%%%%%%%%%%%%%%%%%%%%%%%%%%%%%%%%%%%%%%%%%%%%%%%%%%%%%%%%%%%%%%%%%%%%%
%%%%%%%%%%%%%%%%%%%%%%%%% PART III    %%%%%%%%%%%%%%%%%%%%%%%%%%%%%%%%%%%
%%%%%%%%%%%%%%%%%%%%%%%%%%%%%%%%%%%%%%%%%%%%%%%%%%%%%%%%%%%%%%%%%%%%%%%%%
%%%%%%%%%%%%%%%%%%%%%%%%%%%%%%%%%%%%%%%%%%%%%%%%%%%%%%%%%%%%%%%%%%%%%%%%%
\vskip 20pt 
\centerline{\titlepart Part III.  A {\titlepartit q}-analogue of Grothendieck's conjecture} 
\vpar
\centerline{\titlepart on {\titlepartit p}-curvatures}
%%%%%%%%%%%%%%%%%%%%%%%%%%%%%%%%%%%%%%%%%%%%%%%%%%%%%%%%%%%%%%%%%%%%%%%%%
%%%%%%%%%%%%%%%%%%%%%%%%%%%%%%%%%%%%%%%%%%%%%%%%%%%%%%%%%%%%%%%%%%%%%%%%%
%%%%%%%%%%%%%%%%%%%%%%%%%%%%%%%%%%%%%%%%%%%%%%%%%%%%%%%%%%%%%%%%%%%%%%%%%
%%%%%%%%%%%%%%%%%%%%%%%%%%%%%%%%%%%%%%%%%%%%%%%%%%%%%%%%%%%%%%%%%%%%%%%%%

%%%%%%%%%%%%%%%%%%%%%%%%%%%%%%%%%%%%%%%%%%%%%%%%%%%%%%%%%%%%%%%%%%%%%%%%%
%%%%%%%%%%%%%%%%%%%%%%%%%%%%   6  %%%%%%%%%%%%%%%%%%%%%%%%%%%%%%%%%%%%%%%
%%%%%%%%%%%%%%%%%%%%%%%%%%%%%%%%%%%%%%%%%%%%%%%%%%%%%%%%%%%%%%%%%%%%%%%%%

\newchapter{Arithmetic {\tiit q}-difference modules and regularity}
\chaplabel{regolarita}

%%%%%%%%%%%%%%%%%%%%%%%%%%%%%%%%%%%%%%%%%%%%%%%%%%%%%%%%%%%%%%%%%%%%%%%%%

\vpar
We now establish some notation that will be maintained until the end of the paper:
\par\noindent
$K$= a number field.
\par\noindent
${\c V}_K$= the ring of integers of $K$.
\par\noindent
$|~|_v$= a $v$-adic absolute value of $K$. 
In the non-archimedean case we normalize $|~|_v$ as follows:
$$
|p|_v=p^{-[K_v:{\b Q}_p]/[K:{\b Q}]}\;,
$$
where $K_v$ is the $v$-adic completion of $K$ and $v\vert p$. 
Similarly, in the archimedean case we normalize $|~|_v$ by setting 
$$|x|_v=\cases{\cr
|x|_\reali^{1/[K:{\b Q}]} & if $K_v=\reali$\cr\cr
|x|_{\b C}^{2/[K:{\b Q}]} & if $K_v={\b C}$\cr\cr}\ ,$$
where $|~|_\reali$ and $|~|_{\b C}$ are the usual absolute values 
of $\reali$ and and of ${\b C}$ respectively.
\par\noindent
$\Sg_f$= the set of finite places $v$ of $K$. 
\par\noindent
$\varpi_v$= uniformizer $\in{\c V}_K$ associated to the finite place $v$. 
\par\noindent
$k_v$= residue field of $K$ with respect to a finite place $v$. 
\par\noindent
$\Sg_\infty$=the set of archimedean places of $K$.

%%%%%%%%%%%%%%%%%%%%%%%%%%%%%%%%%%%%%%%%%%%%%%%%%%%%%%%%%%%%%%%%%%%%%%%%%
%%%%%%%
\newsubchapter{On cyclic subgroups of $\o\raz^\times$}
\smallskip
\centerline{\tit and their reduction modulo almost every prime}

%%%%%%%
%%%%%%%%%%%%%%%%%%%%%%%%%%%%%%%%%%%%%%%%%%%%%%%%%%%%%%%%%%%%%%%%%%%%%%%%%

\vpar
{\it We fix an element $q$ of $K$ which is not zero and not   
a root of unity.} 
For each $v\in\Sg_f$ such that $|q|_v=1$, we define $\kappa_{v}$ to be 
the multiplicative order of the image of $q$ in the residue field of $K$ 
with respect to $v$. 
We refer to [BHV] for the most recent results on the distribution of 
$(\kappa_v)_v$. 

\vpar
We recall that the {\it Dirichlet density $d(S)$ of a set $S$} of finite places of a
number field $K$({\it cf.}\ for instance [N, VII, \S 13]) is defined by  
$$
d(S)=\limsup_{s\rightarrow 1^+}
{\sum_{v\in S}p^{-sf_v}\over \sum_{v\in\Sg_f}p^{-sf_v}}\ ,
$$
where $f_v=[k_v:{\b F}_p]$, if $v\vert p$. 

\vpar
The proposition below is a particular case of a theorem by Schinzel
[Sc, Th. 2]. We prefer to give a direct proof here.

\P
{Let $S\subset\Sg_f$ be a set of finite places of $K$ of Dirichlet density 1 
and let $a,b$ be two elements of $K^\times=K\smallsetminus\{0\}$ 
such that for all $v\in S$, the reduction of $b$ modulo $\varpi_v$ 
belongs to the cyclic group generated by the reduction of $a$ modulo $\varpi_v$. 
Then $b\in a^\interi$.}
\label{colmez}

\C
{Let $a,b$ be two elements of $K$, which are not roots of unity,  
such that for almost all $v\in\Sg_f$
the order of $a$ modulo $\varpi_v$ and the order of $b$ modulo $\varpi_v$ coincide. 
Then either $a=b$ or $a=b^{-1}$.}
\label{corollariouffa}

\Oss
This shows that $\{q,q^{-1}\}$ is uniquely determined by the
family of integers $(\kappa_v)_v$.

\vpar\noindent
{\bf Proof of corollary \pcorollariouffa.}
We recall that $k_v^\times$ is a cyclic group and that, therefore, its
subgroups are determined by their order.
By the previous proposition, we know that $b=a^n$ and $a=b^m$ for some integers
$n$ and $m$. We deduce that $a^{nm}=a$. Since $a$ is not a root of
unity, we have $mn=1$ and hence either $m=n=1$ or $m=n=-1$.
\QED

\vpar\noindent
{\bf Proof of proposition \pcolmez\ (following an argument of P. Colmez).} 
We fix a rational prime $\ell$. Let $\zeta_{\ell}$ be an $\ell$-th root
of unity. 
We consider the following Galois extensions of $K$: 
$K_1=K(a^{1/\ell},\zeta_{\ell})$, $K_2=K(b^{1/\ell},\zeta_{\ell})$ 
and $K_{12}=K(a^{1/\ell},b^{1/\ell},\zeta_{\ell})$. 
We will prove that $K_1=K_{12}$, and hence that $K_2\subset K_1$, 
by applying the following corollary of the  \v Cebotarev Density theorem:

\vpar
\centerline{\vbox{\hsize 15 true cm\noindent
{\bf [N, VII, (13.6)]} 
{\sl Let $\widetilde K$ be a Galois extension of the number field $K$ 
and let $P(\widetilde K/K)$ be the set of primes of $K$ 
that split totally in $\widetilde K$. Then the Dirichlet density of 
$P(\widetilde K/K)$ is 
$$
d(P(\widetilde K/K))={1\over [\widetilde K:K]}\ .
$$}}}

\vpar\noindent
Let $v\in\Sg_f$ be a prime of $K$ such that $v\vert p$, $p>\ell$, 
and let $\{w_1,\dots,w_r\}\subset\Sg_f$ be 
the set of all primes $w$ of $K_1$ such that $w\vert v$. 
Let $e_i$ be the ramification index of $w_i\vert v$ and 
$f_i$ be the residue degree.
Since $K_1/K$ is a Galois extension we have: $e=e_1=\dots=e_r$ and 
$f=f_1=\dots=f_r$ ({\it cf.}\ [N, IV, page 55]).  
Therefore,  $e=f=1$ if and only if we have 
$[K_1:K]=\sum_{i=1}^r e_i f_i=r$: so $v$ splits totally in $K_1$ if 
and only if $e=f=1$.  
Then $P(K_1/K)$ is the set of all primes $v\in\Sg_f$ of $K$ 
such that:
$$\eqalign{
&\hbox{- $v\vert p$ and $p\equiv 1$ mod $\ell$;}\cr 
&\hbox{- there exists $a^\p\in k_v$ such that 
  ${a^\p}^{\ell}\equiv a$ in $k_v$;}\cr
&\hbox{- $v$ is not ramified in $K_1$.}\cr}
$$
For the same reason, $P(K_{12}/K)$ is the set 
of all primes $v\in\Sg_f$ such that:
$$\eqalign{
&\hbox{- $v\vert p$ and $p\equiv 1$ mod $\ell$;}\cr 
&\hbox{- there exists $a^\p\in k_v$ such that 
  ${a^\p}^{\ell}\equiv a$ in $k_v$;}\cr
&\hbox{- there exists $b^\p\in k_v$ such that 
  ${b^\p}^{\ell}\equiv b$ in $k_v$;}\cr
&\hbox{- $v$ is not ramified in $K_{12}$.}\cr}
$$
Let $v\in P(K_1/K)\cap S$ and let $a^\p\in k_v$ be such that 
$a\equiv {a^\p}^{\ell}$ in $k_v$. By hypothesis 
there exists a positive integer $n(v)$ such that 
$b\equiv a^{n(v)}$ in $k_v$ and hence 
$b\equiv ({a^\p}^{n(v)})^{\ell}$. 
Hence if $vv\in P(K_1/K)\cap S$ is not ramified in $K_{12}$ then $v\in P(K_{12}/K)$.
Taking into account that $S$ has density 1 and that there are only finitely many $v\in\Sg_f$ which ramify 
in $K_{12}$, we have 
$$
d(P(K_{12}/K))={1\over [K_{12}:K]}\geq 
d(P(K_1/K))={1\over [K_1:K]}\ .
$$ 
We conclude that $K_{12}=K_1$ and therefore $K_2\subset K_1$.
\par
We recall the following fact from Kummer theory:

\vpar
\centerline{\vbox{\hsize 15 true cm\noindent
{\bf [N, VII, (3.6)]} 
{\sl Let $n$ be a positive integer which is relatively prime with respect to the
characteristic of the field $K$, and assume that $K$ contains the group of
$n$-th roots of unity. Then the abelian extensions $\widetilde K/K$
of exponents $n$ are in one-to-one correspondence with the subgroups 
$\Gamma\subset K^\times=K\smallsetminus \{0\}$, which contain 
${K^\times}^n$, via the rule 
$\Gamma\longmapsto \widetilde K
=K(\vbox to 8pt{}^n\!\!\!\!\sqrt{\Gamma})$.}}}

\vpar\noindent
This statement, applied to $K_2\subset K_1=K_{12}$ and $n=\ell$,
says that 
$$
b^\interi{K^\times}^{\ell}\subset a^\interi{K^\times}^{\ell}\Rightarrow
b\in a^\interi{K^\times}^{\ell}\ .
$$
Since $\ell$ is arbitrary, we conclude that $b\in a^\interi$.\QED

%%%%%%%%%%%%%%%%%%%%%%%%%%%%%%%%%%%%%%%%%%%%%%%%%%%%%%%%%%%%%%%%%%%%%%%%%
%%%%%%%
\newsubchapter{Unipotent reduction and regularity}
%%%%%%%
%%%%%%%%%%%%%%%%%%%%%%%%%%%%%%%%%%%%%%%%%%%%%%%%%%%%%%%%%%%%%%%%%%%%%%%%%

\vpar
Let ${\c M}=(M,\Phi_q)$ be a $q$-difference module over a
$q$-difference algebra ${\c F}\subset K(x)$ 
essentially of finite type over ${\c V}_K$ ({\it cf.}\ (\pdilatationdefinizione)). 

\Def
{The $q$-difference module $\c M$ over $\c F$ is {\it regular singular} 
if both ${\c M}_{K((x))}$ and ${\c M}_{K((1/x))}$ are 
regular singular $q$-difference modules.}

\vpar
Let $\Sg_{nilp}$ be the set of finite places $v$ of $K$ 
such that ${\c M}$ has unipotent reduction modulo $\varpi_v$. 
The following result is a $q$-analogue of a well known result due to 
Katz ({\it cf.}\ [K1, 13.0]). 

\T
{\tparnoind
1) If $\Sg_{nilp}$ is infinite, then ${\c M}$ is regular singular.
\tparnoind
2) If moreover $\Sg_{nilp}$ has Dirichlet density 1, 
the exponents of $K((x))\otimes_{\c F}M$ with respect to some basis $\u e$ 
(and hence to any basis) over $K((x))$ coincide with $q^\interi$.}
\label{IIIkatz}

\Dim
\par\noindent
{\sl 1)} It is enough to prove the statement at zero. 
Let $\u e$ be a basis of $M$ over ${\c F}$ such that 
$\Phi_q\u e=\u eA(x)$, with $A(x)\in Gl_\mu({\c F})$.
Then $A(x)$ can be regarded as an element of $Gl_\mu\l(K((x))\r)$, 
which means that $A(x)$ has the following form:
$$
A(x)={1\over x^l}\sum_{i\geq 0}A_ix^i\in 
{1\over x^l}Gl_\mu(K\[[x\]])\ ,
$$
for some $l\in\interi$.
If $l=0$, the $q$-difference module $\c M$ is regular singular at 
zero, so let us suppose $l\neq 0$.
For all positive integers $m$, we have 
$$
\Phi_q^m(\u e)
=\u eA(x)A(qx)\cdots A(q^{m-1}x)
=\u e\l({A_0^m\over q^{lm(m-1)\over 2}x^{ml}}+h.o.t.\r)
$$
By hypothesis, for any $v\in\Sg_{nilp}$, there exists a positive integer 
$n(v)\geq 1$  such that we have 
$$
\l(A(x)A(qx)\cdots A(q^{\kappa_v-1}x)
-1\r)^{n(v)}\equiv 0\hbox{\ mod $\varpi_v$}; 
$$
we deduce that $A_0^{\kappa_v}\equiv 0$ modulo $\varpi_v$,
for any $v\in\Sg_{nilp}$, and hence that 
$A_0$ is a nilpotent matrix. 
\par
We suppose that zero is not a regular singularity.
By (\pformclass), there exist an extension $L((t))$ of $K((x))$ 
and $\widetilde q\in L$, with $t^d=x$ and $\widetilde q^d=q$,  
such that we can find a basis $\u f$ of $L((t))\otimes_{\c F}M$ over 
$L((t))$ with the following properties:
$$
\Phi_{\widetilde q }(\u f)=\u fB(t)
$$ 
and 
$$
B(t)={B_k\over t^k}+{B_{k-1}\over t^{k-1}}+\cdots+{B_1\over t^1}
+\widetilde B_0(t)\ ,
$$
with $\widetilde B_0(t)\in M_{\mu\times\mu}(L\[[t\]])$, 
$k\geq 1$ and $B_k\in Gl_\mu(L)$ {\it non} nilpotent and in Jordan normal
form. 
Let $F(t)={F\over t^m}+h.o.t\in Gl_\mu(L((t)))$ be such that 
$\u e=\u f F(t)$. This implies that 
$A(x)=F(t)^{-1}B(t)F(\widetilde q t)$.  
We get a contradiction since the matrix $F\in Gl_\mu(L)$ satisfies $A_0=F^{-1}B_k F$. 
\par\noindent
{\sl 2)} We know by {\sl 1)} that $\c M$ has a regular 
singularity at zero.  
Then there exists a $K((x))$-basis $\u e$ of $K((x))\otimes_{{\c F}}M$ 
such that $\Phi_q(\u e)=\u eA$, with $A\in Gl_\mu(K)$ 
in Jordan normal form.
By remark (\pAo), we can chose $\u e$ such that 
for almost all $v\in\Sg_{nilp}$ there exists $n(v)\geq 1$ 
satisfying the equivalence 
$$
(A^{\kappa_v }-1)^{n(v)}\equiv 0\hbox{\ mod $\varpi_v$}.  
$$
Therefore the matrix $A^{\kappa_v }$ is unipotent modulo $\varpi_v$ for
all $v\in\Sg_{nilp}$.  We deduce that the reduction modulo $\varpi_v$ of the 
eigenvalues of $A$ are $\kappa_v$-th roots of unity for almost all $v\in\Sg_{nilp}$.
This means that the reduction modulo $\varpi_v$ of any eigenvalues of $A$ 
is an element of the cyclic group generated by the reduction of $q$, 
for  all $v\in\Sg_{nilp}$.  The conclusion follows by applying lemma (\pcolmez).
\QED

\P
{Let us assume that the $q$-difference module $\c M$ over $\c F$ 
has the property that for almost all finite places $v$ of $K$ 
the morphism $\Phi_q^{\kappa_v}$ induces the identity on the
reduction of $\c M$ modulo $\varpi_v$.  Then $\c M$ becomes trivial over $K((x))$.}
\label{corfond}

\vpar\noindent
{\bf Remark.}
The $q$-difference module $\c M$ becomes trivial over $K((x))$ 
if and only if there exists a basis $\u e$ of ${\c M}_{K(x)}$ over $K(x)$ such that the
associated $q$-difference system has a fundamental 
matrix of solutions $Y(x)$ with coefficients in $K\[[x\]]$. 
In this case the matrix $G(x)$ defined by $\Dq(\u e)=\u eG(x)$ has no poles at zero. 
This implies that the matrix $G_n(x)$, defined by $\Dq^n(\u e)=\u
eG_n(x)$, have no poles at zero, for any positive integer $n$, and
hence that 
$$
Y(x)=Y(0)\l({\b I}_\mu+\sum_{n\geq 0}{G_n(0)\over[n]_q!}x^n\r)\ .
$$

\Dim
By theorem (\pIIIkatz) we know that $\c M$ is a regular singular
$q$-difference module. By the formal classification
(\pformclass), there exists a $K((x))$-basis $\u f$ of $K((x))\otimes_{{\c F}}M$ 
such that $\Phi_q(\u f)=\u fA$, with $A\in Gl_\mu(K)$ in Jordan normal form.
By (\pAo), we can choose $\u f$ such that for almost all $v\in\Sg_f$
we have 
$$
A^{\kappa_v }-1\equiv 0\hbox{\ mod $\varpi_v$}.
$$
We deduce that $A$ is actually a diagonal matrix and that the eigenvalues of $A$ are in $q^\interi$. 
We can assume $A={\b I}_\mu$ by applying a ``shearing transformation" ({\it cf.}\ [PS, page
154]),  {\it i.e.} a basis change of the form 
$$\pmatrix
{x^{n_1}{\b I}_{\nu_1} & & \cr
&\ddots &\cr
 & &x^{n_r}{\b I}_{\nu_r} }\ ,
$$
where $\nu_1,\dots,\nu_r$ are positive integers such
that $\sum \nu_i=\mu$ and $n_1,\dots,n_r\in\interi$. 
Let $\u e$ be a basis of $M$ over ${\c F}$. Then there exists 
$F(x)\in Gl_\mu\l(K((x))\r)$ such that $\u e=\u fF(x)$. It follows that 
$$
\Phi_q\u e=\u f F(qx)=\u e F(x)^{-1}F(qx)\ .
$$
Then $F(x)$ is a fundamental matrix of solutions 
for the $q$-difference system associated to 
$\c M$ with respect to the basis $\u e$.
After a change of basis of the form $\u e^\p=x^mC\u e$, 
where $m\in\interi$ and $C$ is a constant invertible matrix, 
we obtain a $q$-difference system having a solution 
$Y(x)={\b I}_\mu+\sum_{m\geq 1}Y_m x^m\in Gl_\mu(K\[[x\]])$. 
\QED 

%%%%%%%%%%%%%%%%%%%%%%%%%%%%%%%%%%%%%%%%%%%%%%%%%%%%%%%%%%%%%%%%%%%%%%%%%
%%%%%%%%%%%%%%%%%%%%%%%%%%%%   7  %%%%%%%%%%%%%%%%%%%%%%%%%%%%%%%%%%%%%%%
%%%%%%%%%%%%%%%%%%%%%%%%%%%%%%%%%%%%%%%%%%%%%%%%%%%%%%%%%%%%%%%%%%%%%%%%%

\newchapter{Statement of the {\tiit q}-analogue of} 
\smallskip
\centerline{\ti Grothendieck's conjecture on {\tiit p}-cur\-va\-tures}
\chaplabel{groth}
%%%%%%%%%%%%%%%%%%%%%%%%%%%%%%%%%%%%%%%%%%%%%%%%%%%%%%%%%%%%%%%%%%%%%%%%%

%%%%%%%%%%%%%%%%%%%%%%%%%%%%%%%%%%%%%%%%%%%%%%%%%%%%%%%%%%%%%%%%%%%%%%%%%
%%%%%%%
\newsubchapter{Statement of the theorem}
%%%%%%%
%%%%%%%%%%%%%%%%%%%%%%%%%%%%%%%%%%%%%%%%%%%%%%%%%%%%%%%%%%%%%%%%%%%%%%%%%

\vpar
We recall that $K$ is a number field, ${\c V}_K$ its ring of integers,
$v$ is a finite or an infinite place of $K$, and $q$ an element of $K$, which is not a root of
unity.  The uniformizer of the finite place $v$ is denoted by $\varpi_v$. 
For almost all finite places $v$, let $\kappa_v$ be the 
the multiplicative order of the image of $q$ in the residue field of 
$M$ modulo $\varpi_v$.  Let $\varpi_{q,v}$ be the integer power of
$\varpi_v$ such that $|\varpi_{q,v}|_v=|1-q^{\kappa_v}|_v$. 
\par
We consider a $q$-difference algebra ${\c F}\subset K(x)$ essentially of
finite type over ${\c V}_K$ and a $q$-difference module ${\c M}=(M,\Phi_q)$ over $\c F$. 

\vpar
We want to prove the following theorem, 
which we consider to be the $q$-analogue of the Grothendieck conjecture 
for differential equations with
$p$-curvature zero for almost all finite places:

\T
{Let ${\c M}=(M,\Phi_q)$ be a $q$-difference module over $\c F$, 
such that 
$$
\matrix{\vbox{\hsize 12 true cm\noindent
the operator $\Phi_q^{\kappa_v}$ induces the identity
on the reduction of $\c M$ modulo $\varpi_{q,v}$ for almost
all finite places $v$.}}
\leqno{(\ast)}$$
Then $\c M$ becomes trivial over $K(x)$.}
\label{grothendieck}

\Oss 
We recall that the $q$-difference module $\c M$ over $K(x)$ is trivial 
if and only if 
the following equivalent conditions are satisfied:
\tparnoind
1) there exists an isomorphism of $q$-difference modules 
$M^{\Phi_q}\otimes_K K(x)\cong M$;
\tparnoind
2) there exists a $K(x)$-vector space isomorphism 
$\psi:M \lrarrow K(x)^\mu$ such that for all $m\in M$ we have:
$\psi(\Phi_q(m))=\varphi_q(\psi(m))$, where $\varphi_q$ is defined component-wise on $K(x)^\mu$.
\tparnoind
3) there exists a basis $\u e$ of $M$ over $K(x)$ such that, if $\Dq\u e=\u e G(x)$,  
we can find $Y(x)\in Gl(K(x))$ satisfying the $q$-difference linear
system $\dq Y(x)=Y(x)G(x)$.
\par
It is clear that if a $q$-difference module $\c M$ over $\c F$
becomes trivial over $K(x)$, the hypothesis $(\ast)$ of the theorem above is
satisfied. 

\vpar
By (\pgrothradici) we immediately obtain:

\C
{Let ${\c M}=(M,\Phi_q)$ be a $q$-difference module over $\c F$
such that the reduction of $\c M$ modulo $\varpi_{q,v}$ is trivial for
almost all $v$.  Then $(M,\Phi_q)$ is trivial over $K(x)$.}

\vpar
In (\pgrothendieck), we assumed that $q$ is not a root of
unity: if $q$ is a root of unity, theorem (\pgrothendieck) is an easy
consequence of the results in (\ccurva).  We notice that in this particular case,      
we just need the hypothesis of trivial reduction modulo $\varpi_v$:

\P
{Let $q$ be a primitive $\kappa$-root of unity, with $\kappa\geq 1$.
Then the $q$-difference module ${\c M}=(M,\Phi_q)$ over $\c F$ becomes 
trivial over $K(x)$ if and only if $\c M$ has trivial reduction 
modulo $\varpi_v$ for an infinite number of $v\in\Sg_f$.}

\Dim
Let $\u e$ be a basis of $M$ over $\c F$
and let $\Phi_q^m(\u e)=\u eA_m(x)$ for all $m\geq 1$.
By (\pgrothradici) it is enough to prove that 
$$
A_\kappa(x)={\b I}_\mu\Leftrightarrow A_{\kappa_v}(x)\equiv {\b I}_\mu
\hbox{\ modulo $\varpi_v$, infinitely many $v\in\Sg_f$.}
$$ 
To conclude, it is enough to notice that $\kappa_v=\kappa$ for almost
all $v\in\Sg_f$. 
\QED

%%%%%%%%%%%%%%%%%%%%%%%%%%%%%%%%%%%%%%%%%%%%%%%%%%%%%%%%%%%%%%%%%%%%%%%%%
%%%%%%%
\newsubchapter{Idea of the proof}
%%%%%%%
%%%%%%%%%%%%%%%%%%%%%%%%%%%%%%%%%%%%%%%%%%%%%%%%%%%%%%%%%%%%%%%%%%%%%%%%%

\vpar
The proof of (\pgrothendieck) is inspired by the theory of
$G$-functions, from which we derive the definitions below. 
In the $q$-difference case, they are not as interesting as in the
differential case; in fact, as we will see later, the two invariants
that we are going to define are finite only 
when the $q$-difference module is trivial over $K(x)$.
In any case, they will be useful in some intermediate steps of the proof.

\Def
{Let $y=\sum_{n=0}^\infty a_n x^n\in K\[[x\]]$.  
We set $h(y,n,v)=\sup_{|\u\a|\leq n}\l(\lgp |a_{\u\a}|_{v}\r)$ and 
we define the {\it size of $y$} ({\it cf.}\ {\rm [A1, I, 1.3]}) to be the number 
$$
\sg(y)=\limsup_{n\rightarrow\infty}{1\over n}
\sum_{v\in\Sg_f\cup\Sg_\infty}h(y,n,v)\ .
$$
Let ${\c M}=(M,\Phi_q)$ be a $q$-difference module over a $q$-difference algebra 
${\c F}\subset K(x)$.  We fix a basis $\u e$
of ${\c M}_{K(x)}$ over $K(x)$ and we set as usual $\Dq^n\u e=\u
eG_n(x)$ and $h(M,n,v)=\sup_{0\leq s\leq
  n}\log\l|{G_s(x)\over[s]_q!}\r|_{v,Gauss}$.
We define the {\it size of $\c M$} to be 
$$
\sg({\c M})=\limsup_{n\rightarrow\infty}{1\over n}
\sum_{v\in\Sg_f\atop |1-q^{\kappa_v}|_v<|p|_v^{1/(p-1)}}h(M,n,v)\ .
$$}

\vpar
The proof (to be given explicitly in the next section) 
is organized as follows:
$$\matrix{\cr
\boxit{\vbox{\hsize 4.7 true cm\noindent
     ${\c M}$ satisfies the hypothesis $(\ast)$}}\cr\cr 
\hbox to 8 true cm{\upbracefill}\cr\cr
\matrix{\hskip -0.8 true cm (\pcorfond)\l\Downarrow\vbox to 15pt{}\r. & 
     &\hskip -0.8 true cm (\stepb)\l\Downarrow\vbox to 15pt{}\r.\cr\cr
\boxit{\vbox{\hsize 5.6 true cm\noindent 
     $(M,\Phi_q)$ becomes trivial over $K((x))$}}
     &&
     \vbox{\hsize 6 true cm\noindent
     \centerline{\boxit{\vbox{\hsize 1.8 true cm
     \noindent$\sg({\c M})<\infty$}}}}}\cr\cr
\hbox to 13 true cm{\upbracefill}\cr\cr
\hskip -0.8 true cm (\stepc)\l\Downarrow\vbox to 15pt{}\r.\cr\cr
\boxit{\vbox{\hsize 11true cm\noindent
     There exists a basis $\u e$ of ${\c M}_{K(x)}$ over $K(x)$ such that 
     $\Dq\u e=\u eG(x)$ and there exists an invertible matrix 
     $Y(x)\in Gl_\mu(K\[[x\]])$ such that $\dq Y(x)=Y(x)G(x)$ 
     and the entries of $Y(x)$ have finite size.}}\cr\cr
\hskip -0.8 true cm (\stepd)\l\Downarrow\vbox to 15pt{}\r.\cr\cr
\boxit{\vbox{\hsize 2.9 true cm\noindent $Y(x)\in Gl_\mu(K(x))$}}}
$$

%%%%%%%%%%%%%%%%%%%%%%%%%%%%%%%%%%%%%%%%%%%%%%%%%%%%%%%%%%%%%%%%%%%%%%%%%
%%%%%%%%%%%%%%%%%%%%%%%%%%%%   8  %%%%%%%%%%%%%%%%%%%%%%%%%%%%%%%%%%%%%%%
%%%%%%%%%%%%%%%%%%%%%%%%%%%%%%%%%%%%%%%%%%%%%%%%%%%%%%%%%%%%%%%%%%%%%%%%%

\newchapter{Proof of (\pgrothendieck)}

%%%%%%%%%%%%%%%%%%%%%%%%%%%%%%%%%%%%%%%%%%%%%%%%%%%%%%%%%%%%%%%%%%%%%%%%%

%%%%%%%%%%%%%%%%%%%%%%%%%%%%%%%%%%%%%%%%%%%%%%%%%%%%%%%%%%%%%%%%%%%%%%%%%
%%%%%%%
\newsubchapter{Finiteness of the size of $\c M$}
%%%%%%%
%%%%%%%%%%%%%%%%%%%%%%%%%%%%%%%%%%%%%%%%%%%%%%%%%%%%%%%%%%%%%%%%%%%%%%%%%

\P
{Let ${\c M}$ be a $q$-difference module over ${\c F}\subset K(x)$ 
satisfying $(\ast)$. 
Then
$$
\sg({\c M})<+\infty\ .
$$}
\label{sgrho}

\Dim
By assumption ({\it cf.} (\funipotenti) and (\punipotenza)) 
there exists a basis $\u e$ of $M$ over ${\c F}$
such that $\Dq^m\u e=\u eG_m(x)$ and 
$$
\hbox{$|G_{\kappa_v}(x)|_{v,Gauss}\leq |[\kappa_v]_q|_v$ for almost all $v\in\Sg_f$.}  
$$
For such $v$, by 
(\ftrivG) we have 
$$
\l|{G_n(x)\over[n]_q!}\r|_{v,Gauss}
\leq\l|{G_{\l[{n\over\kappa_v}\r]\kappa_v}(x)\over[n]_q!}\r|_{v,Gauss}
\leq {|[\kappa_v]_q|_v^{\l[{n\over\kappa_v}\r]}\over |[n]_q!|_v}
\leq|p|_v^{-n/\kappa_v(p-1)}\ ,
$$
from which we obtain 
$$
h(M,n,v)\leq n{\log|p|_v^{-1}\over\kappa_v(p-1)}\ .
$$
Let 
$$
T_1=\{v\in\Sg_f:|1-q^{\kappa_v}|_v<|p|_v^{1/(p-1)},~|G_{\kappa_v}(x)|_{v,Gauss}\leq |[\kappa_v]_q|\}
$$
and
$$
T_2=\{v\in\Sg_f:|1-q^{\kappa_v}|_v<|p|_v^{1/(p-1)},~|G_{\kappa_v}(x)|_{v,Gauss}>|[\kappa_v]_q|\}\ .
$$
The assumption $(\ast)$ implies that $T_2$ is finite. 
\par
By (\flimh) we deduce that 
$$
\sg({\c M})\leq
\sum_{v\in T_1}
{\log|p|_v^{-1}\over\kappa_v(p-1)}
+\sum_{v\in  T_2}
\lgp{1\over\chi_v(M)}\ .
$$
Let us consider the set $\widetilde\Sg$ of all finite places $v$ of $K$ such
that $|q|_v=1$.  Then $T_1$ is cofinite in $\widetilde \Sg$; hence, to conclude
that $\sg({\c M})$ is finite, it is enough to prove the lemma:

\L
{Let $q\neq 0$ be an element of $K$ which is not a root of unity, 
$\widetilde\Sg$ be the set of all finite places $v$ of $K$ such
that $|q|_v=1$, and $\kappa_v$ be the multiplicative order of the image of $q$ in the
residue field of $K$ with respect to $v\in\widetilde\Sg$. 
Then
$$
\sum_{v\in\widetilde\Sg}{\log|p|_v^{-1}\over\kappa_v(p-1)}<\infty\ .
$$}

\Dim  
Let $0<\veps<1$. 
We consider the sets:
$$\eqalign{
&S_0=\{v\in\widetilde\Sg,~\kappa_v\geq p-1\}\ ,\cr
&S_1=\{v\in\widetilde\Sg,~\kappa_v<p-1~,\kappa_v^2\geq p^{1+\veps}\}\ ,\cr
&S_2=\{v\in\widetilde\Sg,~\kappa_v<p-1~,\kappa_v^2< p^{1+\veps}\}\ .\cr}
$$
Then, for $\eta=(1-\veps)/(1+\veps)$ and $v\in S_2$, we have 
$$
p>\kappa_v^{2/(1+\veps)}~\Rightarrow
p-1\geq\kappa_v^{2/(1+\veps)}=\kappa_v^{1+\eta}\ ,
$$
and hence we obtain 
$$
\sum_{v\in\widetilde\Sg}{\lg|p|_v^{-1}\over \kappa_v(p-1)}
\leq\sum_{v\in S_0}{\lg|p|_v^{-1}\over(p-1)^2}+
    \sum_{v\in S_1}{\lg|p|_v^{-1}\over p^{1+\veps}}+
    \sum_{v\in S_2}{\lg|p|_v^{-1}\over\kappa_v^{2+\eta}}\ .
$$
The sums over $S_0$ and $S_1$ are clearly convergent.  
Since for almost all $v\in\widetilde\Sg$ we have
$|1-q^{\kappa_v}|_v^{-1}\geq |p|_v^{-1}$, 
to conclude that 
the sum over $S_2$ is convergent it is enough to prove that 
$$
\sum_{v\in\widetilde\Sg}{\lg|1-q^{\kappa_v}|_v^{-1}\over\kappa_v^{2+\eta}}
$$
is convergent for all $\eta>0$.
We recall that $\widetilde\Sg$ is cofinite in $\Sg_f$ and that for all
integers $n\geq 1$ there exists at worst only a finite number of $v\in\Sg_f$ such
that $\kappa_v=n$ (since $|1-q^n|_v=1$ for almost every $v\in\Sg_f$). 
Therefore by the Product Formula, we get 
$$\eqalign{\sum_{v\in\widetilde\Sg}{\lg|1-q^{\kappa_v}|_v^{-1}\over\kappa_v^{2+\eta}}
&=\sum_{n=1}^\infty\sum_{v\in\widetilde\Sg\atop\kappa_v=n}
     {\lg|1-q^n|_v^{-1}\over n^{2+\eta}}\cr 
&=\sum_{n=1}^\infty\sum_{v\in\widetilde\Sg,~\kappa_v\neq n\atop
     {\rm or}~v\in(\Sg_f\cup\Sg_\infty)\smallsetminus\widetilde\Sg}
     {\lg|1-q^n|_v\over n^{2+\eta}}\ .\cr}
$$
For every $v\in\Sg_f$ such that $|q|_v\leq 1$ we have 
$|1-q^n|_v\leq 1$. In particular $|1-q^n|_v= 1$ for almost all $v\in\Sg_f$. 
\par\vbox{\noindent
Therefore we obtain 
$$\eqalign{\sum_{v\in\widetilde\Sg}{\lg|1-q^{\kappa_v}|_v^{-1}\over\kappa_v^{2+\eta}}
&\leq\sum_{n=1}^\infty\sum_{v\in\Sg_f,~|q|_v>1\atop {\rm or}~v\in\Sg_\infty}
     {\lg|1-q^n|_v\over n^{2+\eta}}\cr 
&\leq \sum_{n=1}^\infty\l(\sum_{v\in\Sg_\infty,~|q|_v\leq 1}
     {\lg 2\over n^{2+\eta}}+\sum_{v\in\Sg_f\cup\Sg_\infty\atop |q|_v>1}
     {\lg(1+|q|_v)\over n^{1+\eta}}\r)<\infty\ .\cr}
\numero{arch}$$
\QED}

%%%%%%%%%%%%%%%%%%%%%%%%%%%%%%%%%%%%%%%%%%%%%%%%%%%%%%%%%%%%%%%%%%%%%%%%%
%%%%%%%
\newsubchapter{Finiteness of the size of a fundamental matrix of solutions}
%%%%%%%
%%%%%%%%%%%%%%%%%%%%%%%%%%%%%%%%%%%%%%%%%%%%%%%%%%%%%%%%%%%%%%%%%%%%%%%%%

\vpar
We have already proved that a $q$-difference module ${\c M}$ satisfying
$(\ast)$ becomes trivial over $K((x))$ ({\it cf.}\ (\pcorfond)) and
that it has finite size ({\it cf.}\ (\psgrho)). 

\P {Under the hypothesis $(\ast)$, there exists a basis $\u e$ of ${\c M}_{K(x)}$ 
over $K(x)$ such that the $q$-difference system associated
to  ${\c M}_{K(x)}$ with respect to $\u e$  has an invertible solution matrix
$Y(x)\in Gl_\mu(K\[[x\]])$, whose entries have finite size.}
\label{finitesize}

\Oss
We recall a property of the size of a formal power
series that we will use in the proof below.
We know ({\it cf.}\ for example [DGS, VI, 4.1]) that 
$$
\limsup_{n\rightarrow\infty} {1\over n}h(n,v,y)=\lgp{1\over r_v(y)}\ ,
$$
where $r_v(y)$ is the $v$-adic radius of convergence of $y$.
Then $\sg(y)<\infty$ if and only if there exists a finite set 
$S\subset\Sg_f\cup\Sg_\infty$ such that $y$ has nonzero radius of 
convergence for all $v\in S$ and 
$$
\limsup_{n\rightarrow\infty} {1\over n}
\sum_{v\in \Sg_f\cup\Sg_\infty\smallsetminus S}h(n,v,y)<\infty\ .
$$
\label{osservazionesigma}

\Dim
As in the proof of (\pcorfond), we can find a basis $\u e$ of ${\c M}_{K(x)}$ such that the associated 
$q$-difference system has a solution $Y(x)\in Gl_\mu(K\[[x\]])$ such
that $Y(0)={\b I}_\mu$.
Let $\Dq^n\u e=\u eG_n(x)$.  Since $G_0=Y(0)={\b I}_\mu$,
we conclude that $Y(x)=\sum_{n\geq 0}{G_n(0)\over[n]_q!}x^n$. 
\par
We want to prove that the entries of $Y(x)$ have finite size.
By (\pbuonadef) and (\pqmerdosi), we deduce that 
$$
\limsup_{n\rightarrow\infty}{1\over n}\sum_{v\in\Sg_f\atop |q|_v\leq 1}
      \sup_{0\leq s\leq n}\lgp\l|Y_n\r|_v
\leq\sg({\c M})+\sum_{v\in\Sg_f, |q|_v<1\atop{\rm or}~ 1>|1-q^{\kappa_v}|_v\geq |p|_v^{1/(p-1)}}
      \lgp{1\over\chi_v(M)}<\infty\ .
$$
By the previous remark
it is enough to prove that the entries of $Y(x)$ have nonzero radius of
convergence for all $v\in\Sg_\infty$ and for all $v\in\Sg_f$ such that
$|q|_v>1$.  Since $\c M$ is a regular $q$-difference module over $K(x)$, each
entry of $Y(x)$ is the solution of a regular singular $q$-difference
equation, by (\pnewton).  It follows from [Be, IV] and [BB, IV]
\footnote{$\hskip -3pt{}^{(2)}$}
{\rm In [Be] and [BB] the authors assume $|q|<1$.  This is only a 
matter of convention and their results translate to our situation.}
that the entries of $Y(x)$ have nonzero radius of
convergence for all $v\in\Sg_\infty$ such that $|q|_v\neq 1$ 
and for all $v\in\Sg_f$ such that $|q|_v>1$.
\par
The conclusion of the proof of (\pfinitesize) is the subject of the next section. 

%%%%%%%%%%%%%%%%%%%%%%%%%%%%%%%%%%%%%%%%%%%%%%%%%%%%%%%%%%%%%%%%%%%%%%%%%
%%%%%%%
\newsubchapter{How to deal with the problem of archimedean small divisors}
%%%%%%%
%%%%%%%%%%%%%%%%%%%%%%%%%%%%%%%%%%%%%%%%%%%%%%%%%%%%%%%%%%%%%%%%%%%%%%%%%

\vpar
To conclude that the entries of $Y(x)$ 
have finite radius of convergence for all infinite places $v$ such that
$|q|_v=1$, we recall the following result:

\vpar
\centerline{\vbox{\hsize 15 true cm\noindent
{\bf [Be, 6.1]}
{\sl Let 
${\c L}=\sum_{i=0}^\mu\sum_{j=0}^\nu a_{i,j}x^j\varphi_q^i\in\complessi[x,\varphi_q]$ 
be a $q$-difference operator and let $Q(x)=\sum_{i=0}^\mu a_{i,j_0}x^i$ be a
polynomial such that $j_0=\min\{j=0,\dots,\nu: a_{i,j}\neq 0\}$. 
We suppose that $|q|_\complessi=1$ and that there 
exist two positive real constants $c_1$ and $c_2$ such that
all the roots $u$ of the polynomial $(x-1)Q(x)$ satisfy the inequality 
$|q^n-u|_v\geq c_1n^{-c_2}$ for $n>>0$.  Then a formal power series 
$y\in\complessi\[[x\]]$ which solves ${\c L}y=0$ is convergent.}}}

\vpar\noindent
It is enough to prove that for any finite place $v$ such that 
$|q|_v=1$, there exist two positive real constants 
$c_{1,v}$ and $c_{2,v}$ such that 
$$
|q^n-u|_v\geq c_{1,v}n^{-c_{2,v}}
\numero{qu}$$ 
for $n>>0$.  To verify (\fqu) we will use the following theorem by Baker
\footnote{$\hskip -3pt{}^{(3)}$}
{\rm There is a proof of the additive version of this theorem in [B], 
but we prefer to cite the version in [Se], which is more suitable to our situation.}:  

\vpar
\centerline{\vbox{\hsize 15 true cm\noindent
{\bf [Se, 8.2, Corollary]}
{\sl Let $K$ be a number field, $\a_1,\dots,\a_l\in K$, 
$\be_1,\dots,\be_l\in\interi$ and $v$ a place of $K$.  If 
$\a_1^{\be_1}\cdots\a_l^{\be_l}\neq 1$ then 
$$
\l|\a_1^{\be_1}\cdots\a_l^{\be_l}-1\r|_v\geq\sup(4,\be_1,\dots,\be_l)^{-const}\ ,
$$
where the constant depends only on $v$ and on $\a_1,\dots,\a_l\in K$.}}}

\vpar\noindent
Let $l=2$, $\a_1=q$, $\a_2=u$, $\be_1=n$ and $\be_2=-1$.  
Since for $n>>0$ we have $q^nu^{-1}\neq 1$, we obtain 
$|q^n-u|_v\geq |u|_v n^{-c(u)}$.  Here $c(u)$ is a constant depending
on $u$, $q$ and $v$.
We set 
$$
c_{1,v}=\sup\l(\{|u|_v:\hbox{such that $Q(u)=0$}\}\cup\{1\}\r)
$$
and 
$$
c_{2,v}=\sup\{c(u):\hbox{$u$ such that $Q(u)=0$ or $u=1$}\}\ ;
$$
then we obtain the desired inequality. 
This achieves the proof of (\pfinitesize).
\QED

%%%%%%%%%%%%%%%%%%%%%%%%%%%%%%%%%%%%%%%%%%%%%%%%%%%%%%%%%%%%%%%%%%%%%%%%%
%%%%%%%
\newsubchapter{Conclusion of the proof: a criterion for rationality}
%%%%%%%
%%%%%%%%%%%%%%%%%%%%%%%%%%%%%%%%%%%%%%%%%%%%%%%%%%%%%%%%%%%%%%%%%%%%%%%%%

\vpar
We complete the proof of (\pgrothendieck) by applying 
the following proposition:

\P
{{\bf (Y. Andr{\'e})} Let $y(x)\in K\[[x\]]$ be a formal power series 
solution of a $q$-difference equation 
$$
a_\mu(x)\dq^\mu(y)(x)+a_{\mu-1}(x)\dq^{\mu-1}(y)(x)+\cdots+a_0(x) y(x)=0\ ,
$$
with $a_i(x)\in K(x)$ for all $i=0,\dots,\mu$.  If $\sg(y)<\infty$ then $y(x)$ is the 
Taylor expansion of a rational function $\in K(x)$.}
\label{criterioandre}

\Dim
It is a general property of $q$-difference equations 
that for all $v\in\Sg_f\cup\Sg_\infty$ such that 
$|q|_v\neq 1$, the series $y(x)$ with nonzero radius of convergence 
has infinite radius of meromorphy 
({\it cf.} for instance [BB, 7.2]).  We remark that we can always find
such a $v$ since $q$ is not a root of unity. 
To conclude that $y(x)$ is a rational 
function it is enough to apply the more general result
[A1, VIII, 1.1, Th.], but we prefer to sketch the proof, since it
simplifies under the present assumptions. 
First we prove that $y(x)$ is an algebraic function (steps from 1 to 5), following  
the proof of [A2, 2.3.1] adapted to this particular case. 
Then, in step 6, we prove that $y(x)$ is the expansion of a rational function.

\smallskip\noindent
$\u{\hbox{\it Step 1.}}$
We fix $\eta\in(0,1]$ and an integer $\nu>1$. 
Let 
$$
\ovec Y={}^t(1,y(x),\dots,y(x)^{\nu-1})=\sum_{m\geq 0}\ovec Y_m x^m\in
K\[[x\]]^\nu\ .
$$
Using Siegel's Lemma, one can construct a polynomial vector   
$$
\ovec{P_N}(x)=(P_{N,0}(x),\dots,P_{N,\nu-1}(x))=
\sum_{m\geq 0}\ovec P_N^{(m)}x^m\in K[x]^\nu\ , 
$$
for all $N\in\nat$, such that 
$$
M=\ord_0(\ovec{P_N}\cdot\ovec Y)\geq N\ ;
\leqno{\it i)}$$
$$
\deg_x\ovec{P_N}=\sup_{i=0,\dots,\nu-1}(\deg_x
P_{N,i}(x))\leq {1\over \nu}\l(1+{1\over\eta}\r)N+o(N)\ ;
\leqno{\it ii)}$$
$$
h(\ovec{P_N})=\limsup_{m\rightarrow\infty}
{1\over m}\sum_{v\in\Sg_f\cup\Sg_\infty}\lgp\l(\sup_{m\geq 0}
|\ovec{P_N}^{(m)}|_v\r)\leq\eta N\sg(\ovec Y)+o(N)\ .
\leqno{\it iii)}$$

\smallskip\noindent
$\u{\hbox{\it Step 2.}}$
Let us suppose that $y(x)$ is {\it not} algebraic; 
then $\ovec{P_N}\cdot\ovec Y\neq 0$ for all $N\geq 1$ and
hence $M<\infty$.  We set 
$$
\a={1\over M!}{d^M\over dx^M}(\ovec{P_N}\ovec Y)(0)\ .
$$
Of course, $\a\neq 0$.  For all $v\in\Sg_f$ we have 
$$
\lg|\a|_v\leq \lgp|\ovec{P_N}(x)|_{v,Gauss}+\sup_{m\leq M}\lgp|\ovec Y_m|_v\ .
\numero{prima}$$

\smallskip\noindent
$\u{\hbox{\it Step 3.}}$
Let $V$ be a finite subset of 
$\Sg_f\cup\Sg_\infty$ containing $\Sg_\infty$ and at least one 
place $v$ such that the radius of meromorphy $M_v(y)$ is infinite,
since $q$ is not a root of unity. 
For all $v\in V$ the formal power series $y(x)$ is the germ at zero of
a meromorphic function, therefore we can write 
$y(x)=f_v(x)/g_v(x)$, where $f_v(x)$ and $g_v(x)$ are $v$-adic
analytic functions converging for $|x|_v<M_v(y)$. We can suppose 
that $g_v(0)=1$.  We set:
$$\eqalign{
&\ovec Z_v(x)=(g_v(x)^{\nu-1}, g_v(x)^{\nu-2}f_v(x),\dots,f_v(x)^{\nu-1}),\cr
&\psi_v(x)=\ovec{P_N}(x)\cdot\ovec Z_v(x);\cr}
$$
from which it follows that 
$$
\ovec{P_N}(x)\cdot\ovec Y(x)={1\over g_v(x)^{\nu-1}}\psi_v(x)\ .
$$ 
We deduce that:
$$
\a={1\over M!}{d^M\over dx^M}(\psi_v)(0)\ .
$$

\smallskip\noindent
$\u{\hbox{\it Step 4.}}$
Let us fix $m_v<M_v(y)$ for all $v\in V$.  By Cauchy's estimates we
obtain 
$$\eqalign{\lg|\a|_v
&\leq -M\lg m_v+\lg\l(\sup_{|x|_v=m_v}|\psi_v(x)|_v\r)\cr
&\leq -M\lg m_v+\lgp\l(\sup_{m\leq N}|\ovec{P_N}^{(m)}|_v\r)
 +m_v\deg_x\ovec{P_N}+o(N)\ .\cr}
\numero{seconda}$$

\smallskip\noindent
$\u{\hbox{\it Step 5.}}$
Summing (\fprima) for $v\in(\Sg_\infty\cup\Sg_f)\smallsetminus V$ 
and (\fseconda) for $v\in V$, by the Product Formula we obtain 
$$\eqalign{M\sum_{v\in V}\lg m_v
&\leq \sum_{v\in\Sg_f\cup\Sg_\infty}\lgp\l(\sup_{m\leq N}|\ovec{P_N}^{(m)}|_v\r)\cr
&\hskip 20pt+\sum_{v\not\in V}\sup_{m\leq M}\lgp|\ovec Y_m|_v
+\deg_x\ovec{P_N}\sum_{v\in V}m_v+o(N)\ ;}
$$
dividing by $M\geq N$ and taking
the $\limsup$ for $N\rightarrow\infty$ we have 
$$
\sum_{v\in V}\lg m_v\leq (\eta+1)\sg(y)+
{1\over \nu}\l(1+{1\over\eta}\r)\sum_{v\in V}m_v\ .
$$
Finally we can take the limit for $\nu\rightarrow\infty$
and $\eta\rightarrow 0$ and obtain 
$$
\sum_{v\in V}\lg m_v\leq\sg(y)\ .
$$
Since $\sg(y)<\infty$, we get a contradiction by letting 
$m_v\rightarrow M_v(y)$.
\smallskip
Thus we have proved that $y(x)$ is an algebraic function. 

\smallskip\noindent
$\u{\hbox{\it Step 6.}}$ Now we prove that $y(x)$ is the 
Taylor expansion of a rational function.
Since $y(x)$ is algebraic over $K(x)$ there exists $P(x)\in K[x]$ 
such that $g(x)=P(x)y(x)$ satisfies a relation of the form:
$$
g(x)^s-Q_1(x)g(x)^{s-1}-\dots-Q_s(x)=0\ ,
\numero{intero}$$
with $Q_1(x),\dots,Q_s(x)\in K[x]$. 
Let us fix $v\in\Sg_\infty\cup\Sg_f$ such that $|q|_v\neq 1$, hence such that 
$y(x)$ has infinite $v$-adic radius of meromorphy. 
The relation (\fintero) implies that $g(x)$ is an entire analytic function. 
Proving that $y(x)$ is a rational function is equivalent to prove that 
$g(x)=P(x)y(x)$ is a polynomial.
We deduce by (\fintero) that there exist two real positive constants 
$c_1$ and $c_2$ such that for $|x|_v>>1$ we have 
$$
|g(x)|_v\leq c_1|x|_v^{c_2}\ ,
$$  
which implies that $g(x)\in K[x]$.
%The field $\complessi_v(x,y)$ is the field of meromorphic functions 
%over a surface $E$. The immersion
%$\complessi_v(x)\hrarrow\complessi_v(x,y)$ 
%determines an algebraic covering $p:E\lrarrow {\b A}_{\complessi_v}^1$. 
%Since $y(x)$ is a meromorphic function over $E$, $p$ is an {\'e}tale
%covering of ${\b A}_{\complessi_v}^1$.  Such a covering is necessarily trivial, 
%hence $\complessi_v(x,y)=\complessi_v(x)$.
\QED

\vpar
This completes the proof of theorem (\pgrothendieck).

%%%%%%%%%%%%%%%%%%%%%%%%%%%%%%%%%%%%%%%%%%%%%%%%%%%%%%%%%%%%%%%%%%%%%%%%%
%%%%%%%
\newsubchapter{A corollary}
%%%%%%%
%%%%%%%%%%%%%%%%%%%%%%%%%%%%%%%%%%%%%%%%%%%%%%%%%%%%%%%%%%%%%%%%%%%%%%%%%

\vpar
We point out that the following corollary is a consequence the proof
of (\pgrothendieck):

\C
{Let $S$ be a subset of $\Sg_f$ having Dirichlet density 1 
and $\c M$ be a $q$-difference module over a $q$-difference algebra ${\c F}\subset K(x)$ 
essentially of finite type over ${\c V}_v$. 
We assume that for all $v\in S$ the operator $\Phi_q^{\kappa_v}$
induces the identity over the reduction of $\c M$ modulo $\varpi_{q,v}$. 
We assume moreover that:
$$
\sum_{v\in\Sg_f\smallsetminus S\atop|1-q^{\kappa_v}|_v<|p|_v^{1/(p-1)}}
\lg{1\over\chi_v(M)}<\infty\ .
$$ 
Then $\c M$ becomes trivial over $K(x)$.}
\label{almostgroth}

\Dim
We notice that in the proof of (\psgrho) we have actually shown that:
$$
\limsup_{n\rightarrow\infty}{1\over n}
\sum_{v\in S\atop |1-q^{\kappa_v}|_v\leq |p|_v^{1/(p-1)}}
h(M,n,v)<\infty\ .
$$
Let
$$
T=\{v\in\Sg_f\smallsetminus S:|G(x)|_{v,Gauss}\leq 1\ \hbox{and}\ |1-q^{\kappa_v}|_v<|p|_v^{1/(p-1)}\}\ .
$$
If we prove that 
$$
\limsup_{n\rightarrow\infty}{1\over n}
\sum_{v\in T}
h(M,n,v)
\leq\sum_{v\in T}
\lg{1\over\chi_v(M)}
\numero{nonnilp}
$$
then we obtain $\sg({\c M})<\infty$, since 
$|1-q^{\kappa_v}|_v<|p|_v^{1/(p-1)}$  and 
$|G(x)|_{v,Gauss}\leq 1$
for almost all $v\in\Sg_f$.  
Then we can complete the proof
as we did the proof of (\pgrothendieck). 
\par\vbox{\noindent
To prove (\fnonnilp) we need only notice that for all
$v\in\Sg_f\smallsetminus S$ 
such that $|1-q^{\kappa_v}|_v<|p|_v^{1/(p-1)}$ and
$|G(x)|_{v,Gauss}\leq 1$ 
and for all $n\geq 1$, we have 
$$
h(M,n,v)\leq
\lg\l(|[\kappa_v]_q|_v^{-\l[{n\over \kappa_v}\r]}
|p|_v^{-\l[{n\over \kappa_v}\r]{1\over p-1}}\r)
\leq n\lg{1\over\chi_v(M)}\ .
$$
\QED}

%%%%%%%%%%%%%%%%%%%%%%%%%%%%%%%%%%%%%%%%%%%%%%%%%%%%%%%%%%%%%%%%%%%%%%%%%
%%%%%%%%%%%%%%%%%%%%%%%% PART IV %%%%%%%%%%%%%%%%%%%%%%%%%%%%%%%%%%%%%%%%
%%%%%%%%%%%%%%%%%%%%%%%%%%%%%%%%%%%%%%%%%%%%%%%%%%%%%%%%%%%%%%%%%%%%%%%%%
%%%%%%%%%%%%%%%%%%%%%%%%%%%%%%%%%%%%%%%%%%%%%%%%%%%%%%%%%%%%%%%%%%%%%%%%%
\vskip 20pt 
\centerline{\titlepart Part IV.  A {\titlepartit q}-analogue of Katz's conjectural
description} 
\vpar
\centerline{\titlepart of the generic Galois group}
%%%%%%%%%%%%%%%%%%%%%%%%%%%%%%%%%%%%%%%%%%%%%%%%%%%%%%%%%%%%%%%%%%%%%%%%%
%%%%%%%%%%%%%%%%%%%%%%%%%%%%%%%%%%%%%%%%%%%%%%%%%%%%%%%%%%%%%%%%%%%%%%%%%
%%%%%%%%%%%%%%%%%%%%%%%%%%%%%%%%%%%%%%%%%%%%%%%%%%%%%%%%%%%%%%%%%%%%%%%%%
%%%%%%%%%%%%%%%%%%%%%%%%%%%%%%%%%%%%%%%%%%%%%%%%%%%%%%%%%%%%%%%%%%%%%%%%%

%%%%%%%%%%%%%%%%%%%%%%%%%%%%%%%%%%%%%%%%%%%%%%%%%%%%%%%%%%%%%%%%%%%%%%%%%
%%%%%%%%%%%%%%%%%%%%%%%%%%%%  9 %%%%%%%%%%%%%%%%%%%%%%%%%%%%%%%%%%%%%%%%%
%%%%%%%%%%%%%%%%%%%%%%%%%%%%%%%%%%%%%%%%%%%%%%%%%%%%%%%%%%%%%%%%%%%%%%%%%

\newchapter{Definition of the generic {\tiit q}-difference Galois group}
\chaplabel{defgal}

%%%%%%%%%%%%%%%%%%%%%%%%%%%%%%%%%%%%%%%%%%%%%%%%%%%%%%%%%%%%%%%%%%%%%%%%%

%%%%%%%%%%%%%%%%%%%%%%%%%%%%%%%%%%%%%%%%%%%%%%%%%%%%%%%%%%%%%%%%%%%%%%%%%
%%%%%%%
\newsubchapter{Some algebraic constructions}
%%%%%%%
%%%%%%%%%%%%%%%%%%%%%%%%%%%%%%%%%%%%%%%%%%%%%%%%%%%%%%%%%%%%%%%%%%%%%%%%%

\vpar
We consider the following algebraic constructions 
on the category of $q$-difference modules over a fixed $q$-difference algebra 
$\c F$ over a field $K$:

\vpar\noindent
$\u{\hbox{\it Dual $q$-difference module.}}$
Let ${\c M}=(M,\Phi_q)$ be a $q$-difference module over ${\c F}$. 
Let us consider the dual ${\c F}$-module $\V M=Hom_{{\c F}}(M,{\c F})$
of $M$: 
$\V M$ is naturally a $q$-difference module equipped with 
the $q$-difference operator $\V \Phi_q={}^t\Phi_q^{-1}$, defined by
$$
<\V \Phi_q(\V m),m>=<\V m,\Phi_q^{-1}(m)>\ ,\ 
\hbox{for all $\V m\in\V M$ and $m\in M$.}
$$
The $q$-difference module ${\c\V M}=(\V M,\V \Phi_q)$ over $\c F$ is the 
dual $q$-difference module of $\c M$.

\vpar\noindent
$\u{\hbox{\it Tensor product of $q$-difference modules.}}$
Let ${\c M}^\p=(M^\p,\Phi_q^\p)$ and ${\c M}^{\p\p}=(M^{\p\p},\Phi_q^{\p\p})$ be 
two $q$-difference modules of finite rank over ${\c F}$.
The tensor product $M^\p\otimes_{{\c F}}M^{\p\p}$ has the natural structure of
a $q$-difference module defined by:
$$
\Phi_q(m^\p\otimes m^{\p\p})
=\Phi_q^\p(m^\p)\otimes\Phi_q^{\p\p}(m^{\p\p})\ ,\ 
\hbox{for all $m^\p\in M^\p$ and $m^{\p\p}\in M^{\p\p}$.}
$$
The $q$-difference module ${\c M}^\p\otimes_{\c F}{\c M}^{\p\p}=(M^\p\otimes_{{\c F}}M^{\p\p},
\Phi_q^\p\otimes\Phi_q^{\p\p})$ over $\c F$ is the tensor product 
of ${\c M}^\p$ and ${\c M}^{\p\p}$.

\vpar
We denote by $<\!\!{\c M}\!\!>^\otimes$  the full subcategory of the category 
of the $q$-difference modules over $\c F$
containing all the subquotients of the $q$-difference modules obtained
as finite sums of the form $\oplus_{i,j}T^{i,j}({\c M})$, where 
$T^{i,j}({\c M})={\c M}^{\otimes^i}\otimes{\c\V M}^{\otimes^j}$.

%%%%%%%%%%%%%%%%%%%%%%%%%%%%%%%%%%%%%%%%%%%%%%%%%%%%%%%%%%%%%%%%%%%%%%%%%
%%%%%%%
\newsubchapter{Definition of the Galois group of a {\titit q}-difference module}
%%%%%%%
%%%%%%%%%%%%%%%%%%%%%%%%%%%%%%%%%%%%%%%%%%%%%%%%%%%%%%%%%%%%%%%%%%%%%%%%%

\vpar
{\it Let $K$ be a field and $q$ be a nonzero element of $K$ which is not a root of unity.} 
Let ${\c M}=(M,\Phi_q)$ be a $q$-difference module over $K(x)$.
\par
Sometimes, the tensor category $<\!\!{\c M}\!\!>^\otimes$
comes equipped with a $K$-linear fiber functor 
$$
\omega:<\!\!{\c M}\!\!>^\otimes\lrarrow \{\hbox{finite dimensional
  $K$-vector spaces}\}\ .
$$
In fact, such a fiber functor always exists after replacing $K$ by
some finite extension.  Here is an explicit construction.  The 
$q$-difference module ${\c M}$ admits a ``model'' $\widetilde {\c M}$ over some
$q$-difference algebra ${\c F}\subset K(x)$ essentially of finite type over ${\c V}_K$.
This provides a corresponding ``model''  $\widetilde {\c N}$
for any ${\c N}$ in $<\!\!{\c M}\!\!>^\otimes$. 
Consider a $K$-valued point $x$ of $\c F$ (which exists after passing
to a finite extension of $K$).  Then the fiber at $x$ provides a
``fiber functor''.  The corresponding tannakian group may be
interpreted as the Picard-Vessiot group of a $q$-difference system, as
considered in [PS] (this interpretation is not used in the sequel).  

\vpar
It follows ({\it cf.}\ [DM, 2.11] and [A2, III, 2.1.1]) 
that there exists an algebraic closed subgroup $Gal({\c M},\omega)$ of $GL(\omega({\c M}))$, 
such that $\omega$ induces a tensor equivalence of categories between 
$<\!\!{\c M}\!\!>^\otimes$ and the category of finite type representations 
of $Gal({\c M},\omega)$ over $K$. 

\Def
{The algebraic group $Gal({\c M},\omega)$ is the {\it Galois group 
of ${\c M}$ pointed at $\omega$}.}

\vpar
We recall the following results:

\L
{{\rm [A2, III, 2.1.1]}
The algebraic group $Gal({\c M},\omega)$ is the subgroup of $Gl(\omega({\c M}))$ 
which stabilizes $\omega({\c N})$ for all sub-objects $\c N$ of a finite sum 
$T^{i,j}({\c M})={\c M}^{\otimes^i}\otimes{\c\V M}^{\otimes^j}$.}

\Oss
We notice that $Gal({\c M},\omega)$ is a stabilizer in the sense of algebraic groups.

\L
{{\rm [A2, III, 2.1.4]}
The group $Gal({\c M},\omega)$ is trivial if and only if 
$\c M$ is a trivial $q$-difference module over $K(x)$ ({\it cf.}\ (\pdeftriviale)).}
\label{trivialeuno}

%%%%%%%%%%%%%%%%%%%%%%%%%%%%%%%%%%%%%%%%%%%%%%%%%%%%%%%%%%%%%%%%%%%%%%%%%
%%%%%%%
\newsubchapter{Definition of the generic Galois group of a {\titit q}-difference module}
\subchaplabel{defGalois}
%%%%%%%
%%%%%%%%%%%%%%%%%%%%%%%%%%%%%%%%%%%%%%%%%%%%%%%%%%%%%%%%%%%%%%%%%%%%%%%%%

\vpar
Let ${\c M}=(M,\Phi_q)$ be a $q$-difference module over $K(x)$. 
Let us consider the forgetful fiber functor ``underlying vector space'' 
$$
\eta:<\!\!{\c M}\!\!>^\otimes\lrarrow 
\{\hbox{$K(x)$-vector spaces}\}\ . 
$$
The functor $\u{Aut}^\otimes(\eta)$ defined on the category of  
commutative $K(x)$-algebras is representable by 
an algebraic group $Gal({\c M},\eta)$ over $K(x)$. 

\Def
{The algebraic group $Gal({\c M},\eta)$ is the {\it generic 
Galois group} of $\c M$.}

\Oss 
The generic Galois group $Gal({\c M},\eta)$ 
admits the following concrete description: it is the closed subgroup of
$GL(M)$ which stabilizes all the 
$q$-difference sub-modules in finite sums $\oplus_{i,j}T^{i,j}({\c M})$
({\it cf.}\ [A2, \S III, 2.2]), in the sense of algebraic groups.
Since $GL(M)$ is a noetherian algebraic variety,  $Gal({\c M},\eta)$ 
is defined as the stabilizer of a finite number of $q$-difference sub-modules 
${\c N}_1,\dots,{\c N}_r$ of some finite sums 
$\oplus_{i,j}T^{i,j}({\c M})$: this is equivalent to demanding that 
$Gal({\c M},\eta)$ be the stabilizer of the maximal exterior power of 
$\oplus_{i=1}^r{\c N}_i$ ({\it cf.}\ [W, A.2]). 
This shows that $Gal({\c M},\eta)$ can be defined as the stabilizer of a $q$-difference 
sub-module of rank 1 of a finite sum $\oplus_{h}T^{i_h,j_h}({\c M})$.

\vpar\noindent
{\bf Warning.}
A sub-$K(x)$-vector space of a finite sum $\oplus_{i,j}T^{i,j}(M)$ stabilized by 
the generic Galois group $Gal({\c M},\eta)$ is not necessarily a $q$-difference module.

\Oss
If $\omega$ is a fiber functor over $<\!\!{\c M}\!\!>^\otimes$, with values in $K$-spaces,  
the functor $\u{Isom}^\otimes(\omega\otimes_K 1_{K(x)},\eta)$ is representable by a 
$K(x)$-group scheme $\Sg({\c M},\omega)$, which is a torsor over 
$Gal({\c M},\omega)\otimes_K K(x)$ ({\it cf.}\ [DM, 3.2] and [A2, III, 2.2]), such that 
$Gal({\c M},\eta)=Aut_{Gal({\c M},\omega)\otimes_K K(x)}\Sg({\c M},\omega)$. 

\L
{A $q$-difference module $\c M$ over $K(x)$ is trivial if and only if $Gal({\c M},\eta)$
is the trivial group.}
\label{lemmatriviale}

\Dim
If a fiber functor exists, this follows from the previous remark and
(\ptrivialeuno).  In general, $\omega$ exists after replacing $K$ by a
finite extension, and the result follows by an easy Galois descent. 
\QED

%%%%%%%%%%%%%%%%%%%%%%%%%%%%%%%%%%%%%%%%%%%%%%%%%%%%%%%%%%%%%%%%%%%%%%%%%
%%%%%%%%%%%%%%%%%%%%%%%%%%%%  10  %%%%%%%%%%%%%%%%%%%%%%%%%%%%%%%%%%%%%%%%
%%%%%%%%%%%%%%%%%%%%%%%%%%%%%%%%%%%%%%%%%%%%%%%%%%%%%%%%%%%%%%%%%%%%%%%%%

\newchapter{An arithmetic description of the generic Galois group}
\chaplabel{curvature}

%%%%%%%%%%%%%%%%%%%%%%%%%%%%%%%%%%%%%%%%%%%%%%%%%%%%%%%%%%%%%%%%%%%%%%%%%

\vpar
Let $K$ be a number field and ${\c V}_K$ the ring of integers of $K$.
We denote by $\Sg_f$ the set of all finite places $v$ of $K$, 
by $\varpi_v\in{\c V}_K$ the uniformizer associated to $v$, and by 
${\c V}_v$ the discrete valuation ring of $K$ associated to $v$. 
\par
We choose an element $q\in K$ which is {\it not} a root of unity.    
For every finite place $v$ such that $q$ is a unit of ${\c V}_v$, we denote by 
$\kappa_v$ the order of the cyclic group generated by the image of $q$ in the residue field 
of ${\c V}_v$, {\it i.e.}:
$$
\kappa_v=\min\{\hbox{$m\in\interi$: $m>0$ and $1-q^m\in\varpi_v{\c V}_v$}\}\ .
$$ 
Let $\varpi_{q,v}$ be the power of $\varpi_v$ 
satisfying $1-q^{\kappa_v}\in\varpi_{q,v}{\c V}_v$ and 
$1-q^{\kappa_v}\not\in\varpi_v\varpi_{q,v}{\c V}_v$. 
We set $k_{q,v}={\c V}_K/\varpi_{q,v}{\c V}_K$.

%%%%%%%%%%%%%%%%%%%%%%%%%%%%%%%%%%%%%%%%%%%%%%%%%%%%%%%%%%%%%%%%%%%%%%%%%
%%%%%%%
\newsubchapter{Algebraic groups ``containing $\Phi_q^{\kappa_v}$ 
for almost all $v$''}
%%%%%%%
%%%%%%%%%%%%%%%%%%%%%%%%%%%%%%%%%%%%%%%%%%%%%%%%%%%%%%%%%%%%%%%%%%%%%%%%%

\vpar
Let ${\c M}=(M,\Phi_q)$ be a $q$-difference module over $K(x)$. 
One can always find a $q$-difference algebra ${\c F}\subset K(x)$ 
essentially of finite type over the ring of integers ${\c V}_K$ of $K$
and a $q$-difference module $\widetilde{\c M}=(\widetilde M,\widetilde \Phi_q)$ over $\c F$ such that 
${\c M}$ is isomorphic to $\widetilde {\c M}_{K(x)}$.

\vpar\noindent
{\bf Remark.}
Since $q^{\kappa_v}\equiv 1$ modulo $\pi_{q,v}$, $\Phi_q^{\kappa_v}$
induces a $({\c F}\otimes_{{\c V}_K}k_{q,v})$-linear morphism on 
$\widetilde M\otimes_{{\c V}_K}k_{q,v}$.

\vpar
Let $G$ be a closed algebraic subgroup of $GL(M)$.  By Chevalley's
theorem, $G$ is the stabilizer of a one-dimensional sub-$K(x)$-vector
space $L$ in a finite sum $\oplus_hT^{i_h,j_h}(M)$. 
Up to enlarging $\c F$, there exists an ${\c F}$-free module $\widetilde L$ such that 
$L\cong\widetilde L\otimes_{\c F}K(x)$. 

\Def
{The closed algebraic subgroup $G$ of $GL(M)$ {\it contains 
$\Phi_q^{\kappa_v}$ for almost all $v\in\Sg_f$} if,
for almost every finite place $v$ of $K$, 
$\widetilde L\otimes_{{\c V}_K} k_{q,v}$ is stable by
$\Phi_q^{\kappa_v}$ in 
$\oplus_hT^{i_h,j_h}(\widetilde M)\otimes_{{\c V}_K} k_{q,v}$.}
\label{contenerecurvature}

\Oss
We notice that the notion of an algebraic group over $K(x)$ containing
$\Phi_q^{\kappa_v}$ for almost all $v\in\Sg_f$ is well defined: 

\smallskip\noindent
{\it Independence of the choice of $\c F$ and $\widetilde L$:}
let $\widetilde L^\p$ and ${\c F}^\p$ be a different choice for $\widetilde L$ and
${\c F}$ in the previous definition.  Then there exists a third
${\c V}_K$-algebra ${\c F}^{\p\p}$ of the same form with  
${\c F},{\c F}^\p\subset {\c F}^{\p\p}$. 
So by extension of scalars, we may suppose that 
${\c F}={\c F}^\p$ and that $\widetilde L$ and $\widetilde L^\p$ are two different 
$\c F$-lattices of $L$.  By enlarging ${\c F}$, we may suppose that
there exists an $\c F$-linear isomorphism 
$\psi:\widetilde L\lrarrow \widetilde L^\p$.
For almost all $v\in\Sg_f$ the morphism $\psi$ induces a 
$({\c F}\otimes_{{\c V}_K}k_{q,v})$-linear isomorphism
$\widetilde L\otimes_{{\c V}_K}k_{q,v}\lrarrow \widetilde L^\p\otimes_{{\c V}_K}k_{q,v}$
commuting with the action of $\Phi_q^{\kappa_v}$. 

\smallskip\noindent
{\it Independence of the choice of $L$:}
this follows from the fact that $\widetilde L$ is a direct factor of an 
${\c F}$-lattice of $\oplus_hT^{i_h,j_h}(M)$, hence any 
${\c F}\otimes_{{\c V}_K}k_{q,v}$-linear automorphism 
$\oplus_hT^{i_h,j_h}(\widetilde M)\otimes_{{\c V}_K} k_{q,v}$
stabilizing $\widetilde L\otimes_{{\c V}_K}k_{q,v}$ 
comes from an $\c F$-linear automorphism 
of $\oplus_hT^{i_h,j_h}(\widetilde M)$ stabilizing 
$\widetilde L$. 

\L
{The smallest closed algebraic 
subgroup $G_{\Phi^\kappa}({\cal M})$ of $GL(M)$ which contains $\Phi_q^{\kappa_v}$ for almost all
$v\in\Sg_f$ is well defined.}

\Dim
Let $G_1$ and $G_2$ be two closed algebraic subgroups of $GL(M)$
containing $\Phi_q^{\kappa_v}$ for almost
all $v$.  Let $G_1$ (resp.\ $G_2$) be defined as the stabilizer 
of a line $L_1$ (resp.\ $L_2$) in some $\oplus_hT^{i_h,j_h}(M)$. 
Then the intersection of $G_1$ and $G_2$ is the algebraic group
stabilizing the lines $L_1$ and $L_2$, or equivalently the line
$\wedge^2(L_1\oplus L_2)$ ({\it cf.}\ [W, A2]). 
This implies that the intersection of two algebraic subgroups of $GL(M)$ containing the
$\Phi_q^{\kappa_v}$ for almost all $v\in\Sg_f$ is still an algebraic
subgroup of $GL(M)$ containing the
$\Phi_q^{\kappa_v}$ for almost all $v\in\Sg_f$, in the sense of 
definition (\pcontenerecurvature).
\par 
Moreover $GL(M)$ is an algebraic variety of finite dimension, hence
any descending chain of closed algebraic subgroups containing the
$\Phi_q^{\kappa_v}$ for almost all $v\in\Sg_f$ is stationary. 
\QED

\L
{Let ${\c N}=(N,\Phi_q)$ be an object of $<\!\!{\c M}\!\!>^\otimes$. Then the natural morphism
$$
\Omega_{\Phi^\kappa}:G_{\Phi^\kappa}({\c M})\lrarrow G_{\Phi^\kappa}({\c N})
$$
is surjective.}
\label{lemmaGF}

\Dim
Since the action of $\Phi_q$ on $\c N$ is induced by the action of
$\Phi_q$ on $\c M$, the image of the natural morphism 
$\Omega_{\Phi^\kappa}:G_{\Phi^\kappa}({\c M})\lrarrow GL({\c N})$ contains 
$\Phi_q^{\kappa_v}$ for almost all $v$, and hence contains
$G_{\Phi^\kappa}({\c N})$.
\par
Let us choose a line $L$ in some finite sum $\oplus_hT^{i_h,j_h}(N)$,
such that $G_{\Phi^\kappa}({\c N})$ is the stabilizer of $L$. 
Let $\widetilde L$ be an $\c F$-lattice
of $L$, defined over a $q$-difference algebra ${\c F}\subset K(x)$ 
essentially of finite type over ${\c V}_K$.  Since $\c N$ is an object 
of $<\!\!{\c M}\!\!>^\otimes$, $L$ is a line in a suitable 
subquotient of a finite sum $\oplus_lT^{i_l,j_l}(M)$.  Since 
$\widetilde L\otimes_{{\c V}_K}k_{q,v}$ is stable by 
$\Phi_q^{\kappa_v}$ for almost all $v$, 
$L$ is stabilized by $G_{\Phi^\kappa}({\c M})$, by construction of 
$G_{\Phi^\kappa}({\c M})$. 
It follows that the image of $\Omega_{\Phi^\kappa}$ is precisely 
$G_{\Phi^\kappa}({\c N})$. 
\QED

%%%%%%%%%%%%%%%%%%%%%%%%%%%%%%%%%%%%%%%%%%%%%%%%%%%%%%%%%%%%%%%%%%%%%%%%%
%%%%%%%
\newsubchapter{Statement of the main theorem}
%%%%%%%
%%%%%%%%%%%%%%%%%%%%%%%%%%%%%%%%%%%%%%%%%%%%%%%%%%%%%%%%%%%%%%%%%%%%%%%%%

\newpar
{Main theorem }
{\sl The algebraic group $Gal({\c M},\eta)$ is the smallest closed subgroup of 
$GL(M)$ containing $\Phi_q^{\kappa_v}$ for almost all $v\in\Sg_f$.}
\label{galoiskatz}

\vpar\noindent
{\bf Example.}  
Let us consider the $q$-difference equation $y(qx)=q^{1/2}y(x)$,
associated to the $q$-difference module:
$$\matrix{
\Phi_q:& K(x)&\lrarrow & K(x)\cr
&f(x)&\longmapsto & q^{1/2}f(qx)}\ .
$$
The $q$-difference module $(K(x),\Phi_q)$ is trivial over
$K(x^{1/2})$, hence the generic Galois group of $(K(x),\Phi_q)$ 
is the group $\mu_2=\{1,-1\}$.
For all $v$ such that $|q|_v=1$ and such that the image of $q^{1/2}$ is an
element of the cyclic group generated by the image of $q$ in
$k_{q,v}$, the module $(K(x),\Phi_q)$ has $\kappa_v$-curvature zero.
For every other $v$ such that $|q|_v=1$, we have $\phi_q^{\kappa_v}\not\equiv 1$
and $\phi_q^{2\kappa_v}\equiv 1$ over $k_{q,v}$, which means that 
$\phi_q^{\kappa_v}\equiv -1$.  So the Galois
group is the smallest algebraic subgroup of the multiplicative
group $K(x)^\times\cong GL(K(x))$ containing $\Phi_q^{\kappa_v}$ for almost
all $v$.
\label{esempiotriv}

\vpar
The proof of the last statement relies on the $q$-analogue of
Grothendieck's conjecture on $p$-curvatures ({\it cf.}\ (\pgrothendieck)). 

\vpar
A part of the statement is very easy to prove:

\P
{The algebraic group $Gal({\c M},\eta)$ contains $\Phi_q^{\kappa_v}$ for almost all
$v\in\Sg_f$.}
\label{inclusione}

\Dim
The algebraic group $Gal({\c M},\eta)$ can be defined as the
stabilizer of a $q$-difference module $L$ of rank one over $K(x)$. 
The choice of an $\c F$-lattice $\widetilde M$ of $M$ 
determines an $\c F$-lattice $\widetilde L$ of $L$ of rank one. 
The reduction over $k_{q,v}$ of $\widetilde L$ is stable by the
morphism induced by $\Phi_q^{\kappa_v}$ since $\widetilde L$ is a $q$-difference
module, hence stable under $\Phi_q$.
\QED

%%%%%%%%%%%%%%%%%%%%%%%%%%%%%%%%%%%%%%%%%%%%%%%%%%%%%%%%%%%%%%%%%%%%%%%%%
%%%%%%%
\newsubchapter{Proof of the main theorem}
%%%%%%%
%%%%%%%%%%%%%%%%%%%%%%%%%%%%%%%%%%%%%%%%%%%%%%%%%%%%%%%%%%%%%%%%%%%%%%%%%

\vpar
Let ${\c M}=(M,\Phi_q)$ be a $q$-difference module over $K(x)$ and let 
$Gal({\c M},\eta)$ be its generic Galois group. 
We denote by $G_{\Phi^\kappa}({\c M})$ the smallest algebraic subgroup of $GL(M)$ containing 
$\Phi_q^{\kappa_v}$ for almost all $v$.
Our purpose is to prove that 
$Gal({\c M},\eta)=G_{\Phi^\kappa}({\c M})$ ({\it cf.}\ (\pgaloiskatz)). 
\par
We recall that we have already proved that 
$G_{\Phi^\kappa}\subset Gal({\c M},\eta)$ in (\pinclusione).

\vpar
We choose a $K(x)$-vector space $L$ of dimension 1 in a finite sum of the form 
$\oplus_hT^{i_h,j_h}(M)$, such that $G_{\Phi^\kappa}({\c M})$ is the stabilizer of $L$. 
\par
We denote by ${\c W}=(W,\Phi_q)$ the smallest $q$-difference sub-module of $\oplus_hT^{i_h,j_h}(M)$
containing $L$. 

\vpar
Let ${\c F}\subset K(x)$ be a $q$-difference algebra essentially of
finite type over ${\cal V}_K$ and $\widetilde M$ an $\c F$-lattice of $M$, stable by $\Phi_q$.
Let $m$ be a basis of the $\c F$-lattice $\widetilde L$ of $L$
determined by $\widetilde M$. Then $m$ is a cyclic vector for a
suitable $\c F$-lattice $\widetilde W$ of $W$. 
For almost all $v\in\Sg_f$, $\widetilde L\otimes_{{\c V}_K}k_{q,v}$ is stable with respect
to the morphism induced by $\Phi_q^{\kappa_v}$, which means that:
$$
\Phi_q^{\kappa_v}(m)\equiv \a_v(x)m\  
\hbox{in $\widetilde L\otimes_{{\c V}_K}k_{q,v}$, with $\a_v(x)\in {\c F}\otimes_{{\c V}_K}k_{q,v}$.}
$$
If $\nu$ is the rank of $W$, we obtain:
$$\eqalign{
&\Phi_q^{\kappa_v}(m,\Phi_q(m),\dots,\Phi_q^{\nu-1}(m))\cr
&\hskip 20 pt \equiv(m,\Phi_q(m),\dots,\Phi_q^{\nu-1}(m))
\pmatrix{\a_v(x)&&0\cr
         &\ddots&\cr
         0&&\a_v(q^{\nu-1}x)}\  
\hbox{in $\widetilde W\otimes_{{\c V}_K}k_{q,v}$.}}
$$
We deduce that the reduction modulo $\varpi_{q,v}$ 
of the sub-$\c F$-module of $\widetilde W$ generated by $\Phi_q^i(m)$,
for any $i=0,\dots,\nu-1$, is stable by $\Phi_q^{\kappa_v}$, for almost all $v$. 
This implies that the $K(x)$-vector space generated by $\Phi_q^i(m)$, 
for any $i=0,\dots,\nu-1$, is stable by $G_{\Phi^\kappa}({\c M})$.
Let us call $U$ the sub-$K(x)$-vector space of 
$W$ generated by $(\Phi_q(m),\dots,\Phi_q^{\nu-1}(m))$.  Then 
$W=L\oplus U$ is a decomposition of $W$ in subspaces stable by 
$G_{\Phi^\kappa}({\c M})$.
Let us consider the dual decomposition of $\V W$: $\V W=\V L\oplus\V U$. 
It follows that $G_{\phi^\kappa}({\c M})$ is the group fixing the line
$L\otimes\V L$ in $W\otimes\V W$ ({\it cf.}\ for instance the proof of theorem [D, 3.1]). 
\par
Let us consider the line $L\otimes\V L$ instead of the line $L$ to
define  $G_{\phi^\kappa}({\c M})$ as a stabilizer. 
Then we are in the following situation: $G_{\phi^\kappa}({\c M})$ is the group fixing the line
$L$ and ${\c W}=(W,\Phi_q)$ is the smallest $q$-difference module 
containing $L$.  
The $\c F$-lattice $\widetilde L$ of $L$ is a direct factor in a
suitable ${\c F}$-lattice of $\oplus_hT^{i_h,j_h}(M)$, hence 
$\widetilde L\otimes_{{\c V}_K}k_{q,v}$ is fixed by
$\Phi_q^{\kappa_v}$, for almost all $v\in\Sg_f$: 
$$
\Phi_q^{\kappa_v}(m)\equiv m\  
\hbox{in $\widetilde L\otimes_{{\c V}_K}k_{q,v}$, for all
  $m\in\widetilde L$.}
$$
Let us fix a cyclic vector $m\in\widetilde L$ for $\widetilde
W$. Then we have:
$$
\Phi_q^{\kappa_v}(m,\Phi_q(m),\dots,\Phi_q^{\nu-1}(m))
\equiv(m,\Phi_q(m),\dots,\Phi_q^{\nu-1}(m)){\b I}_\nu\  
\hbox{in $\widetilde W\otimes_{{\c V}_K}k_{q,v}$.}
$$
By (\pgrothendieck), the $q$-difference module 
$\c W$ is trivial and hence $Gal({\c W},\eta)=1$.     
Since ${\c W}\in<\!\!{\c M}\!\!>^\otimes$, we have a natural morphism
$$
Gal({\c M},\eta)\lrarrow Gal({\c W},\eta)=1\ ,
$$
which proves that $Gal({\c M},\eta)$ stabilizes each line of $W$. 
In particular $Gal({\c M},\eta)$ stabilizes $L$, hence 
$Gal({\c M},\eta)=G_{\phi^\kappa}({\c M})$. 
This completes the proof. 

%%%%%%%%%%%%%%%%%%%%%%%%%%%%%%%%%%%%%%%%%%%%%%%%%%%%%%%%%%%%%%%%%%%%%%%%%
%%%%%%%%%%%%%%%%%%%%%%%%%%%%  11  %%%%%%%%%%%%%%%%%%%%%%%%%%%%%%%%%%%%%%%%
%%%%%%%%%%%%%%%%%%%%%%%%%%%%%%%%%%%%%%%%%%%%%%%%%%%%%%%%%%%%%%%%%%%%%%%%%

\newchapter{Examples of calculation of generic Galois groups}
\chaplabel{esempi}

%%%%%%%%%%%%%%%%%%%%%%%%%%%%%%%%%%%%%%%%%%%%%%%%%%%%%%%%%%%%%%%%%%%%%%%%%

\vpar
We conclude with some examples of arithmetic calculations of the generic Galois
group. 

\vpar\noindent
{\bf Example 11.1.}
\par\noindent
Let us consider the $q$-difference equation 
$$
\pmatrix{y_1(qx)\cr y_2(qx)}=
\pmatrix {1 & a(x)\cr 0 &b(x)}
\pmatrix{y_1(x)\cr y_2(x)}\ ,
\leqno{(11.1.1)}$$
with $a(x)\neq 0$.  We suppose moreover that $a(x)$ does not have a zero
at $x=0$.  Then for all positive integers $n$ we obtain:
$$
\pmatrix{y_1(q^nx)\cr y_2(q^nx)}=
\pmatrix {1 & \ast\cr 0 &b(q^{n-1}x)\cdots b(x)}
\pmatrix{y_1(x)\cr y_2(x)}\ .
$$
Let us consider the $q$-difference modules ${\c M}$ generated by
$(e_1,e_2)$ with 
$$
\Phi_q(e_1,e_2)=(e_1,e_2)\pmatrix{1& 0\cr a(x)& b(x)}\ .
$$
Then the $q$-difference linear system (11.1.1) is associated to ${\c M}$
with respect to the basis $\u e$.  We distinguish several cases:
\par\noindent
1) $\u{\hbox{\it $y(qx)=b(x)y(x)$ has a solution in $K(x)$}}$ \hfill\break
Then the generic Galois group of ${\c M}$ is:
$$
Gal({\c M},\eta)=\l\{\pmatrix{1&0\cr c(x)&1}~:~c(x)\in K(x)\r\}\ .
$$
\par\noindent
2) $\u{\hbox{\it $y(qx)=b(x)y(x)$ has a solution in 
an extension $K(q^{1/d})(x^{1/d})$ of $K(x)$}}$, \hfill\break
for a suitable integer $d>1$.  We choose $d$ minimal with respect to this
property.  We obtain:
$$
Gal({\c M},\eta)=\l\{\pmatrix{1&0\cr c(x)&\zeta}~:~c(x)\in K(x),\zeta\in\mu_d\r\}\ .
$$
\par\noindent
3) $\u{\hbox{\it none of the previous conditions is satisfied}}$.\hfill\break 
We find the algebraic group:
$$
Gal({\c M},\eta)=\l\{\pmatrix{1&0\cr c(x)&d(x)}~:~c(x),d(x)\in K(x), d(x)\neq 0\r\}\ .
$$

\vpar\noindent
{\bf Example 11.2.}
\par\noindent
Let us consider the $q$-difference linear system of order two:
$$
\pmatrix{y_1(qx)\cr y_2(qx)}=
\pmatrix {0 & r(x)\cr 1 & 0}
\pmatrix{y_1(x)\cr y_2(x)}\ , 
$$
with $r(x)\in K(x)$ and $r(x)\neq 0$.  We can easily calculate by induction
that for all positive integers $n$ we have:
$$
\pmatrix{y_1(q^{2n}x)\cr y_2(q^{2n}x)}=
\pmatrix { r(q^{2n-1}x)\cdots r(q^3x)r(qx)&0\cr 0 & r(q^{2n-2}x)\cdots r(q^2x)r(x)}
\pmatrix{y_1(x)\cr y_2(x)}
$$
and
$$\eqalign{\pmatrix{y_1(q^{2n+1}x)\cr y_2(q^{2n+1}x)}
%&=\pmatrix { r(q^{2n}x)\cdots r(q^4x)r(q^2x)&0\cr 0 & r(q^{2n-1}x)\cdots r(q^3x)r(qx)}
%  \pmatrix {0 & r(x)\cr 1 & 0}\pmatrix{y_1(x)\cr y_2(x)}\cr
&=\pmatrix {0& r(q^{2n}x)\cdots r(q^2x)r(x) \cr r(q^{2n-1}x)\cdots r(q^3x)r(qx) &0}
  \pmatrix{y_1(x)\cr y_2(x)}\ .}
$$
It follows that for the generic Galois group of the $q$-difference module 
${\c M}$ of rank 2 such that for a fixed basis $(e_1,e_2)$ we have:
$$
\Phi_q(e_1,e_2)=(e_1,e_2)\pmatrix{0& 1\cr r(x)&0}
$$
there are two possibilities:
\par\noindent
1) $\u{\hbox{\it $y(q^2x)=r(x)y(x)$ has a solution in $K(x)$}}$\hfill\break
Then $\Phi_q^{2\kappa_v}\equiv 1$ modulo $\varpi_{q,v}$ for almost all $v$ 
and the generic Galois group of ${\c M}$ is represented as the algebraic 
linear subgroup $Gl_2(K(x))$ of the form:
$$
Gal({\c M},\eta)=\l\{{\b I}_2,-{\b I}_2,\pmatrix{0&1\cr 1&0},
\pmatrix{0&-1\cr -1&0}\r\}\ .
$$
\par\noindent
2) $\u{\hbox{\it $y(q^2x)=r(x)y(x)$ has a solution in an extension 
$K(q^{2/d})(x^{1/d})$ of $K(x)$}}$,\hfill\break 
for a suitable integer $d>1$. We choose $d$ minimal with respect to this
property.  Then:
$$
Gal({\c M},\eta)=\l\{\pmatrix{\zeta_1&0\cr
  0&\zeta_2}~:~\zeta_1,\zeta_2\in\mu_d\r\}
\cup\l\{\pmatrix{0&\zeta_3\cr
  \zeta_4&0}~:~\zeta_3,\zeta_4\in\mu_d\r\}\ .
$$
\par\noindent
3) $\u{\hbox{\it none of the previous conditions is satisfied}}$.\hfill\break 
then the generic Galois group of ${\c M}$ is represented as the infinite algebraic 
linear subgroup $Gl_2(K(x))$ of the form:
$$
\eqalign{Gal({\c M},\eta)
&=\l\{\pmatrix{a(x)&0\cr 0&b(x)}~:~a(x),b(x)\in K(x), a(x)b(x)\neq 0\r\}\cr
&\hskip 50 pt\cup\l\{\pmatrix{0&c(x)\cr d(x)&0}~:~c(x),d(x)\in K(x), c(x)d(x)\neq 0\r\}\ .}
$$

%%%%%%%%%%%%%%%%%%%%%%%%%%%%%%%%%%%%%%%%%%%%%%%%%%%%%%%%%%%%%%%%%%%%%%%%%%%%
%%%%%%%%%%%%%%%%%%%%%%%%%%%%%%%%%%%%%%%%%%%%%%%%%%%%%%%%%%%%%%%%%%%%%%%%%%%%
%%%%%%%%%%%%%%%%%%%%%%%%%%%%%%%%%%%%%%%%%%%%%%%%%%%%%%%%%%%%%%%%%%%%%%%%%%%%
%%%%%%%%%%%%%%%%%%%%%%%%%%%%%%%%%%%%%%%%%%%%%%%%%%%%%%%%%%%%%%%%%%%%%%%%%%%%
%%%%%%%%%%%%%%%%%%%%%%%%%%%%%%%%%% APPENDIX  %%%%%%%%%%%%%%%%%%%%%%%%%%%%%%%
%%%%%%%%%%%%%%%%%%%%%%%%%%%%%%%%%%%%%%%%%%%%%%%%%%%%%%%%%%%%%%%%%%%%%%%%%%%%
%%%%%%%%%%%%%%%%%%%%%%%%%%%%%%%%%%%%%%%%%%%%%%%%%%%%%%%%%%%%%%%%%%%%%%%%%%%%
%%%%%%%%%%%%%%%%%%%%%%%%%%%%%%%%%%%%%%%%%%%%%%%%%%%%%%%%%%%%%%%%%%%%%%%%%%%%
%%%%%%%%%%%%%%%%%%%%%%%%%%%%%%%%%%%%%%%%%%%%%%%%%%%%%%%%%%%%%%%%%%%%%%%%%%%%

\vskip 30pt\noindent
\centerline{\titlepart Appendix.  A {\titlepartit q}-analogue 
of Schwarz's list}

\def\parag{A}\chapno=0\subchapno=0\equano=0

\vpar\vpar
{\it Let $a,b,c,q$ be complex numbers. We suppose that 
$q$ is not zero and not a root of unity.} 

\vpar
We consider the basic hypergeometric function:
$$
{}_2\phi_1(a,b,c;q,x)=\sum_{n\geq 0}
{(a;q)_n(b;q)_n\over(c;q)_n(q;q)_n}x^n\ ,
$$
where $(a;q)_n=(1-a)(1-aq)\cdots (1-aq^{n-1})$. 
It is defined if $c\not\in q^{\interi_{\leq 0}}$ 
or if $c\in q^{\interi_{\leq 0}}$ and either 
$a\in q^{\interi_{\leq 0}}$, $ac^{-1}\in q^{\interi_{\geq 0}}$ or 
$b\in q^{\interi_{\leq 0}}$, $bc^{-1}\in q^{\interi_{\geq 0}}$. 
\par
It is a $q$-analogue of the Gauss hypergeometric series
$$
{}_2F_1(\a,\be,\ga;x)=\sum_{n\geq 0}
{(a)_n(b)_n\over(c)_n(q)_n}x^n\ ,
$$
where $(\a)_n=\a(\a+1)\cdots(\a+n-1)$ is the Pochhammer symbol. 
If $c$ is a nonpositive integer ${}_2F_1(\a,\be,\ga;x)$ is defined if and only if either 
$a\in\interi$, $c\leq a\leq 0$ or $b\in\interi$, $c\leq b\leq 0$.  
\par 
The series ${}_2\phi_1(a,b,c;q,x)$ is a 
solution of the basic hypergeometric $q$-difference equation
$$
\varphi_q^2y(x)-{(a+b)x-(1+cq^{-1})\over abx-cq^{-1}}
\varphi_q y(x)+
{x-1\over abx-cq^{-1}}y(x)=0\ , 
\leqno{({\cal H}_{a,b,c})}$$
which is defined as soon as 
neither $a=c=0$ nor $b=c=0$. 

\vpar
Our purpose is to make a ``list'' of all the parameters 
$(a,b,c)$ such that (${\cal H}_{a,b,c}$) has a basis of algebraic solutions 
({\it i.e.} in a finite extension of $\complessi(x)$).  
It may be thought of as an analogue of Schwarz's list 
for hypergeometric differential equations having solution ${}_2F_1(\a,\be,\ga;x)$:
$$
y^{\p\p}(x)+{\ga-(\a+\be+1)x\over x(1-x)}y^\p(x)-
{\a\be\over  x(1-x)}y(x)=0\ ,
\leqno{({\cal E}_{\a,\be,\ga})}$$
where $\a$, $\be$, $\ga$ are complex parameters. 
We also establish the list of parameters $a,b,c$ such that 
(${\cal H}_{a,b,c}$) has a basis of solutions in $\complessi(x)$. 
The complete solution to this problem is actually 
closely linked with the solution of the analogous problem 
for (${\cal E}_{a,b,c}$). 
\par 
We obtain:

\vpar
\centerline{\boxit{\hsize 13 true cm\noindent
$$\eqalign{
\hbox{(${\cal H}_{a,b,c}$) has a basis}&\hbox{ of solutions in $\complessi(x)$}\cr
\Leftrightarrow &
   \cases{\hbox{1) there exists $\a,\be,\ga\in\interi$, such that 
                $a=q^\a$, $b=q^\be$, $c=q^\ga$;}\cr
          \hbox{2) (${\cal E}_{\a,\be,\ga}$) has a basis of 
                solutions in $\complessi(x)$.}}\cr
\hbox{(${\cal H}_{a,b,c}$) has a basis}&\hbox{ of algebraic solutions}\cr
\Leftrightarrow 
     &\hbox{ either it has a basis of solutions in $\complessi(x)$}\cr
     &\hbox{ or the following conditions are satified:}\cr
     &\hskip 5pt\matrix{\hbox{1) $a,b,c\in q^\raz$;}\hfill\cr
              \hbox{2) either $a, bc^{-1}\in q^\interi$ or 
                    $b, ac^{-1}\in q^\interi$;}\hfill\cr
              \hbox{3) $ab^{-1},c\not\in q^\interi$.}\hfill}\cr}
$$}}

\vpar
We can state a more precise result. Let
$$
{\cal Z}=\l(\interi_{>0}\times\interi_{\leq 0}\r)
\cup\l(\interi_{\leq 0}\times\interi_{>0}\r)\ .
$$
We know by [G, Ch. III] that: 

\P
{Let $\a,\be,\ga\in\interi$. The following assertions 
are equivalent:
\tparnoind 
1) (${\cal E}_{\a,\be,\ga}$) has a basis of solutions in $\complessi(x)$; 
\tparnoind
2) $|1-\ga|$, $|\ga-\a-\be|$ and $|\a-\be|$ are the 
lengths of the sides of a triangle; 
\tparnoind
3) the following conditions are satisfied: 
\itemitem{$\ast$} either $(\a,\a+1-\ga)\in{\cal Z}$ 
or $(\be,\be+1-\ga)\in{\cal Z}$,
\itemitem{$\ast$} either $(\a,\be)\in{\cal Z}$ or 
$(\a+1-\ga,\be+1-\ga)\in{\cal Z}$.}

\vpar
Condition {\sl 2)} must be interpreted in the most elementary 
way: if $n_1,n_2,n_3$ are positive integers such that $n_1\leq n_2\leq n_3$, 
they can be the lengths of a triangle if they satify the condition 
$n_1+n_2\geq n_3$.

\vpar
We are thus led to prove the following proposition:

\P
{The $q$-difference equation (${\cal H}_{a,b,c}$) 
has a basis of solutions in $\complessi(x)$ if and only if 
the following conditions are satisfied:
\itemitem{$\ast$} there exists $\a,\be,\ga\in\interi$, 
                  such that  $a=q^\a$, $b=q^\be$, $c=q^\ga$;
\itemitem{$\ast$} either $(\a,\a+1-\ga)\in{\cal Z}$ or 
$(\be,\be+1-\ga)\in{\cal Z}$,
\itemitem{$\ast$} either $(\a,\be)\in{\cal Z}$ or 
$(\a+1-\ga,\be+1-\ga)\in{\cal Z}$.}
\label{razionali}

\Oss
A basis of solutions of (${\cal H}_{a,b,c}$) at zero 
is given generically by 
$$
\cases{{}_2\phi_1(a,b,c;q,x)\cr 
e_{qc^{-1}}(x){}_2\phi_1(aqc^{-1},bqc^{-1},q^2c^{-1};q,x)}\ ,
$$
where $e_{qc^{-1}}(x)$ is a solution of $y(qx)=(qc^{-1})y(x)$. 
A  basis of solutions at $\infty$ is given by
$$
\cases{\displaystyle e_a\l({1\over x}\r){}_2\phi_1\l(a,aqc^{-1},aqb^{-1};q,{cq\over abx}\r)\cr\cr 
\displaystyle e_b\l({1\over x}\r){}_2\phi_1\l(b,bqc^{-1},bqa^{-1};q,{cq\over abx}\r)}\ .
$$
We notice that they are not always well defined when 
$a,b,c\in q^\interi$. 
\par
Let $\a,\be,\ga\in\interi$ be such that  $a=q^\a$, $b=q^\be$, 
$c=q^\ga$. 
Then $e_{qc^{-1}}(x)=x^{1-\ga}$, $e_a(1/x)=(1/x)^a$ and 
$e_b(1/x)=(1/x)^\be$. 
The previous proposition says that 
(${\cal H}_{a,b,c}$) has a basis of rational solutions if and only if 
one of the two basis of solutions above is well 
defined: in this case one solution is necessarily a polynomial, the other 
one is a rational function which can be written explicitly using Heine's formula 
[GR, 1.4.6]: 
$$
{}_2\phi_1(a,b,c;q,x)=
{((abc^{-1})x;q)_\infty\over(x;q)_\infty}
{}_2\phi_1(ca^{-1},cb^{-1},c;q,abc^{-1}x)\ , 
$$
or Heine's contiguity relations [GR, Ex. 1.10]. 

\vpar
In the next section we will give a proof of (\prazionali).  We will discuss the 
case of algebraic solutions later.

%%%%%%%%%%%%%%%%%%%%%%%%%%%%%%%%%%%%%%%%%%%%%%%%%%%%%%%%%%%%

\newsubchapter{Logarithmic singularities}

%%%%%%%%%%%%%%%%%%%%%%%%%%%%%%%%%%%%%%%%%%%%%%%%%%%%%%%%%%%%

\vpar
Let us consider a $q$-difference equation of order two
$$
y(q^2x)+P(x)y(qx)+Q(x)y(x)=0\ ,
\numero{equazione}$$
with $P(x), Q(x)\in \complessi(x)$. 
One says that it is {\it regular singular at zero} if $P(x)$ 
has no pole at zero and if $Q(x)$ has neither 
a pole nor a zero at zero. 

\vpar
Let $e_c(x)$ be a solution of the $q$-difference 
equation $y(qx)=cy(x)$, with 
$c\in\complessi$ ({\it cf.}\ [S2, 0.1]). 
Suppose that (\fequazione) has a solution of the 
form $e_c(x)\sum_{n\geq 0}a_nx^n$. 
Then $c$ satisfies the equation
$$
c^2+P(0)c+Q(0)=0\ . 
$$
The two roots $c_1$ and $c_2$ of this 
equation are called the {\it exponents}
of (\fequazione) at zero.
\par 
When $c_1$, $c_2$ satisfy 
the condition $c_1\not\in c_2q^\interi$, 
we know by [S2, 1.1.4, Th.] that (\fequazione) 
has a basis of solutions of the form 
$e_{c_1}(x)\sum_{n\geq 0}a_nx^n$, 
$e_{c_2}(x)\sum_{n\geq 0}b_nx^n$. 
Moreover, $\sum_{n\geq 0}a_nx^n$ and 
$\sum_{n\geq 0}b_nx^n$ are convergent complex power series; 
if $|q|_\complessi\neq 1$ they are Taylor expansions of 
meromorphic functions on $\complessi$.
If the condition $c_1\not\in c_2q^\interi$ is not satisfied 
this is in general not true: the two solutions may involve 
$q$-logarithms, which are solutions of the 
$q$-difference equation
$y(qx)=y(x)+1$ ({\it cf.}\ [S2, 0.1]). 

\Def
{We say that a $q$-difference equation 
which is regular singular at zero 
does {\it not} have a {\it logarithmic singularity at zero} 
if it has a basis of solutions of the form 
$e_{c_1}(x)\sum_{n\geq 0}a_nx^n$, 
$e_{c_2}(x)\sum_{n\geq 0}b_nx^n$. } 

\vpar
One can give an analogous definition for the point $\infty$, by 
using the variable change $t={1\over x}$. 
We have:

\L
{The equation 
(${\cal H}_{a,b,c}$) has a basis of solutions in $\complessi(x)$ 
if and only if there exists  $\a,\be,\ga\in\interi$ 
such that  $a=q^\a$, $b=q^\be$, $c=q^\ga$, and 
zero and $\infty$ are not logarithmic singularities.}
\label{log}

\Dim
We notice that the exponents of (${\cal H}_{a,b,c}$) at zero 
(resp. $\infty$) are $1$ and $qc^{-1}$ (resp. $a$ and $b$), 
and that the equation 
(${\cal H}_{a,b,c}$) is actually defined over 
the field $\raz(a,b,c,q)(x)$. 
\par
Let us suppose that there exists  $\a,\be,\ga\in\interi$ 
such that  $a=q^\a$, $b=q^\be$, $c=q^\ga$, and 
zero and $\infty$ are not logarithmic singularities.
Then (${\cal H}_{a,b,c}$) has a basis of solutions 
at zero of the form
$$
\hbox{$u_0(x),x^{1-\gamma}v_0(x)$, with 
$u_0(x),v_0(x)\in\raz(q)\[[x\]]$,}
$$ 
and a basis of solutions at $\infty$ of the form
$$
\hbox{${1\over x^\a} u_\infty\l({1\over x}\r),
{1\over x^\be}v_\infty\l({1\over x}\r)$, with 
$u_\infty\l({1\over x}\r),
v_\infty\l({1\over x}\r)\in\raz(q)\[[{1\over x}\]]$.}
$$ 
Let $|~|_q$ be the $q$-adic norm over $\raz(q)$. Then 
$u_0(x),v_0(x)$ (resp. $u_\infty(1/x),v_\infty(1/x)$) 
have infinite radius of meromorphy 
at zero (resp. $\infty$) with respect to the norm $|~|_q$, 
since $|q|_q<1$. 
\par
Let us consider the Birkhoff connection matrix 
({\it cf.}\ [S2, \S2])
$$
P(x)=\pmatrix{(1/x)^\a u_\infty(1/x)&(1/x)^\be v_\infty(1/x)\cr
              (1/qx)^\a u_\infty(1/qx)&(1/qx)^\be 
               v_\infty(1/qx)}^{-1}
     \pmatrix{u_0(x) &x^{1-\gamma}v_0(x)\cr
              u_0(qx)&(qx)^{1-\gamma}v_0(qx)}\ .
$$ 
The entries of $P(x)$ are $q$-adically meromorphic elliptic functions 
over ${\b P}^1_{\raz(q)}\setminus\{0,\infty\}$.
Moreover, since zero and $\infty$ are not 
logarithmic singularities and the exponents are 
integral powers of $q$, the 
singularities of $P(x)$ are in the $q$-orbits of 
the singularities 
of the coefficients of the equation (${\cal H}_{a,b,c}$). 
It follows that $P(x)$ induces a meromorphic elliptic 
function over $\raz(q)^\times/q^\interi$, 
having at worst one pole 
at $q^\interi$, hence the matrix $P(x)$ must be constant. 
Therefore $u_0(x)$, $v_0(x)$ are meromorphic functions over 
${\b P}^1_{\raz(q)}$, {\it i.e.}, they are rational functions. 
\par
The other implication of the equivalence is clear. 
\QED

\vpar
The following two lemmas achieve the proof of 
(\prazionali). Their proof is inspired by the
proof of the analogous statements in [G, Ch. III]
for the hypergeometric differential equation. 

\vpar
Let 
$$
q^{\cal Z}=\l(q^{\interi_{>0}}\times q^{\interi_{\leq 0}}\r)
\cup\l(q^{\interi_{\leq 0}}\times q^{\interi_{>0}}\r)\ .
$$

\L
{Let $c=q^\ga$ with $\ga\in\interi$. 
Then (${\cal H}_{a,b,c}$) has a logarithmic 
singularity at zero if and only if 
either $(a,qac^{-1})\in q^{\cal Z}$ or 
$(b,qbc^{-1})\in q^{\cal Z}$.}
\label{lemmauno}

\Dim
Let $\ga=1$. Then we know by [S2, 1.1.4, Th.] that 
zero is a logarithmic singularity. 
\par
Let us suppose that $\ga\leq 0$. Under this 
assumption, the series 
$$
x^{1-\ga}{}_2\phi_1(aqc^{-1},bqc^{-1},q^2c^{-1};q,x)
$$
is a well defined solution of (${\cal H}_{a,b,c}$). 
Zero is not a logarithmic singularity for (${\cal H}_{a,b,c}$) 
if and only if there exists a second solution 
of the form $\sum_{n\geq 0}c_nx^n$, with $c_0=1$. 
A direct calculation shows that 
the coefficients $c_i$ must satisfy the relation: 
$$
(1-aq^n)(1-bq^n)c_n=(1-q^{\ga+n})(1-q^{n+1})c_{n+1}\ .
\numero{ricorrenza}$$
Let us suppose that neither $a=q^\a$, with $0\geq \a\geq \ga$, 
nor $b=q^\be$, with $0\geq \be\geq \ga$. 
Then $c_1,\dots, c_{-\ga}$ can be determined inductively: they 
are necessarily nonzero complex numbers. 
For $n=-\ga$, we get 0 on the right hand 
side of (\fricorrenza), while 
the left hand side is not zero. Therefore 
we cannot find a solution of the form $\sum_{n\geq 0}c_nx^n$, with $c_0=1$, 
and zero is a logarithmic singularity.
\par
On the other hand, if either $a=q^\a$, 
with $0\geq \a\geq \ga$, 
or $b=q^\be$, with $0\geq \be\geq \ga$, the series  
${}_2\phi_1(a,b,c;q,x)$ is a well defined 
solution of (${\cal H}_{a,b,c}$), satisfying 
${}_2\phi_1(a,b,c;q,0)=1$.
\par
If $\ga\leq 0$, we conclude that zero is not a 
logarithmic singularity if and only if either 
$a=q^\a$, with $0\geq \a\geq \ga$, 
or $b=q^\be$, with $0\geq \a\geq \ga$. 
This completes the proof in the case $\ga\leq 0$. 
\par
Let $\ga\geq 2$. A series of the form 
$y(x)=x^{1-\ga}z(x)$ is a solution of (${\cal H}_{a,b,c}$)
if and only if $z(x)$ is a solution of 
(${\cal H}_{a^\p,b^\p,c^\p}$), with 
$a^\p=aq^{1-\ga}$, $b^\p=bq^{1-\ga}$ 
and $c^\p=q^{2-\ga}$. So we can deduce 
this case from the case $\ga\leq 0$. 
\QED

\L
{Let $ab^{-1}\in q^\interi$. Then  
(${\cal H}_{a,b,c}$) has a logarithmic singularity 
at $\infty$ if and only if 
either $(a,b)\in q^{\cal Z}$ or 
$(qac^{-1},qbc^{-1})\in q^{\cal Z}$.}
\label{lemmadue}

\Dim
Let $t=cq/abx$. Then $y(t)=e_a(abc^{-1}q^{-1}t)z(t)$ is a solution of 
(${\cal H}_{a,b,c}$) if and only if $z(t)$ is a solution of 
(${\cal H}_{a^\p,b^\p,c^\p}$), with $a^\p=a$, $b^\p=qac^{-1}$ 
and $c^\p=qab^{-1}\in q^\interi$: this fact 
can be deduced by the remark that, for generic parameters $a,b,c$,  
a solution of (${\cal H}_{a,b,c}$) at $\infty$ is given by 
$e_a\l(abc^{-1}q^{-1}t\r){}_2\phi_1\l(a,aqc^{-1},aqb^{-1};q,t\r)$. 
We conclude by (\plemmauno).
\QED

%%%%%%%%%%%%%%%%%%%%%%%%%%%%%%%%%%%%%%%%%%%%%%%%%%%%

\newsubchapter{The case of algebraic solutions}

%%%%%%%%%%%%%%%%%%%%%%%%%%%%%%%%%%%%%%%%%%%%%%%%%%%%

\P
{The $q$-difference equation (${\cal H}_{a,b,c}$) 
has a basis of algebraic solutions if and only if 
one of the following conditions is satisfied:
\tparnoind
1) (${\cal H}_{a,b,c}$) has a basis of rational solutions; 
\tparnoind
2) $a,b,c\in q^\raz$, $c, ab^{-1}\not\in q^\interi$ 
and either $a, bc^{-1}\in q^\interi$ or 
$b, ac^{-1}\in q^\interi$.}
\label{algebrici}

\Oss
If the condition {\sl 2)} is satisfied, 
the equation (${\cal H}_{a,b,c}$) has a 
well-defined basis of solutions of
the form 
$$
\l({1\over x}\r)^\a {}_2\phi_1\l(a,aqc^{-1},aqb^{-1};q,{cq\over abx}\r),
\ \l({1\over x}\r)^\be {}_2\phi_1\l(b,bqc^{-1},bqa^{-1};q,{cq\over abx}\r)\ ,
\numero{basealgebrica}$$ 
where 
${}_2\phi_1(a,aqc^{-1},aqb^{-1};q,cq/abx)$ and 
${}_2\phi_1(b,bqc^{-1},bqa^{-1};q,cq/abx)$ 
are rational functions and 
$\a,\be\in\raz$ are such that $a=q^\a$ and $b=q^\be$.

\vpar
The remark above proves one side of the equivalence 
in (\palgebrici). 
On the other side we have:

\L
{The $q$-difference equation (${\cal H}_{a,b,c}$) 
has only algebraic solutions if and only if
the following conditions are satisfied:
\item{$\ast$} $a,b,c\in q^\raz$;
\item{$\ast$} (${\cal H}_{a,b,c}$) has no logarithmic 
singularities at zero or at $\infty$.
\tparnoind
The solutions are then of the form $x^\delta u(x)$, 
where $\delta\in\raz$ and 
$u(x)\in\complessi(x)$.}

\vpar
The proof of this lemma is similar to the proof of lemma 
\plog\ and it is based on the remark that 
the equation $y(qx)=cy(x)$ has no algebraic 
solutions if and only if $c\in\complessi^\times\setminus q^\raz$ 
(the Galois group being infinite). 
 
By lemmas (\plemmauno) and (\plemmadue) we deduce that 
the $q$-difference equation (${\cal H}_{a,b,c}$) 
has only algebraic solutions if and only 
if the following conditions are satisfied:
\item{$\ast$} $a,b,c\in q^\raz$;
\item{$\ast$} either $ab^{-1}\not\in q^\interi$ or 
    $(a,b)\in q^{\cal Z}$ or $(aq/b,bq/a)\in q^{\cal Z}$;
\item{$\ast$} either $c\not\in q^\interi$ or 
    $(a,qac^{-1})\in q^{\cal Z}$ or $(b,qbc^{-1})\in q^{\cal Z}$.
\item{$\ast$} either $a, bc^{-1}\in q^\interi$ or 
    $b, ac^{-1}\in q^\interi$.
\par\noindent
The last condition follows from the fact that the 
solutions are of the form $x^\delta u(x)$, where $\delta\in\raz$ and 
$u(x)\in\complessi(x)$. 
\par
If $ab^{-1}\not\in q^\interi$ and $c\not\in q^\interi$ 
then the solutions are in a Kummer extension 
$\complessi(x^{1/n})$ of $\complessi(x)$, otherwise thay are rational functions.

%%%%%%%%%%%%%%%%%%%%%%%%%%%%%%%
%%%%%%%%%%%%%%%%%%%%%%%%%%%%%%% BIBLIOGRAFIA
%%%%%%%%%%%%%%%%%%%%%%%%%%%%%%%

\vskip 30pt\noindent
\centerline{\titlepart References.}
\vpar

\def\libri#1#2#3#4#5
{\item {[{#1}]}{#2}: ``{\it #3}", {#4}.\medskip}
\def\articoli#1#2#3#4#5
{\item {[{#1}]}{#2}: ``{#3}", {#4}.\medskip}

\libri
{A1}
{Andr{\'e} Y.}
{G-functions and Geometry}
{Aspects of Mathematics E13, Vie\-weg, Braunschweig/Wies\-baden, 1989}

\articoli
{A2}
{Andr{\'e} Y.}
{Sur la conjecture des $p$-courbures de Grothendieck-Katz}
{Institut de Math{\'e}\-matiques de Jussieu, preprint 134, octobre 1997}

\articoli{A3}
{Andr{\'e} Y.}
{S{\'e}ries Gevrey de type arithm{\'e}tique.  II. Transcendance sans
transcendance}  
{Annals of Math, 151, 2 (2000), 741-756}

\articoli{A4}
{Andr{\'e} Y.}
{Diff{\'e}rentielles non-commutatives et th{\'e}orie de Galois
diff{\'e}rentielle ou aux diff{\'e}rences}  
{to appear in Annals Scient. {\'E}c. Norm. Sup.}

\libri{AB}
{Andr{\'e} Y., Baldassarri F.}
{De Rham Cohomology of Differential Modules on Algebraic Varieties}
{Progress in Mathematics, 189. Birkhäuser Verlag, Basel, 2001}

\articoli{B}
{Baker A.}
{The theory of linear forms in logarithms} 
{Transcendence theory: advances and applications (Proc. Conf., Univ. Cambridge, Cambridge, 1976), pp. 1--27.
Academic Press, London, 1977}

\articoli
{Be1}
{B{\'e}zivin J.P.}
{Les suites $q$-r{\'e}currentes lin{\'e}aires} 
{Compositio Mathematica 80, 285-307, (1991)}

\articoli
{Be2}
{B{\'e}zivin J.P.}
{Sur les {\'e}quations fonctionnelles aux $q$-diff{\'e}rences}
{Aequationes math. 43, (1992)}

\articoli
{BB}
{B{\'e}zivin J.P., Boutabaa A.}
{Sur les {\'e}quations fonctionnelles $p$-adique aux $q$-diff{\'e}rences}
{Collect. Math. 43, 2 (1992)}

\articoli{BHV}
{Bilu Yu., Hanrot G., Voutier P. M.} 
{Existence of primitive divisors of Lucas and Lehmer numbers. With an appendix by M. Mignotte.} 
{J. Reine Angew. Math. 539 (2001), 75-122}

\articoli{Bo}
{Bost J.B.}
{Algebraic leaves of algebraic foliations over number fields} 
{Publ. Math. IHES 93 (2001), 161-221} 

\articoli{CC}
{Chudnovsky D.V., Chudnovsky G.V.}
{Applications of Pad{\'e} approximations to the Gro\-then\-dieck conjecture
on linear differential equations}
{Lect. Notes Math. 1135 (1985), 52-100}

\libri{Ch}
{Cohen P.}
{Skew field constructions}
{Cambridge Univ. Press, 1977}

\articoli{D}
{Deligne P.}
{Hodge cycles on abelian variety}
{Springer Lecture Notes 900 (1982), 9-100}

\articoli{DM}
{Deligne P., Milne J.}
{Tannakian categories}
{Springer Lecture Notes 900 (1982), 101-228}

\articoli{Du}
{Duval A.}
{Lemmes de Hensel et Factorisation Formelle pour les Op{\'e}rateurs aux Dif\-f{\'e}ren\-ces}
{Funkcialaj Ekvacioj, 26 (1983), 349-368}

\libri{DGS}
{Dwork B., Gerotto G., Sullivan F.}
{An Introduction to {\rm G}-functions}
{Annals of Mathematical Studies 133, Princeton University Press, 
Princeton N.J., 1994}

\articoli{E}
{Etingof P. I.}
{Galois groups and connection matrices of $q$-difference equations}
{Electronic Research Announcements of the A.M.S., Vol. I (1),1995}

\libri{GR}
{Gasper G., Rahman M.}
{Basic hypergeometric series}
{Cambridge Univ. Press, Cambridge, 1990}

\libri{G}
{Goursat E.}
{Le\ce{c}ons sur les s\'eries hyperg\'eom\'etriqies et sur quelques fonctions 
qui s'y rattachent}
{Hermann, Paris, 1936}

\articoli{H}
{Hendriks P. A.}
{Algebraic Aspects of Linear Differential and Difference Equations}
{Ph. D. thesis, University of Groningen, 1996}

\articoli{Ho}
{Honda T.}
{Algebraic differential equations}
{I.N.D.A.M. Symp. Math. XXIV (1981), 169-204}

\articoli{K1}
{Katz N. M.}
{Nilpotent connections and the monodromy theorem. Applications of a result of Turrittin} 
{Publ. Math. IHES 39 (1970), 175-232}

\articoli{K2}
{Katz N. M.}
{Algebraic solution of differential equations ($p$-curvature and Hodge Filtration)}
{Invent. Math. 18, (1972), 1-118}

\articoli{K3}
{Katz N. M.}
{A conjecture in the arithmetic theory of differential equations}
{Bull. S.M.F. 110, (1982), 203-239, and Bull. S.M.F. 111, (1982), 347-348}

\articoli{K4}
{Katz N. M.}
{On the calculation of some differential Galois groups}
{Invent. math. 87, 13-61 (1987)}

\libri{N}
{Neukirch J.}
{Algebraic number theory}
{Grund\-leh\-ren der math. Wissenschaften 322, Springer-Verlag, 1999}

\articoli{P}
{Praagman C.}
{The formal classification of linear difference equation}
{Proc. Kon. Ned. Ac. Wet. Ser. A, 86, 1983}

\libri{PS}
{van der Put M., Singer M. F.}
{Galois Theory of Difference Equations}
{Lecture Notes in Mathematics 1666, Springer, 1997}

\articoli{Ra}
{Ramis J-P.}
{About the growth of entire functions solutions of linear
  $q$-difference equations}
{Ann. de la Fac. des Sci. de Toulouse, S\'erie 6, Vol. I, 1 (1992), 53-94}

\libri{S1}
{Sauloy J.}
{Th{\'e}orie de Galois des {\'e}quations aux $q$-diff{\'e}rences
fuchsiennes}
{Th{\`e}se doctorale de l'Universit{\'e} Paul Sabatier de Toulouse}

\libri{S2}
{Sauloy J.}
{Syst\`eme aux $q$-diff\'erences singuliers r\'eguliers~:
classification, matrice de connexion et monodromie}
{Ann. Inst. Fourier (Grenoble), 50 (2000), no. 4, 1021-1071}

\articoli{Sc}
{Schinzel A.}
{Abelians binomials, power residues and exponential congruences}
{Acta Arithmetica, 27, 1975, 397-420}

\libri{Se}
{Serre J.P.}
{Lectures on the Mordell-Weil theorem}
{Aspects of Mathematics E15, Vie\-weg, Braunschweig/Wiesbaden, 1989}

\libri{W}
{Waterhouse W.C.}
{Introduction to Affine Group Schemes}
{ Springer, Graduate Texts in Mathematics 66, 1979}

\end